\documentclass{article}
\usepackage[margin=1.2in]{geometry}
\usepackage{setspace}
\RequirePackage{amsthm,amsmath,amsfonts,amssymb}
\RequirePackage[authoryear]{natbib}
\RequirePackage[colorlinks,citecolor=blue,urlcolor=blue]{hyperref}
\theoremstyle{plain}
\newtheorem{theorem}{Theorem}[section]
\newtheorem{proposition}[theorem]{Proposition}
\newtheorem{lemma}[theorem]{Lemma}
\newtheorem{corollary}[theorem]{Corollary}
\theoremstyle{definition}

\newtheorem{assumption}{Assumption}
\newtheorem{example}{Example}
\newtheorem{remark}{Remark}

\usepackage{amssymb}

\usepackage{multirow}
\usepackage{colortbl}
\usepackage{threeparttable}

\usepackage{microtype}
\usepackage{graphicx}
\usepackage{subfigure}
\usepackage{booktabs}
\usepackage{hyperref}

\usepackage{tikz}
\usetikzlibrary{calc}
\usetikzlibrary{arrows.meta}

\usepackage{mathtools}

\usepackage[capitalize,noabbrev]{cleveref}

\usepackage{todonotes}

\usepackage[algo2e,lined,boxed,commentsnumbered,ruled]{algorithm2e}

\usepackage{enumitem}
\usepackage{nicefrac}       

\usepackage{comment}

\definecolor{mygray}{gray}{0.8}

\newcommand{\numberthis}{\addtocounter{equation}{1}\tag{\theequation}}

\allowdisplaybreaks

\DeclareMathOperator*{\argmax}{argmax}
\DeclareMathOperator*{\argmin}{argmin}

\newcommand{\lrset}[1]{\left\{ #1 \right\}}
\newcommand{\lrp}[1]{\left( #1 \right)}
\newcommand{\lrs}[1]{\left[ #1 \right]}
\newcommand{\abs}[1]{\left|{#1}\right|}

\newcommand{\KLinf}{\operatorname{KL_{inf}}}
\newcommand{\infKL}{\operatorname{inf-KL}}
\newcommand{\KL}{\operatorname{KL}}

\newcommand{\KLinfDH}{\operatorname{KL^{DH}_{inf}}}
\newcommand{\KLtildeinf}{\operatorname{KL^\dagger_{inf}}}
\newcommand{\KLtildeinvinf}{{\operatorname{KL^{\dagger,-1}_{inf}}}}
\newcommand{\KLtildeinfDH}{\operatorname{KL^{\dagger, DH}_{inf}}}

\newcommand{\tauone}{\tau^{\text{err}}}
\newcommand{\tautwo}{\tau^{\text{gap}}}

\newcommand{\Exp}[2]{\mathbb{E}_{#1}\lrs{#2}}
\newcommand{\V}{\operatorname{Var}}

\newcommand{\R}{\mathbb{R}}
\newcommand{\N}{\mathbb{N}}

\newcommand{\bP}{\mathbb{P}}
\newcommand{\A}{\mathbb{A}}

\newcommand{\calP}{\mathcal P}
\newcommand{\calK}{\mathcal K}
\newcommand{\calE}{\mathcal E}

\newcommand{\calX}{\mathcal X}
\newcommand{\calD}{\mathcal D}

\newcommand{\supp}{\operatorname{supp}}

\newcommand{\kl}{\operatorname{kl}}

\newcommand{\Ber}{\operatorname{Ber}}
\newcommand{\shubhada}[1]{{\color{cyan} #1}}

\DeclarePairedDelimiter\floor{\lfloor}{\rfloor}

\title{
On Stopping Times of Power-one Sequential Tests:\\ Tight Lower and Upper Bounds}
\date{}
\author{Shubhada Agrawal$^1$, Ashwin Ram$^2$, Aaditya Ramdas$^3$\\
$^1$Indian Institute of Science, $^{2,3}$Stanford University\\
$^1$\texttt{shubhada@iisc.ac.in}, $^{2,3}$\texttt{\{ashram03,aramdas\}@stanford.edu}\\
\smallskip \\
\today
}

\begin{document}
\maketitle

\begin{abstract}
We present two general lower bounds for stopping times of sequential tests between arbitrary composite nulls $\mathcal P$ and alternatives $\mathcal Q$. The first lower bound is for the ``Wald setting'' where the type-1 error level $\alpha$ approaches zero for a fixed alternative $Q \in \mathcal Q$, and equals $\log(1/\alpha)$ divided by a certain infimum KL divergence between $\mathcal P$ and $Q$, termed $\KLinf$. The second lower bound applies to the ``Farrell setting'', where $\alpha$ is fixed and $\KLinf$ approaches $0$ along a sequence of alternatives such that the required expected sample size along that sequence is of order at least $\operatorname{\KL^{-1}_{inf}} \log \log \operatorname{\KL^{-1}_{inf}}$. Our main contribution is the generality of these bounds, which hold in non-parametric, composite settings, without requiring a dominating reference measure, substantially generalizing the known parametric results. We also provide sufficient conditions for matching upper bounds and show that these are met in several nontrivial non-parametric cases.
\end{abstract}


\section{Introduction}
Suppose that we observe a stream of i.i.d.\ data $X_1, X_2,\dots$ from some distribution. Given two nonintersecting sets $\mathcal P$ and $\mathcal Q$ of probability distributions on $\R$, we study the problem of testing the null $\mathcal P$ against the alternative $\mathcal Q$. We consider the classical setup of $\alpha$-correct power-one (or one-sided) sequential tests: under the null,  the probability of error is restricted to be at most a given small positive constant $\alpha$, and under the alternative, the goal is to stop after observing as few samples as possible. 

We prove lower and upper bounds on the stopping time $\tau_\alpha$ of any test, that hold for general distribution families. Importantly, beyond the i.i.d.\ assumption made in the first line, our lower bound in the Wald setting requires no additional assumptions on $\mathcal P$ and $\mathcal Q$; most importantly, there need not exist a single common dominating reference measure for both classes, and so densities may not exist and likelihood ratios have to be handled with care. This generality allows us to discuss interesting non-parametric classes of distributions like bounded distributions on a known interval and $1$-sub-Gaussian distributions, as two common examples encountered in the modern literature. 

Further, while most of the existing literature focuses on the ``Wald setting'', where one fixes an alternative $Q$ and studies the expected (or almost sure) stopping time as $\alpha\to0$, we also focus on an interesting ``Farrell setting'' where $\mathcal P$ and $\mathcal Q$ are not separated, and $\alpha$ is fixed. In this regime, we derive an iterated-logarithmic lower bound along certain sequences of alternatives. This hardness is a consequence of several alternatives approaching a common null boundary component, and holds both in parametric and non-parametric settings; examples are provided in the paper. We also provide general sufficient conditions for matching upper bounds in both regimes.


Let us elaborate with two well-known examples. For the Wald setting, if $\mathcal P = \{P\}$ and $\mathcal Q=\{Q\}$ are singletons, consider Wald's one-sided sequential probability ratio test (SPRT) 
$$
\tau_\alpha = \inf\lrset{n: \frac{dQ}{dP}(X_1, \dots, X_n) \geq \frac{1}{\alpha} },
$$
where $\frac{dQ}{dP}$ denotes the likelihood ratio (Radon-Nikodym derivative) between $Q$ and $P$. Then, we have $P(\tau_\alpha < \infty) \leq \alpha$ by Ville's inequality \citep{ville1939etude}, and
\[
\lim_{\alpha \to 0} \frac{\Exp{Q}{\tau_\alpha}}{\log({1}/{\alpha})} = \frac{1}{\KL(Q,P)},
\] 
where $\KL$ is the Kullback-Leibler divergence, and no test can achieve a lower right-hand side. 

For the Farrell setting, consider  $\mathcal P = \lrset{N(m,1): m \le 0}$, $\mathcal Q = \lrset{ N(m,1): m > 0}$, where $N(m,1)$ denotes the Gaussian distribution with mean $m$ and unit variance. \cite{robbins1974expected} derived a mixture SPRT that requires $\frac{2}{\Delta^2}\log\frac{1}{\alpha}$ samples when $\alpha\rightarrow 0$ and $Q=N(\Delta,1)$ is fixed,  but also requires $O(\frac{1}{\Delta^2} \log\log\frac{1}{\Delta})$ samples for a sequence of alternatives in the regime $\Delta \rightarrow 0$, the latter tightness being due to a lower bound by~\cite{farrell1964asymptotic}.

Although such special cases have been studied in the literature, in this work we develop a much more general theory for sequential testing using power-one tests. 


\paragraph*{Contributions} We now make the key contributions of the paper more precise:
\begin{enumerate}
    \item\label{one} 

    For any $Q\in \mathcal Q$, we prove that as $\alpha\to 0$, $\tau_\alpha/(\log(1/\alpha))$ is lower bounded, both $Q$-almost surely and in $Q$-expectation, by  $1/\KLinf(Q,\calP)$   (Theorem~\ref{th:lbsmallalpha}), where $\KLinf(Q,\calP) := \inf\nolimits_{P\in\calP}\KL(Q,P)$. We present sufficient conditions for sequential tests to exactly match this lower bound in the $\alpha\rightarrow 0$ setting (Theorem~\ref{th:generalasupperboundcond}).

    \item \label{two} We then consider composite alternatives $\mathcal Q$ that are not separated from $\mathcal P$, meaning $$\inf_{Q\in\mathcal Q,P\in\mathcal P}\KL(Q,P)=0,$$ despite $\mathcal Q\cap\mathcal P=\emptyset$ (eg., testing negative versus positive mean). While the previous bound keeps $Q$ fixed and lets $\alpha\to0$, we fix $\alpha$ and consider an exponentially spaced (in KL) sequence of alternatives ``directly'' approaching the null (formally defined later). Theorem~\ref{th:lowerbound2} shows that along every such sequence, all but an asymptotically sparse set of alternatives require an expected sample size of order at least $\operatorname{\KL^{-1}_{inf}}\log\log\operatorname{\KL^{-1}_{inf}}$. Corollary~\ref{cor:lilseq} further shows that whenever at least one such direct sequence exists, every $\alpha$-correct power-one test incurs the aforementioned iterated-logarithmic lower bound along some sequence with $\KLinf\to0$. This substantially generalizes the classical result of \citet{farrell1964asymptotic} beyond parametric families, without requiring a common dominating reference measure. We also present sufficient conditions for sequential tests to attain the matching $\operatorname{\KL^{-1}_{inf}}\log\log\operatorname{\KL^{-1}_{inf}}$ order (Theorem~\ref{generalupperboundcond2}).
    
    
    \item\label{three} We finally show that the aforementioned sufficient conditions can be met, and thus that the two derived lower bounds are tight, for a wide range of parametric and non-parametric problems (eg: bounded or sub-Gaussian data). The tightness of both lower bounds is known in specific settings; our upper bounds show tightness in much greater generality. 
\end{enumerate}

\paragraph*{Limitations} Despite the above generality, the two main limitations of our work are that we restrict ourselves to the i.i.d.\ setting, and that we do not examine higher-order notions of optimality in either Wald or Farrell regimes. We think that these are reasonably broad open directions for future work, out of scope of the current paper.

\paragraph*{Paper outline} In Section~\ref{sec:prelim}, we formally introduce the setup and notation, and review the relevant literature. We prove the two lower bounds in Section~\ref{sec:lb}; the Wald regime is discussed in Section~\ref{sec:lb_delta}, and the Farrell regime is presented in Section~\ref{sec:lbsmallgapfixedconf}.

We demonstrate the tightness of the corresponding lower bounds in several representative settings in Section~\ref{sec:ub}; we study non-parametric composite nulls against parametric point and composite alternatives in Sections~\ref{sec:ub_nonpcompnullparampointalt} and~\ref{sec:ub_nonpcompnullpcompalt}, respectively. Section~\ref{sec:testingmeanofbdd} examines the non-parametric composite-vs-composite setting of testing the mean of bounded distributions. Section~\ref{sec:ub_finitelygenhypo} extends our analysis to a setting where hypotheses are generated by a finite set of constraints. We construct a single test simultaneously achieving the two lower bounds in the respective asymptotic regimes in Section~\ref{sec:jointoptimality}, and conclude in Section~\ref{sec:conc}.

All proofs omitted from the main text are in the appendices. In Appendix~\ref{app:kltildeinf}, we also state and prove several properties and concentration results related to $\KLinf$ (introduced in the next section), which may be of independent interest. 

\section{Preliminaries: setup and background} \label{sec:prelim}
Given two non-intersecting sets of probability measures $\cal P$ and $\cal Q$ on a filtered measurable space $(\Omega, \mathcal F)$, consider sequentially testing $\cal P$ (null) against $\cal Q$ (alternative). Let $X_1, X_2, \dots$ denote an i.i.d. sequence of random variables from some distribution in $\mathcal P \cup \mathcal Q$. Given $\alpha > 0$, a one-sided sequential test is simply a stopping time at which one rejects the null (and if the stopping time is infinite, we never reject the null). An $\alpha$-correct, power-one sequential test is a stopping time $\tau_\alpha$ that satisfies:
\begin{align*}  
    &P\lrs{\tau_\alpha < \infty} \le \alpha, \quad \forall P \in \mathcal P,  \tag{$\alpha$-correct} \\
    &Q\lrs{ \tau_\alpha < \infty } = 1, \quad \forall Q\in\mathcal Q. \tag{Power one}
\end{align*}
Let $\mathcal T_\alpha \equiv \mathcal T_\alpha(\mathcal P,\mathcal Q)$ denote the set of all $\alpha$-correct, power-one tests.
The main focus of this work is to derive tight lower bounds on $\Exp{Q}{\tau_\alpha}$, where the expectation is taken under the joint distribution of the observations generated by the unknown  $Q \in \mathcal Q$.



\paragraph*{KL-inf: the separation notion} 
Given two probability measures, the KL divergence between $Q$ and $P$ (denoted $\KL(Q,P)$) is defined as 
$$\KL(Q,P) = \mathbb{E}_{Q}\left[\log\frac{dQ}{dP}(X)\right],$$ 
if $Q \ll P$, and $\infty$ otherwise.
Next, given the null $\cal P$ and $Q \in \cal Q$, we define $$\KLinf(Q,\mathcal P) := \inf\{\KL(Q,P): P \in \mathcal P\}.$$ 
 $\KLinf$ has been studied for specific classes of $\cal P$ and $\cal Q$ in the stochastic multi-armed bandit literature, where optimal algorithms use a certain empirical version of it as a test statistic at each time. We refer the reader to \citet[\S 4]{thesis} for an exposition.

We say that $Q$ is separated from $\cal P$ if $\KLinf(Q,\mathcal P) > 0$. Observe that if $Q = N(m_Q,1)$ and $\mathcal P = \{P\} = \lrset{N(m_P,1)}$, then $\KLinf(Q,\mathcal P) = \KL(Q,P) = 0.5(m_P - m_Q)^2$. The lower bounds that we prove are proportional to $\operatorname{KL^{-1}_{inf}}$ or $\operatorname{KL^{-1}_{inf}} \log\log\operatorname{KL^{-1}_{inf}}$. 
Finally, we say that $\cal P$ and $\cal Q$ are separated if $\inf\nolimits_{Q\in\mathcal Q}~\KLinf(Q,\mathcal P) > 0$.


\paragraph*{Test supermartingales and Ville's inequality} A common approach to constructing an $\alpha$-correct sequential test is to find a nonnegative stochastic process ${M}$ that is a $P$-supermartingale with mean $\Exp{P}{M_1} \le 1$ for all $P \in \cal P$. Such processes ${M}$ are known as $\calP$-test supermartingales. Sequential tests satisfying the above requirements can then be derived by setting $\tau_\alpha = \min\{n: M_n \ge \frac{1}{\alpha}\}$, and $\alpha$-correctness follows from Ville's inequality \citep{ville1939etude}. However, such test supermartingales in general do not suffice: there exist $\mathcal P$ for which every $\calP$-test supermartingale is a nonincreasing process. This means it has  power at most $\alpha$ against any $\cal Q$ even if $\cal Q$ is separated from $\cal P$. It turns out that a more general concept --- an e-process --- is both necessary and sufficient for recovering all sequential tests~\citep{ramdas2022testing,ramdas2020admissible}, and delivering power-one tests whenever they exist, even when $\cal P$ supermartingales may be powerless. We define this next.


\paragraph*{Numeraire e-variable, and e-process} An \emph{e-variable} for $\cal P$ is a nonnegative random variable such that $\sup_{P\in\mathcal P}\Exp{P}{E} \le 1$ (cf. \citet[\S 1]{ramdas2024hypothesis}).
For testing $\cal P$ against a point alternative $Q$, the \emph{numeraire e-variable} is a $Q$-a.s. unique and strictly positive $\cal P$-e-variable $E^*$ such that $\Exp{Q}{E/E^*}\le 1$ for every other $\mathcal P$-e-variable $E$ (cf. \cite{larsson2024numeraire}). The numeraire also satisfies $\Exp{Q}{\log(E^*/E)}\geq 0$ for all other e-variables $E$, a property known as log-optimality. 

While e-variables play an important role in nonsequential testing problems, their sequential analog called \emph{e-processes} are fundamental for sequential testing. A \emph{$\cal P$-e-process} is a nonnegative stochastic process $\{E_n\}$ adapted to the filtration $\cal F$ such that $\Exp{P}{E_\tau} \leq 1$ for all $P \in \mathcal P$ and all stopping times $\tau$ \citep[\S 7]{ramdas2024hypothesis}. Test supermartingales are e-processes (by the optional stopping theorem) but not vice versa.
E-processes are an important generalization of nonnegative test supermartingales; see \cite{ramdas2022testing} for details.
Given a $\calP$-e-process $\lrset{E_n}_{n\ge 1}$, we can define an $\alpha$-correct sequential test as: 
\begin{equation}\label{eq:generalstoppingtime} \tau_{\alpha} = \min\lrset{n: E_n \ge {1}/{\alpha}}. \end{equation}
As earlier, $\alpha$-correctness of $\tau_{\alpha}$ also follows from Ville's inequality \citep{ville1939etude}.
 
\paragraph*{Related work} The field of sequential hypothesis testing began with the introduction of the sequential probability ratio test (SPRT) in the seminal work of \cite{wald1945sequential}, which addressed simple hypothesis testing in a parametric setting. Subsequently, \cite{wald1948optimum} established the optimality of Wald's SPRT by showing that it minimizes the expected number of samples among all tests with a given type-1 error and power.

The framework of $\alpha$-correct power-one tests and their relationship to the law of the iterated logarithm was heavily studied by Robbins and collaborators, for example in \cite{darling1967iterated,darling1967confidence, robbins1968iterated, robbins1970statistical, robbins1970boundary, lai1976confidence}. Most of these works primarily focused on the related problems of constructing confidence sequences and boundary-crossing probabilities, which are essential tools for designing and analyzing sequential tests and their stopping times.

\cite{farrell1964asymptotic} investigated the relationship between power-one tests and the law of the iterated logarithm, proving a lower bound on the expected number of samples required by such tests in parametric settings, in what we call the Farrell regime. Our Farrell-type result substantially generalizes this phenomenon: we show that this iterated-logarithm obstruction is not intrinsically parametric; it is a consequence of many alternatives approaching a common null boundary component, and holds for broad composite and non-parametric hypothesis classes without requiring a common dominating reference measure.

\cite{robbins1974expected} were the first to develop optimal sequential tests for both simple and composite hypothesis testing problems concerning the means of distributions from an exponential family. They derived asymptotic expansions for the expected number of samples required by their proposed tests and showed that these matched the two lower bounds established earlier by \cite{wald1945sequential} and \cite{farrell1964asymptotic}. Later, \cite{pollak1975approximations} presented further asymptotic formulas for the expected sample size of $\alpha$-correct power-one sequential tests for exponential families, analyzing the effects of distributional misspecification and demonstrating the robustness of these expansions. To demonstrate the tightness of the general lower bounds we develop in this work, we prove matching sample complexity upper bounds for a wide range of sequential tests, including composite and non-parametric ones.

Some recent works that have derived bounds on the expected sample complexity of specific sequential tests include \cite{shin2021non-parametric,chugg2023auditing,shekhar2023non-parametric, chen2025optimistic}, as well as some parallel works to ours by  \cite{durand2025power, waudby2025universal}. However, these bounds are either not sharp in one or both of the limiting regimes that we consider in this work, or are studied only for specific classes of problems. 

There are related works studying specific problems in the ``property testing'' literature; these are often discrete settings, predominantly studied in theoretical computer science and information theory. 
For example, \cite{oufkir2021sequential} derive a high-probability lower bound on the number of samples required to test the equality of discrete distributions supported on finitely many points --- an instance of the general hypothesis testing problem considered in this work. Their sample complexity lower bound, when samples are generated from a distribution that is $\epsilon$ far (in total variation distance) from the point null under consideration, is of the form $\Omega(\frac{1}{\epsilon^2} \log \log \frac{1}{\epsilon})$. They also design a sequential test that achieves their lower bound, improving upon the lower and upper bounds of \cite{daskalakis2017optimal} in the same setting.  

In the stochastic multi-armed bandit literature, \cite{lai1985asymptotically} used a change-of-measure argument to derive a lower bound on the number of samples required by a ``good'' bandit algorithm when the arm distributions are parameterized by a single scalar parameter. This result was extended to multiple parametric families by \cite{burnetas1996optimal}. In a related bandit setting, closer to the power-one sequential hypothesis testing framework studied in this paper, \cite{kaufmann2016complexity, garivier2019explore} developed an information-theoretic inequality capturing the change-of-measure arguments. Building on this, \cite{garivier16a} proved a lower bound on the expected number of samples required by a pure-exploration bandit algorithm for arms drawn from an exponential family, and proposed an algorithm that matches this bound, demonstrating its asymptotic tightness in the sense of \cite{wald1945sequential}. These results were later extended to non-parametric settings with bounded-support and heavy-tailed distributions in \cite{agrawal2020optimal, agrawal2021optimal}.  \cite{pmlr-v35-jamieson14,chen2015optimal} derived lower bounds for the pure-exploration bandit algorithms with sub-Gaussian arm distributions, which are tight in the asymptotic regime of \cite{farrell1964asymptotic}.

Despite being originally written in a bandit context, the methods from the above papers can be adapted to sequential testing and yield certain upper and lower bounds for sequential tests in specific settings. Our work substantially generalizes all the above results, proving lower and upper bounds in greater generality than in the union of the above papers. The Wald lower bound applies without any structural assumptions on the underlying hypothesis classes, while the Farrell lower bound extends the classical iterated-logarithm phenomenon to broad composite and non-parametric settings under a mild direct-approach condition, formally described in Section~\ref{sec:lbsmallgapfixedconf}. We revisit some of these classical and recent developments later in the paper, where we compare their proof techniques and results with those in this work.

\paragraph*{Notation} We denote by $\R$ and $\N$ the collection of real numbers and positive natural numbers, respectively. We always use $\cal P$ and $\cal Q$ for the null and alternative hypothesis, $\alpha$ for the error probability,  $M_1[0,1]$ to denote the collection of probability measures supported in $[0,1]$, $M_1(\R)$ for the collection of probability measures supported in $\R$, and $M_1(\calX)$ for probability measures supported in an abstract space $\calX$. For $p, q \in (0,1)$, we use $\kl(p,q)$ to denote the KL divergence between Bernoulli distributions with mean $p$ and $q$, respectively. Given two probability measures $P$ and $Q$ with mean $m_P$ and $m_Q$, respectively, we use $P\ll Q$ ($P\not\ll Q$) to indicate that $P$ is (not) absolutely continuous with respect to $Q$. We write $Q[\cdot]$ and $\Exp{Q}{\cdot}$ to denote the probability and expectation under the probability measure $Q$, and $\V[Q]$ to denote its variance. We say that a measurable set $A$ is $\mathcal P$-negligible if $P(A) = 0$ for all $P\in\mathcal P$. A pointwise property of a measurable function $f$ is said to hold $\mathcal P$-quasi surely if the set where it fails is $\mathcal P$-negligible. {Wherever needed, we endow the space of probability measures with the topology induced by L\'evy metric, which metrizes weak convergence. For probability measures $\{Q_k\}_{k\ge 0}$ and $Q$, we write $Q_k \Rightarrow Q$ to denote weak convergence of $\{Q_k\}$   to $Q$.}

Next, we interpret $\frac{1}{0}$ as unbounded or $\infty$, $\log(0)$ as $-\infty$, $\log(\infty)$ as $\infty$, and $\frac{1}{\infty}$ as $0$. For $a, b \in \R$, we use $a\vee b$ as a shorthand for $\max\{a,b\}$ and $a\wedge b$ for $\min\{a,b\}$. For $K\in\N$, we write $[K]$ for the set $\lrset{1,\dots,K}$. Finally, we use the Bachmann-Landau notation to describe the asymptotic behavior of certain functions. For functions $f(\cdot)$ and $g(\cdot)$, we write $f(x) = O(g(x))$ as $x\rightarrow \infty$ if there exists a constant $c > 0$ and $x_0$ such that $|f(x)| \le c g(x)$ for all $x \ge x_0$, or equivalently, $\limsup\nolimits_{x\rightarrow\infty}~{|f(x)|}/{g(x)} < \infty$. We write $f(x) = \Omega(g(x))$ if $g(x) =O(f(x))$, and use $f(x) = o(g(x))$ to denote $\lim\nolimits_{x\rightarrow\infty}\frac{f(x)}{g(x)} = 0$.

\section{Lower bounds for general hypothesis testing problems}\label{sec:lb}

\subsection{Small error, fixed gap regime}\label{sec:lb_delta}
Recall that  $\KLinf(Q, \mathcal P) := \inf~\lrset{\KL(Q,P) | P \in \mathcal P}.$

\begin{theorem}\label{th:lbsmallalpha}
    For any $\mathcal P$, $Q \in \mathcal Q$ and $\alpha \in (0, 1)$,  every $\tau_\alpha \in \mathcal T_\alpha$ satisfies
    \[ 
    \Exp{Q}{\tau_\alpha} \ge \frac{\log({1}/{\alpha})}{\KLinf(Q,\mathcal P)}, \numberthis\label{eq:lbexpectedineq} 
    \]
    and also
    \[    Q\left(\liminf_{n\to\infty}\frac{\tau_{e^{-n}}}{n}\ge \frac{1}{\KLinf(Q,\calP)}\right)=1.\numberthis\label{eq:aslbongrid}
    \]
    In addition, if the family of stopping times $(\tau_\alpha)_{\alpha\in(0,1/2)}$ is nested, meaning $0<\alpha'\le\alpha<1/2$ implies that $\tau_{\alpha'} \ge \tau_{\alpha}$ pointwise on every sample path, then the following holds:
    \[
    Q\left(\liminf_{\alpha\to 0}\frac{\tau_\alpha}{\log(1/\alpha)}\ge \frac{1}{\KLinf(Q,\calP)}\right)=1. \numberthis\label{eq:aslbineq}
    \]
\end{theorem}
There are connections between our lower bound~\eqref{eq:lbexpectedineq} and others in the sequential analysis literature, and those in the stochastic multi-armed bandit literature, both for the expected regret in the regret-minimization setup (cf. \citet{lai1985asymptotically, burnetas1996optimal}) and the expected sample complexity in the pure exploration setup (cf. \citet{garivier16a, agrawal2020optimal}). 
The hypothesis testing problem at hand can often be viewed as a one-armed bandit problem, and lower bound techniques from the pure exploration literature can also be extended to arrive at the bound in~\eqref{eq:lbexpectedineq}. In the regret-minimization setting, \citet{agrawal2021regret,panda2026regret} show that the corresponding $\KLinf$-based lower bound can be achieved for a broad non-parametric class of heavy-tailed distributions. However, unlike past work, our result places no assumptions on $\calP$ and $\mathcal Q$, and also proves statements~\eqref{eq:aslbongrid} and~\eqref{eq:aslbineq}. Hence, Theorem~\ref{th:lbsmallalpha} can be viewed as a generalization of these previous approaches. More recently, \citet{sethi2026asymptotically} developed analogous lower and upper bounds in the Wald regime for finite-state ergodic Markov chains, with the information complexity given by a stationary-weighted $\KL$ projection. 

Finally, note that the lower bound in~\eqref{eq:aslbineq} holds only for a nested family of tests. To prove~\eqref{eq:aslbineq}, we first derive the bound in~\eqref{eq:aslbongrid} that holds along an exponentially-spaced grid of $\alpha$, and without requiring any assumptions on the stopping times at different confidence levels. The nesting assumption is used to pass from this grid-based statement to the full limit as $\alpha\to 0$. We now sketch the proof for Theorem~\ref{th:lbsmallalpha} and refer the reader to Appendix~\ref{app:th:lbsmallalpha_proof} for a complete proof. 

{
\begin{proof}[Proof sketch]
The expected lower bound in~\eqref{eq:lbexpectedineq} follows from the stopping-time
data-processing inequality: for every $P\in\calP$ with
$\KL(Q,P)<\infty$,
\[
\Exp{Q}{\tau_\alpha}\KL(Q,P)
\ge
\kl\left(Q[\tau_\alpha<\infty],
           P[\tau_\alpha<\infty]\right)
\ge \log(1/\alpha),
\]
and taking the infimum over $P\in\calP$. Next, to see the almost-sure bound on the grid $\alpha=e^{-n}$ in~\eqref{eq:aslbongrid}, fix
$c<1/\KLinf(Q,\calP)$, choose $P_\epsilon\in\calP$ with
$\KL(Q,P_\epsilon)\le \KLinf(Q,\calP)+\epsilon$, and consider the event $\{\tau_{e^{-n}}\le \lfloor cn\rfloor\}$. Clearly, its probability under $P_\epsilon$ is at most $e^{-n}$. A similar change-of-measure argument as in \cite{lai1985asymptotically}, who prove a lower bound for expected regret in the stochastic multi-armed bandit setup with parametric distributions, then gives a summable upper bound on its $Q$-probability, and Borel--Cantelli Lemma \citep{williams1991probability} rules out stopping by time $cn$ infinitely often under $Q$, implying~\eqref{eq:aslbongrid}.  Finally, nesting allows interpolating $\alpha$ between consecutive points, that is, $e^{-(n+1)}<\alpha\le e^{-n}$, which proves the full $\alpha\to0$ statement in~\eqref{eq:aslbineq}. 
\end{proof}
}


\subsection{Small gap, fixed error regime}\label{sec:lbsmallgapfixedconf}
While the lower bound of Theorem~\ref{th:lbsmallalpha} captures the dominant term in $\alpha$ in the sample complexity, in this section, we consider $\alpha \in (0,0.5)$ as a fixed constant. Instead, we consider the setting when there exists a sequence of distributions in $\mathcal Q$ that gets arbitrarily close to $\mathcal P$ in the sense of the function $\KLinf(\cdot, \mathcal P)$ introduced earlier. As in the previous section, $\KLinf(Q, \mathcal P)$ will play a central role in characterizing the sample complexity under $Q$. We will derive a lower bound in the regime when $\inf_{Q \in \cal Q} \KLinf(Q, \mathcal P) = 0$, meaning that there exists a sequence $Q_n \in \mathcal Q$ such that 
$$0 < \KLinf(Q_n, \mathcal P) \rightarrow 0 \quad \text{ as } \quad n \rightarrow \infty.$$
For this reason, the lower bound in this section holds only for problems with composite alternatives that ``kiss'' the null, meaning that despite $\calP \cap \mathcal Q = \emptyset$, they are not separated (eg., the null has nonpositive means, and the alternative has strictly positive means). 
For $\Delta \in (0,1/e)$,  let 
$$F(\Delta) := \frac{1}{\Delta^2} \log\log \frac{1}{\Delta}.$$ 

In order to state our lower bound, it is useful to standardize how \emph{fast} a sequence of alternatives ${\bf Q}:= (Q_1,Q_2,\dots) \in \mathcal Q^\infty$ approaches $\calP$. We will consider sequences whose  $i$-th element is at a $\KL$-distance of $[e^{-2(i+1)},e^{-2i})$ from $\calP$.\footnote{The normalization of the exponential buckets is inessential. For every fixed $a\in (0,1]$, one may replace the bucket $[\exp({-2(i+1)}), \exp({-2i}))$ by $[a\exp({-2(i+1)}), a \exp({-2i}))$, throughout. The same proof applies with only fixed constants changed, and the conclusions of this section are unchanged.}
Formally, consider the set
\[
\mathcal S := \{ {\bf Q} : {\KLinf(Q_i,\mathcal P)} \in [e^{-2(i+1)}, e^{-2i}) ~\forall i\}. \numberthis\label{eq:setS}
\]

To obtain the lower bound in this section, we focus on sequences that ``directly approach'', or remain close to, a common null component called an \emph{anchor}. For $C \ge 1$, define
\[ \mathcal S_C := \{{\bf Q} \in \mathcal S: \exists P_{\bf Q} \in \calP \text{ such that } \KL(Q_i, P_{\bf Q}) \le C\KLinf(Q_i, \calP)~\forall i\}. \]
We say that $\bf Q$ \emph{directly approaches the null} or that $\bf Q$ is a \emph{direct sequence} if ${\bf Q} \in \mathcal S_C$ for some $C \ge 1$. The anchor $P_{\bf Q}$ is allowed to depend on the sequence $\bf Q$, while the constant $C$ cannot. 

For $c > 0$, $N\in\N$, $\mathbf Q\in\mathcal S_C$, and $\tau_\alpha \in \mathcal T_\alpha$, define the counting function
\[ \nu(\tau_\alpha, c, N,{\bf Q}) \!:=\!  \sum\limits_{i=1}^N \! {\bf 1}\lrp{ \Exp{Q_i}{\tau_\alpha} < cF\lrp{\sqrt{\KLinf(Q_i,\mathcal P)}} }.\] 
In words, $\nu$ counts the number of ``easy instances'' within the first $N$ elements of  $\bf Q$. We say that $\tau_\alpha$ finds a sequence $\bf Q$ ``hard'' if it has very few easy instances, in the sense that for some constant $c$, we have $\nu(\tau_\alpha, c, N,{\bf Q})=o(N)$. One consequence of our main result is:
\begin{quote} 
    For any $C \ge 1$, every sequence ${\bf Q} \in \mathcal S_C$ is hard for every $\alpha$-correct power-one test $\tau_\alpha \in \mathcal T_\alpha$. 
\end{quote}
In fact, we prove an even stronger result below.


\begin{theorem}\label{th:lowerbound2}
    Consider $\cal P,\cal Q$ such that $\inf_{Q\in\mathcal Q}\KLinf(Q,\mathcal P) = 0$. Fix  $\alpha\in (0,0.5)$, $C \ge 1$,  and  $\gamma >0$. If $\mathcal S_C \neq \emptyset$ and $\mathcal T_\alpha \ne \emptyset$, then 
    \[
     \limsup\limits_{N\rightarrow \infty}~ \sup\limits_{ \tau_\alpha \in \mathcal T_\alpha}~ \sup_{{\bf Q} \in \mathcal S_C} ~ \frac{\nu(\tau_\alpha, \frac{\gamma}{4C}, N,{\bf Q})}{N^\gamma} = 0.
    \]
\end{theorem}
Informally, the above theorem states that for every $\tau_\alpha$ and a direct sequence $\bf Q$, the number of easy instances grows not only as $o(N)$ but also smaller than any fractional power of $N$. Hence, any $\alpha$-correct sequential test requires $\Omega(F(\sqrt{\KLinf(Q,\calP)}))$ samples in expectation, for most distributions $Q$ along every sequence that has several alternative elements approaching a common null boundary component. 


\begin{figure}[t]
\centering
\begin{tikzpicture}[
    x=1.05cm,
    y=0.9cm,
    point/.style={
        circle,
        fill=black,
        inner sep=1.5pt
    },
    projected point/.style={
        circle,
        fill=black!65,
        inner sep=1.2pt
    },
    direct path/.style={
        thick,
        -{Stealth[length=2.2mm,width=1.5mm]}
    },
    nondirect path/.style={
        thick,
        dashed,
        -{Stealth[length=2.2mm,width=1.5mm]}
    },
    projection/.style={
        thin,
        densely dotted,
        black!60
    },
    sequence label/.style={
        font=\small
    }
]

\def\rho{0.475}
\def\Qscale{7.6}
\def\Rscale{5.8}
\def\Rheight{2.9}

\begin{scope}[xscale=0.8,yscale=0.8]
\path[fill=black!4]
    (-2.35,3.35)
    .. controls (-0.65,3.30) and (0,3.00) .. (0,2.40)
    -- (0,-2.20)
    .. controls (0,-2.85) and (-0.80,-3.20) .. (-2.00,-3.45)
    .. controls (-3.25,-2.05) and (-3.35,2.00) .. (-2.35,3.35)
    -- cycle;

\draw[thick]
    (-2.35,3.35)
    .. controls (-0.65,3.30) and (0,3.00) .. (0,2.40)
    -- (0,-2.20)
    .. controls (0,-2.85) and (-0.80,-3.20) .. (-2.00,-3.45);

\draw[thick,dashed]
    (-2.35,3.35)
    .. controls (-3.35,2.00) and (-3.25,-2.05) .. (-2.00,-3.45);

\node[font=\Large] at (-1.65,-0.15) {$\mathcal P$};
\end{scope}

\coordinate (anchor) at (0,0);
\node[point] at (anchor) {};
\node[below left=0pt] at (anchor) {$P_{\mathbf Q}$};

\draw[direct path] (8.05,0) -- ($(anchor)+(0.05,0)$);

\foreach \i in {1,...,5} {
    \pgfmathsetmacro{\xq}{\Qscale*pow(\rho,\i-1)}
    \coordinate (Q\i) at (\xq,0);
    \node[point] at (Q\i) {};
}

\node[below=4pt] at (Q1) {$Q_1$};
\node[below=4pt] at (Q2) {$Q_2$};
\node[below=4pt] at (Q3) {$Q_3$};
\node[below=4pt] at (Q4) {$Q_4$};
\node[below=4pt] at (Q5) {$Q_5$};

\node[sequence label] at (5.65,-0.62)
    {direct sequence $\mathbf Q$};

%
%
\draw[nondirect path]
    plot[
        domain=1:0.055,
        samples=100,
        variable=\t
    ]
    ({\Rscale*\t*\t},{\Rheight*\t});

\foreach \i in {1,...,5} {
    \pgfmathsetmacro{\tr}{pow(\rho,0.5*(\i-1))}
    \pgfmathsetmacro{\xr}{\Rscale*\tr*\tr}
    \pgfmathsetmacro{\yr}{\Rheight*\tr}
    \coordinate (R\i) at (\xr,\yr);
    \node[point] at (R\i) {};
}

\node[below=2pt] at (R1) {$R_1$};
\node[above=4pt]       at (R2) {$R_2$};
\node[above=4pt]       at (R3) {$R_3$};
\node[right=4pt]       at (R4) {$R_4$};
\node[below right=1pt] at (R5) {$R_5$};

\node[sequence label,align=left] at (5.65,2)
    {non-direct sequence $\mathbf R$};

\foreach \i in {2,3,4,5} {
    \coordinate (Pi\i) at (0,0 |- R\i);
    \draw[projection] (Pi\i) -- (R\i);
    \node[projected point] at (Pi\i) {};
}

\node[left=4pt] at (Pi3) {$P_i^\star$};

\end{tikzpicture}

\caption{
The sequence $\mathbf Q$ directly approaches the null with anchor $P_{\mathbf Q}$. The sequence $\mathbf R$ is a non-direct sequence. Although the elements $R_i$ converge to $P_{\mathbf Q}$, their $\KL$-projections $P_i^\star$ move along the boundary, and no fixed $P\in\calP$ remains within a constant factor of $\KLinf(R_i, \calP)$. In particular, $\KL(R_i,P_{\mathbf Q})/
 \KLinf(R_i,\mathcal P)\to\infty$.}
\label{fig:direct-and-nondirect}
\end{figure}

\begin{remark} The only real restriction of the above theorem is that $\mathcal S_C$ be nonempty, which we should typically expect to be true for large enough $C$. 
In special cases, one can even take $C=1$, for example, when the null is a singleton. One can also often take $C=1$ in many regular parametric and non-parametric examples. For example, testing unit-variance Gaussians with nonpositive versus positive mean, one can take a sequence of Gaussians with positive means that progressively shrink towards zero, with $N(0,1)$ being the anchor. 

The condition of existence of a direct sequence is also non-vacuous in the fully non-parametric bounded-mean testing problem, considered in Section~\ref{sec:testingmeanofbdd}, where for $m_0 \in (0,1)$, 
$$\calP:=\{P \in M_1[0,1]: m_P = m_0\} \quad \text{ and }\quad \mathcal Q := \{Q\in M_1[0,1]: m_Q < m_0\}.$$ 
For $p\in (0,1)$, let $\Ber(p) := (1-p)\delta_0+p\delta_1$ denote the Bernoulli distribution with mean $p$. Clearly, for $q<m_0$, $\Ber(q)\in\mathcal Q$.  For $P\in\mathcal P$ and $q < m_0$, from data-processing inequality,
\[\KL(\Ber(q),P)\geq \KL(\Ber(q),\Ber({m_0})).\]
Consequently, $\KLinf(\Ber(q),\mathcal P)=\KL(\Ber(q),\Ber({m_0}))$, with $\Ber({m_0})$ serving as a common optimizer for every $q<m_0$.  

Since $q\mapsto\KL(\Ber(q),\Ber({m_0}))$ is continuous and decreases to zero as $q\uparrow m_0$, for any $a\in(0,\min\{1,-\log(1-m_0)\})$, one may choose $q_i\uparrow m_0$ such that $\KL(\Ber({q_i}),\Ber({m_0}))=ae^{-2i}$.  Thus, the sequence $(\Ber({q_i}))_{i\geq1} \in \mathcal S_1$, and hence, directly approaches the null with $\Ber({m_0})$ being the anchor.  

Figure~\ref{fig:direct-and-nondirect} gives a schematic picture of a direct sequence and a non-direct sequence.
\end{remark}

\begin{remark} A completely assumption-free Farrell lower bound is not true. To see this, we consider the following two-point distributions on $\R$ supported on mutually-disjoint points. 

For each $i\ge 1$, let 
    \[ P_i := 0.5 \delta_{2i} + 0.5 \delta_{2i+1}, \qquad \text{ and } \qquad Q_i = (1-q_i) \delta_{2i} + q_i \delta_{2i+1}, \]
    where $q_i \in (0.5,1)$ is chosen so that $\KL(Q_i, P_i) = e^{-2(i+1)}$. Let $\calP := \{P_i: i\ge 1\}$ and $\mathcal Q := \{Q_i : i \ge 1\}$. Then the sequence ${\bf Q} = (Q_1, Q_2, \dots)$ lies in the set $\mathcal S$, defined in~\eqref{eq:setS}.  Nevertheless, it is not a direct sequence as no fixed $P\in\calP$ can be an anchor. In fact, there is no direct sequence in $\mathcal Q$. We now show that the Farrell lower bound is not valid for this problem. 
    
    Note that since the supports are disjoint, $X_1$ itself uniquely determines the index $i = \floor{X_1/2}$ of the data-generating distribution, reducing the problem to a point-vs-point setting. Consider 
    $$\tau_\alpha :=\inf\lrset{n\geq1: L_n^{(\lfloor X_1/2\rfloor)}\geq \frac{1}{\alpha}}, \qquad \text{where} \qquad  L_n^{(i)} :=\prod_{t=1}^n\frac{dQ_i}{dP_i}(X_t).$$
    
    Under $Q_i$, the log increments of $L^{(i)}_n$ have positive mean, given by $\mathrm{KL}(Q_i,P_i)=e^{-2(i+1)}$, so $\tau_\alpha<\infty$ almost surely. $\alpha$-correctness follows from Ville's inequality. Next, since the log-likelihood increment is bounded above by $\log(2q_i)\leq\log 2$, we get $$\log L_{\tau_\alpha\wedge n}^{(i)}\le \log(1/\alpha)+\log 2 \quad \text{ and } \quad \mathbb{E}_{Q_i}[{\log L_{\tau_\alpha\wedge n}^{(i)}}] =\KL(Q_i,P_i)\Exp{Q_i}{\tau_\alpha\wedge n},$$ where the second equality follows from Wald's identity. Letting $n\to\infty$ and applying the monotone convergence theorem \citep{williams1991probability} yields $\mathbb{E}_{Q_i}[\tau_\alpha] \le (\log(1/\alpha)+\log 2)/\KL(Q_i,P_i)$. 
\end{remark}

\begin{corollary}\label{cor:lilseq} Fix $\alpha\in (0,0.5)$ and for $\cal P,\cal Q$ with $\inf_{Q\in\mathcal Q}\KLinf(Q,\mathcal P) = 0$, suppose  there exists $C\ge 1$ such that $\mathcal S_C \ne \emptyset$. Then, for every $\tau_\alpha \in \mathcal T_\alpha$, 
    \[ \limsup\limits_{Q\in\mathcal Q:\KLinf(Q,\calP) \rightarrow 0}~ \frac{\Exp{Q}{\tau_\alpha}}{F(\sqrt{\KLinf(Q,\calP)})} > 0. \]
\end{corollary}
Theorem~\ref{th:lowerbound2} is significantly stronger than Corollary~\ref{cor:lilseq}. The corollary asserts the existence of a hard sequence of alternatives approaching the null, along which every $\alpha$-correct power-one test incurs the Farrell-type lower bound. In contrast, Theorem~\ref{th:lowerbound2} shows that every direct sequence is hard. For completeness, we present a proof for the corollary in Appendix~\ref{app:proof_cor:lilseq}. 

Corollary~\ref{cor:lilseq} itself substantially generalizes the lower bound of \citet{farrell1964asymptotic}, who considered one-sided testing in a single-parameter exponential family (SPEF). Indeed, for a regular nondegenerate SPEF with the null and alternative lying on opposite sides of a boundary parameter $\theta_0$, that is, $\calP = \{P_\theta: \theta \le \theta_0\}$ and $\mathcal Q = \{P_\theta: \theta > \theta_0\}$, the map
\[ \theta \to \KLinf(P_{\theta}, \calP) = \KL(P_\theta, P_{\theta_0}) \]
is continuous and strictly increasing for $\theta > \theta_0$. Hence, one can choose one alternative from every bucket, as desired; the resulting sequence directly approaches the null, with anchor $P_{\theta_0}$. The corollary applies whenever there exists at least one direct sequence, without otherwise requiring a parametric structure or a common dominating reference measure.


Along the direct sequences covered by Theorem~\ref{th:lowerbound2} and Corollary~\ref{cor:lilseq}, the lower bound
$$\Omega\lrp{\frac{1}{\KLinf(Q,\mathcal P)} \log\log\frac{1}{{\KLinf(Q,\mathcal P)}}},$$ 
dominates the earlier lower bound of ${\log(1/\alpha)}/{\KLinf(Q,\mathcal P)}$ when focusing on the $\KLinf$ term. The two lower bounds dominate in complementary regimes of $\alpha \to 0$ and $\KLinf\to 0$. 

Other than~\cite{farrell1964asymptotic}, a few other recent works develop similar bounds for specific problems. \cite{chen2015optimal} developed a bound similar to that of Theorem~\ref{th:lowerbound2} for Gaussian mean-testing. More recently, \cite{oufkir2021sequential, daskalakis2017optimal} developed a lower bound similar to that in Corollary~\ref{cor:lilseq} (with total variation distance instead of $\KLinf$ or $\KL$) for testing equality of discrete distributions supported on finitely many points. The setting they consider has a point null and  composite alternative. To our knowledge, no previous result establishes a lower bound in the generality that we do.



We now prove Theorem~\ref{th:lowerbound2}. 
  For brevity, let $\Delta_{Q,\mathcal P} := \sqrt{\KLinf(Q,\calP)}$. We will use a proof by contradiction, for which we need a few additional results, which we first present. The following change-of-measure type lemma relates the probability of a given event under different distributions. 

\begin{lemma}\label{lem:changeofmeasure}
    Consider any two probability measures $Q$ and $P$. Let $\tau$ denote a stopping time, and let $\calE$ be an $\mathcal F_\tau$-measurable event such that ${Q}[\calE] \ge 0.5$ and ${P}[\calE] \le 0.5$. Then, ${P}[\calE] \ge 0.25 e^{-2\Exp{Q}{\tau} \KL(Q, P)}.$
\end{lemma}
We refer the reader to Appendix~\ref{app:proof_changeofmeasure} for a proof, which follows from an application of the data processing inequality (Lemma~\ref{lem:DP}), and the monotonicity of the $\KL$ divergence between two Bernoulli distributions in the mean parameter.  We will construct an event that has a non-negligible probability under any distribution $Q\in \mathcal Q$, but is rare under $\calP$. We will then use a version of Lemma~\ref{lem:changeofmeasure} to arrive at a contradiction. To this end, for a universal constant $d_l > 0$ (to be chosen later), and $\epsilon > 0$ such that $\alpha < 0.5-\epsilon$, define events: 
\[ \calE_U := \lrset{\tau_\alpha < \infty}, \quad \text{and} \quad  \calE(\Delta) := \calE_U \bigcap \lrset{ {d_l}/{\Delta^2} \le \tau_\alpha \le \Exp{Q}{\tau_\alpha}/\epsilon }. \numberthis\label{eq:eventedelta} \]
Note that the existence of $\epsilon > 0$ is guaranteed since $\alpha < 0.5$. 

\begin{lemma}\label{lem:highprobevent}
    Fix $\epsilon > 0$ such that $\alpha < 0.5-\epsilon$, and $Q\in\mathcal Q$ with $\Delta_{Q,\mathcal P} > 0$.  Then, for $d_l < 0.5 \kl(0.5 - \epsilon,\alpha)$, we have $Q[\calE(\Delta_{Q,\mathcal P})] \ge 0.5$.
\end{lemma}

Lemma~\ref{lem:highprobevent} above argues that for $Q\in\mathcal Q$, the event $\calE(\Delta_{Q,\calP})$ defined in~\eqref{eq:eventedelta} occurs with a non-negligible probability under $Q$. Note that since $\tau_\alpha$ is a power-one test, $Q[\calE_U] = 1$. Thus, it only remains to argue that the second event in the intersection in~\eqref{eq:eventedelta} is not rare under $Q$. To this end, observe that the upper bound on $\tau_\alpha$ holds with probability at least $1-\epsilon$ by Markov's inequality. From Theorem~\ref{th:lbsmallalpha}, we also expect that the lower bound on $\tau_\alpha$ holds with a constant probability, since it is smaller than $\log(1/\alpha)/\KLinf(Q,\calP)$, at least in the limit of $\alpha\rightarrow 0$. The proof in Appendix~\ref{app:proof_highprobevent} formalizes all these arguments. 

Since $P[\calE(\Delta_{Q,\mathcal P})] \leq P(\calE_U) \leq \alpha < 0.5$ for all $P \in \calP$, 
Lemma~\ref{lem:changeofmeasure} then gives  a lower bound on $P[\calE(\Delta(Q,\calP))]$ for $P\in\calP$. However, this may still not lead to any contradiction. In  Lemma~\ref{lem:sumprobbound} below, we show that if $\calE_U$ is a union of $n$ disjoint events $\{\calE(\Delta_i)\}_{i=1}^n$, each of which occurs with a non-zero probability under $P\in\calP$, then inequality~\eqref{eq:tobecontradicted} is enforced, which we later use to prove Theorem~\ref{th:lowerbound2} via contradiction.

\begin{lemma}\label{lem:sumprobbound}
    Fix $\epsilon > 0$ and $\alpha < 0.5-\epsilon$, and let $d_l < 0.5\kl(0.5 - \epsilon,\alpha)$. Fix any $\tau_\alpha \in \mathcal T_\alpha$, and constants $c > 0$ and $C\ge 1$. Consider a finite sequence (of $n \ge 1$ elements) $Q_i \in \mathcal Q$ and define $\Delta^2_i := \KLinf(Q_i,\mathcal P)$. Suppose the following hold. 
    \begin{enumerate}
    \item $\Delta_i > 0$ are strictly decreasing, and there exists a common anchor $P_0 \in \calP$ such that $\KL(Q_i, P_0) \le C \Delta^2_i$ for all $i\in [n]$. \item The events $ \calE(\Delta_i) = \calE_U \cap \{ 
    {d_l}/{\Delta^2_i}\le \tau_\alpha < \Exp{Q_i}{\tau_\alpha}/\epsilon\}$ are disjoint, and $\Exp{Q_i}{\tau_\alpha} \le c F(\Delta_i)$ for all $i\in [n]$.
    \end{enumerate}
    Then, 
    \[ \sum\limits_{i=1}^n e^{-2c CF(\Delta_i) \Delta_i^2} \le 4 \alpha. \numberthis\label{eq:tobecontradicted} \]
\end{lemma}
\begin{proof}
    For any $P\in \mathcal P$, we have \begin{align*}
    P\lrs{\calE_U} \ge P\lrs{\cup_{i=1}^n \calE(\Delta_i) } = \sum\limits_{i=1}^n P\lrs{ \calE(\Delta_i) }&\overset{(a)}{\ge} \frac{1}{4} \sum\limits_{i=1}^n e^{-2 \Exp{Q_i}{\tau_\alpha}\KL(Q_i,P)}\\
    &\ge \frac{1}{4} \sum\limits_{i=1}^n e^{-2 c F(\Delta_i)  \KL(Q_i,P)}, 
    \end{align*}
    where inequality $(a)$ follows from Lemma~\ref{lem:changeofmeasure}, and needed Lemma~\ref{lem:highprobevent} to conclude that ${Q_i}[\calE(\Delta_i)] \ge \frac{1}{2}$ for each $i$ along with the observation that $P[\calE(\Delta_i)] \le \alpha < 0.5$.
    Now, since $P[\calE_U]\le \alpha$, we have
    \[ \alpha \ge \frac{1}{4} \sum\limits_{i=1}^n e^{-2 c F(\Delta_i)  \KL(Q_i,P)}. \]
    Since the above inequality holds for an arbitrary $P \in \calP$, it holds for $P_0$ in particular. Thus,  
    \begin{align*}
        \alpha 
        &\ge \frac{1}{4} \sum\limits_{i=1}^n e^{-2 c F(\Delta_i)  \KL(Q_i,P_0)} \ge \frac{1}{4} \sum\limits_{i=1}^n e^{-2 c C F(\Delta_i)  \KLinf(Q_i,\calP)},
    \end{align*}
    where the last inequality follows since $ \KL(Q_i,P_0) \le C\KLinf(Q_i,\calP)$, for all $i\in [n]$. We get the desired inequality by rearranging the above.
\end{proof}


Lemma~\ref{lem:suffdisjointcond} below gives sufficient conditions to construct the disjoint events $\calE(\Delta_{i})$ needed in the previous lemma. 

\begin{lemma}\label{lem:suffdisjointcond}
    Fix $\epsilon > 0$ such that $\alpha < 0.5 - \epsilon$,   $\tau_\alpha \in \mathcal T_\alpha$, and a universal constant $c> 0$. Let $\Delta^2_i := \KLinf(Q_i, \mathcal P)$ for $Q_i \in \mathcal Q$ and $i\in [n]$ be a finite length sequence that satisfies
    \begin{enumerate}
        \item $\frac{1}{e} > \Delta_1 > \Delta_2 > \dots > \Delta_n \ge \beta > 0$, for some $\beta > 0$, 
        \item for all $i\in[n]$, $\Exp{Q_i}{\tau_\alpha} \le c F(\Delta_i)$,
        \item for $L_i := \ln \Delta^{-1}_i$, we have 
        \[ L_{i+1} - L_i > \frac{1}{2} \ln\ln\ln \beta^{-1} + c_1, \quad \text{where} \quad c_1 = \frac{\ln c + \ln \frac{1}{\epsilon} - \ln d_l}{2}. \]
        Then the events ${\calE(\Delta_1), \calE(\Delta_2), \dots, \calE(\Delta_n)}$ are disjoint.
    \end{enumerate}
\end{lemma}
\begin{proof}
    For events $\calE(\Delta_i)$ to be disjoint, we only need to show that the corresponding intervals of $\tau_\alpha$ are disjoint. We will show that this is the case for two adjacent intervals corresponding to $\Delta_{i}$ and $\Delta_{i+1}$, that is,
    \[ \lrs{\frac{d_l}{\Delta_i^2}, \frac{1}{\epsilon} \Exp{Q_i}{\tau_\alpha}} \quad \text{and} \quad \lrs{ \frac{d_l}{\Delta_{i+1}^2}, \frac{1}{\epsilon} \Exp{Q_{i+1}}{\tau_\alpha}} \]
    are disjoint. In fact, since $\frac{\Exp{Q_i}{\tau_\alpha}}{\epsilon} \le  \frac{c}{\epsilon} F(\Delta_i)$, it suffices to show
    \[ \frac{c}{\epsilon}  F(\Delta_i) < \frac{d_l}{\Delta_{i+1}^2}. 
    \]
    The above is equivalent to showing
    \(  \ln \frac{1}{\epsilon} + \ln c + 2 L_i + \ln \ln L_i < \ln d_l + 2 L_{i+1}, 
    \)
    which is further equivalent to showing
    \[ L_{i+1} - L_i > \frac{\ln \frac{1}{\epsilon} + \ln c - \ln d_l }{2} + \frac{1}{2}\ln\ln L_i. \]
    The above follows from the third condition in the lemma statement, with the observation that 
    \[ L_i := \ln \frac{1}{\Delta_i} \le \ln\frac{1}{\beta}  \] 
     since $\beta \le \Delta_i$, for all $i$. This concludes the proof.
\end{proof}

We can now finally prove Theorem~\ref{th:lowerbound2}. Observe from Lemma~\ref{lem:sumprobbound} that for the left-hand side in~\eqref{eq:tobecontradicted} to be  smaller than $4\alpha$, either $n$ needs to be small, or each of the terms in the summation needs to be small. We will construct a long sequence (large $n$) of alternative distributions $Q_i\in \mathcal Q$ (or equivalently, $\Delta_i$) that meets the conditions  in Lemma~\ref{lem:sumprobbound}. Hence,~\eqref{eq:tobecontradicted} holds for this sequence. However, if Theorem~\ref{th:lowerbound2} is violated, the terms $cF(\Delta_i)$ in the exponents of~\eqref{eq:tobecontradicted} are small, and hence~\eqref{eq:tobecontradicted} is violated. This yields Theorem~\ref{th:lowerbound2} by contradiction. 
\subsubsection{Proof of Theorem~\ref{th:lowerbound2}}

\begin{proof}
    Fix $C\ge 1$ and $\gamma > 0$ as in the theorem statement. If $\gamma > 1$, then the theorem follows immediately since $\nu(\tau_\alpha, \tfrac{\gamma}{4C}, N, {\bf Q}) \le N$ uniformly over $\tau_\alpha$ and $\bf Q$. Hence, it only remains to consider $0 < \gamma \le 1$. 
    
    Suppose towards a contradiction that for some $0 < \gamma \le 1$, 
    \[ \limsup\limits_{N\rightarrow \infty}~ \sup\limits_{\tau_\alpha \in \mathcal T_\alpha}~ \sup\limits_{{\bf Q} \in \mathcal S_C}~ \frac{\nu(\tau_\alpha, \frac{\gamma}{4C}, N, {\bf Q})}{N^\gamma} > 0. \numberthis\label{eq:lowerboundcontassump} \]
    We will demonstrate a sequence of ${\bf Q} \in \mathcal S_C$ such that for $\Delta_i := \sqrt{\KLinf(Q_i, \mathcal P)}$, $\calE(\Delta_i)$ are disjoint and $\Exp{Q_i}{\tau_\alpha} \le \tfrac{\gamma}{4C} F(\Delta_i)$, that is, the conditions of Lemma~\ref{lem:sumprobbound} are satisfied. However, along with the assumed~\eqref{eq:lowerboundcontassump}, the inequality~\eqref{eq:tobecontradicted} is violated.
    
    To make the above concrete, first see that if the inequality in~\eqref{eq:lowerboundcontassump} is  satisfied, then there exists a sequence $\{N_i\}_i$ increasing to $\infty$, and a positive constant $\beta > 0$ such that 
    \[ \sup\limits_{\tau_\alpha \in \mathcal T_\alpha}~\sup\limits_{{\bf Q} \in \mathcal S_C}~\frac{\nu(\tau_\alpha, \tfrac{\gamma}{4C}, N_i, {\bf Q})}{N^\gamma_i} > 2\beta, \quad \forall i. \]
    Fix a large enough $N_i$. Then, the above inequality implies that there exists an $\alpha$-correct test (specified by $\tau_\alpha \in \mathcal T_\alpha$), and a direct sequence (specified by ${\bf Q}\in \mathcal S_C$), such that $\nu(\tau_\alpha, \frac{\gamma}{4C}, N_i, {\bf Q}) \ge \beta N^\gamma_i$. 

    Next, since ${\bf Q} \in \mathcal S_C$, $\Delta_r := \sqrt{\KLinf(Q_r, \calP)} \in [e^{-(r+1)}, e^{-r})$ for all $r\ge 1$, and there exists a common anchor $P_{\bf Q} \in \calP$ such that for all $r\ge 1$, $\KL(Q_r, P_{\bf Q}) \le C \KLinf(Q_r, \calP)$. 
    
    We will now carefully pick a subsequence $\{Q_j\}$ of ${\bf Q}$ that satisfies the conditions of Lemma~\ref{lem:sumprobbound}. Let 
        \[ S:= \lrset{j\in \lrset{2, \dots, N_i}: \Exp{Q_j}{\tau_\alpha} < \tfrac{\gamma}{4C} F(\Delta_j)}. \]
    Thus, $|S|  \ge \beta N^\gamma_i-1$. Unfortunately, the constructed sequence $\{Q_j\}_{j\in S}$ may not have disjoint $\calE(\Delta_j)$ as needed in Lemma~\ref{lem:sumprobbound}. To this end, we pick a subsequence, as discussed next.
    
    Let $b\in \N$ be given by
    \[ b = \left\lceil\frac{\ln \frac{\gamma}{4C} + \ln \frac{1}{\epsilon} - \ln d_l + \ln\ln N_i}{2} + 1 \right\rceil. \]
    Keep only $1^{st}$, $(1+b)^{th}$, $(1+2b)^{th}$, \dots elements in $S$ and drop the others (overloading notation, we still call this thinner set $S$). Let us map these remaining elements in $S$ to the corresponding $\{\Delta_i\}_{i=1}^{|S|}$, sorted in decreasing order. Then, 
    ${1}/{e} > \Delta_1 > \Delta_2 > \dots > \Delta_{|S|} \ge e^{-(N_i+1)} > 0.$ 
    
    By construction of $S$, we also have
    \[ \ln \Delta^{-1}_{i+1} - \ln \Delta^{-1}_i > \frac{\ln \frac{\gamma}{4C} + \ln \frac{1}{\epsilon} - \ln d_l}{2} + \frac{1}{2}\ln\ln N_i. \]
    From Lemma~\ref{lem:suffdisjointcond}, we see that the events $\calE(\Delta_i)$ for $i\in [|S|]$ are now disjoint, and satisfy the conditions of Lemma~\ref{lem:sumprobbound}. Moreover, $|S| \ge (\beta N^\gamma_i-1)/b$, which implies (for large enough $N_i$ which we can choose since the sequence $\{N_i\}_i$ increases to $\infty$) that 
    $$|S| \ge \beta N^\gamma_i/ \ln\ln N_i.$$ 
    Further, since $\{Q_j\}_{j \in S}$ is a subsequence of ${\bf Q}$, it inherits the common anchor $P_{\bf Q}$ and satisfies $\KL(Q_j, P_{\bf Q}) \le C \KLinf(Q_j, \calP)$, for all $j \in S$. Now consider the following calculation:
    \begin{align*}
        \sum\limits_{j=1}^{|S|} e^{-\frac{2\gamma}{4C} CF(\Delta_j)\KLinf(Q_j,\mathcal P)} 
        &= \sum\limits_{j=1}^{|S|} e^{-\frac{\gamma}{2} \ln\ln \frac{1}{\sqrt{\KLinf(Q_j,\mathcal P)}}} \\
        &= \sum\limits_{j=1}^{|S|}\lrp{\ln\frac{1}{\sqrt{\KLinf(Q_j,\mathcal P)}}}^{-\frac{\gamma}{2}} \\
        &\ge |S| (N_i+1)^{-\frac{\gamma}{2}} \tag{$\Delta_j = \sqrt{\KLinf(Q_j,\mathcal P)} \ge e^{-(N_i+1)}$, for all $j$} \\
        &\ge \frac{\beta N^\gamma_i}{\ln\ln N_i} (N_i+1)^{-\frac{\gamma}{2}}.
    \end{align*}
Now, since $ 0 <  \gamma \le 1$, we can choose $N_i$ large enough such that 
$$\frac{\beta N^\gamma_i}{\ln\ln N_i} (N_i+1)^{-\frac{\gamma}{2}} > 4 \alpha,$$
which contradicts Lemma~\ref{lem:sumprobbound}.  
\end{proof}

\section{Upper bounds: optimal sequential tests demonstrating tightness of lower bounds}\label{sec:ub}
\cite{robbins1974expected} studied the problem of testing the mean of a single-parameter exponential family of distributions in the two asymptotic regimes we discussed in the previous sections. They established asymptotic expansions for the expected stopping times of the $\alpha$-correct power-one sequential tests that they proposed for the special case of Gaussian distributions with unit variance. In this case, for the point-versus-point setting of $\mathcal P= \{N(0,1)\}$ vs. $\mathcal Q= \{N(m, 1)\}$ for a given known $m\ne 0$,  their sequential test requires $\frac{2}{m^2}\log\frac{1}{\alpha}$ samples in expectation in both the asymptotic regimes: $m\rightarrow 0$ and $\alpha \rightarrow 0$ \cite[Lemma 5]{robbins1974expected}. This demonstrates the tightness of the lower bound in Theorem~\ref{th:lbsmallalpha}. This also establishes that the lower bound in Theorem~\ref{th:lowerbound2} holds only for composite alternatives. 

Next, for composite null and alternative ($\mathcal P: m \le 0$ vs. $\mathcal Q: m > 0$ for Gaussian distributions with unit variance), in the same work, the authors proposed a sequential test that requires $O(\frac{1}{m^2} \log\log\frac{1}{|m|} )$ samples in the $m\rightarrow 0$ setting \cite[Theorem 2]{robbins1974expected} and $\frac{2}{m^2}\log\frac{1}{\alpha}$ \cite[Lemma 5]{robbins1974expected} samples in the $\alpha\rightarrow 0$ setting. This demonstrates the tightness of the lower bounds developed in Theorems~\ref{th:lbsmallalpha} and~\ref{th:lowerbound2} for the (single) parametric setting. 

In this section, we show that the two lower bounds from Theorems~\ref{th:lbsmallalpha} and~\ref{th:lowerbound2} are tight for several non-parametric and composite examples. Towards this, in Theorems~\ref{th:generalasupperboundcond} and~\ref{generalupperboundcond2}, we first present sufficient conditions on sequential tests for the general hypothesis testing problem, to achieve a sample complexity (upper bound) matching the lower bounds in Theorems~\ref{th:lbsmallalpha} and~\ref{th:lowerbound2}, respectively. We then present tests that satisfy these conditions for a wide range of example hypothesis-testing problems. 

\begin{assumption}\label{assmp:E}
Assume there exists a $\calP$-e-process $(E_n)_{n\ge 1}$ that is asymptotically growth rate optimal for $\mathcal Q$, meaning that for every $Q \in \mathcal Q$, 
        \[ \liminf\limits_{n\rightarrow \infty} \lrp{\frac{1}{n} \log E_n} \ge \KLinf(Q,\mathcal P) \quad Q\text{-a.s.} \]
\end{assumption}

This is a rather weak assumption. Such an e-process exists whenever $\KLinf(\cdot,\mathcal P)$ is weakly lower semicontinuous at $Q$ for every $Q \in \mathcal Q$.
\cite{ram2026power} show that this latter condition holds very frequently, even when $\mathcal Q = \calP^c$; some non-parametric examples include $\calP$ being an f-divergence ball around an arbitrary $P_0$ (including KL, TV and Hellinger), many integral probability metric balls (like the kernel maximum mean discrepancy and Kolmogorov-Smirnov) and Wasserstein balls. We leave the details to that paper.

For $Q\in\mathcal Q$ with $\KLinf(Q,\calP) > 0$ and an e-process $(E_n)_{n\ge 1}$, define the random index
\[
N_\eta(Q) := \inf\lrset{n\ge1: \inf\limits_{k\ge n}\lrp{\frac{1}{k}\log E_k} \ge(1-\eta) \KLinf(Q,\calP) },\numberthis\label{eq:warmstartindex}
\]
noting that by convention $\inf\emptyset=+\infty$.

\begin{assumption}\label{ass:warm-start}  An e-process $(E_n)_{n\ge 1}$ that satisfies Assumption~\ref{assmp:E} is said to have \emph{integrable convergence} if $\Exp{Q}{N_\eta(Q)}<\infty$  for every $Q\in\mathcal Q$ with $0<\KLinf(Q,\calP)<\infty$ and every $\eta\in(0,1)$.
\end{assumption}

Intuitively, the integrable convergence assumption is only requiring that the average number of samples needed for the e-process to enter an asymptotic growth phase with rate at least $(1-\eta)\KLinf(Q,\calP)$ is finite for each fixed alternative $Q$.

\begin{theorem}\label{th:generalasupperboundcond}
Let $\calP,\mathcal Q$ be arbitrary subsets of probability measures. Let $\lrset{E_n}_{n\ge 1}$ be a $\calP$-e-process that satisfies Assumption~\ref{assmp:E}. Then, for every $\alpha \in (0,1)$, the stopping time $\tau_{\alpha}$ defined in~\eqref{eq:generalstoppingtime} is $\alpha$-correct for $\calP$. Further, the following statements hold for every $Q\in \mathcal Q$ with $0 < \KLinf(Q,\calP) < \infty$. First, the family $\{\tau_\alpha\}_{\alpha \in (0,1)}$ is a nested family of stopping times that is power one, such that for every $\eta \in (0,1)$,
\[ \tau_{\alpha} \le N_{\eta}(Q) + \left\lceil \frac{\log(1/\alpha)}{(1-\eta)\KLinf(Q,\calP)} \right\rceil,\quad Q\text{-a.s.} \]
Consequently, 
    \[ Q\lrs{\limsup\limits_{\alpha\rightarrow 0} ~\frac{\tau_{\alpha}}{\log({1}/{\alpha})} \le \frac{1}{\KLinf(Q,\mathcal P)}} = 1. \]
Second, if $(E_n)_{n\ge 1}$ also satisfies Assumption~\ref{ass:warm-start}, then it follows that 
\[ \limsup\limits_{\alpha\to 0}~ \frac{\Exp{Q}{\tau_\alpha}}{\log(1/\alpha)} \le \frac{1}{\KLinf(Q,\calP)}. \]
\end{theorem}

\begin{proof}
Recall, $\tau_{\alpha} = \min\{n : E_n \ge {1}/{\alpha}\}  = \min\{n : ~ n \lrp{\tfrac{1}{n}\log E_n} \ge \log({1}/{\alpha})\}.$ Clearly, for $\alpha_1 \le \alpha_2$, we have $1/\alpha_1 \ge 1/\alpha_2$. Since hitting the larger threshold $1/\alpha_1$ cannot occur before hitting the smaller threshold $1/\alpha_2$, we have $\tau_{\alpha_1} \ge \tau_{\alpha_2}$.  

Consider the event 
    $\calE := \{ \liminf\nolimits_{n\rightarrow \infty}~\lrp{\tfrac{1}{n} \log E_n} \ge \KLinf(Q,\mathcal P) \}$. From the given condition, we have $Q[\calE] = 1$. Further, since for all $Q\in \mathcal Q$, $\KLinf(Q,\calP) > 0$, we have that $E_n\to\infty$, $Q$-a.s., and hence, $\tau_\alpha$ is power one. 

Next, for every $\eta \in (0,1)$, define $N_\eta(Q)$
 as in~\eqref{eq:warmstartindex}. Clearly, it is independent of $\alpha$. Further, the following holds pathwise on $\calE$ for all $\alpha \in (0,1)$: 
\begin{align*}
    \tau_{\alpha} 
    &\le N_\eta(Q) +  \min\lrset{ n :~ {n(1-\eta)\KLinf(Q,\mathcal P)} \ge \log\frac{1}{\alpha}}\\
    &\le N_\eta(Q) +  \left\lceil \frac{\log({1}/{\alpha})}{(1-\eta)\KLinf(Q,\mathcal P)} \right\rceil.\numberthis\label{eq:pathwiseubst}
\end{align*}
Since $N_\eta(Q) < \infty$ on  $\mathcal E$, dividing both sides of the above pathwise inequality by $\log(1/\alpha)$ and taking limits, we get that on $\calE$, 
\[ \limsup\limits_{\alpha\rightarrow 0}~ \frac{\tau_{\alpha}}{\log({1}/{\alpha})} \le \frac{1}{(1-\eta)\KLinf(Q,\mathcal P)} \]
Since $Q[\calE] = 1$, 
\begin{align*}
    1 
    &= Q\lrs{\calE, \limsup\limits_{\alpha\rightarrow 0}~ \frac{\tau_{\alpha}}{\log({1}/{\alpha})} \le \frac{1}{(1-\eta)\KLinf(Q,\mathcal P)}}\\
    &\le Q\lrs{\limsup\limits_{\alpha\rightarrow 0}~ \frac{\tau_{\alpha}}{\log({1}/{\alpha})} \le \frac{1}{(1-\eta)\KLinf(Q,\mathcal P)}}.
\end{align*}
The above implies (by intersecting a countable sequence of probability one events obtained by picking $\eta=1/i$ for $i\ge 2$) that
\[ Q\lrs{ \limsup\limits_{\alpha\rightarrow 0}~\frac{\tau_{\alpha}}{\log({1}/{\alpha})} \le \frac{1}{\KLinf(Q,\mathcal P)} } = 1, \]
proving the first statement.

The second part of the theorem follows from taking expectation on both sides of~\eqref{eq:pathwiseubst}, rearranging the terms, and taking the limit as $\alpha\to 0$. 
\end{proof}

\begin{corollary}[Point-vs-point]\label{cor:pvp}
    For any $Q$, $P$ such that $0<\KL(Q,P) < \infty$, there exists an $\alpha$-correct, power-one stopping time $\tau_\alpha$ for testing $\mathcal P: P$ vs. $\mathcal Q: Q$ that satisfies 
    \[ \limsup\limits_{\alpha\rightarrow 0}~ \frac{\tau_\alpha}{\log\frac{1}{\alpha}} \le \frac{1}{\KL(Q,P)},\quad \text{Q-a.s.} \]
\end{corollary}
\begin{proof}
First, note that in the point-vs-point setting, $\KLinf(Q,\calP) = \KL(Q,P)$. Next, observe that the likelihood ratio e-process, defined as $E_n := \prod\nolimits_{i=1}^n\frac{dQ}{dP}(X_i)$, satisfies Assumption~\ref{assmp:E}, as shown below:
\[ \lim\limits_{n\rightarrow \infty} \lrp{\frac{1}{n}\log E_n} = \lim\limits_{n\rightarrow \infty} \lrp{\frac{1}{n} \sum\limits_{i=1}^n \log \frac{dQ}{dP}(X_i) } = \KL(Q,P), \quad Q\text{-a.s.}, \]
where the last equality follows from SLLN. Hence, the desired bound on $\tau_\alpha$ defined by~\eqref{eq:generalstoppingtime} follows from Theorem~\ref{th:generalasupperboundcond}.
\end{proof}

\paragraph*{Why an integrable convergence condition cannot be omitted} Theorem~\ref{th:generalasupperboundcond} shows that the integrable convergence condition (Assumption~\ref{ass:warm-start}) is sufficient for a $\calP$-e-process satisfying Assumption~\ref{assmp:E} to ensure that the associated stopping time has the correct first-order asymptotic expected sample size under every $Q \in \cal Q$. In the following example, we present an $\alpha$-correct test for a Gaussian mean-testing problem, obtained by thresholding a specific e-process which satisfies Assumption~\ref{assmp:E}, but violates the integrable convergence condition. We then show that the resulting stopping time is power one, but has an infinite expected sample size under every alternative distribution, demonstrating the necessity of the integrable convergence condition for obtaining a tight upper bound on $\Exp{Q}{\tau_\alpha}$.

\begin{example}
    Consider $\calP = \{N(0,1)\}$ and $\mathcal Q = \lrset{ N(m,1): m > 0 }$.     Let $(E^{\mathrm{base}}_n)_{n\ge 0}$ be the given base $\calP$-e-process satisfying Assumption~\ref{assmp:E}, and $D$ be $\mathcal F_1$-measurable (completely determined after observing $X_1$), $\mathbb N$-valued random variable. Here, $D$ represents a data-dependent delay. Let $E_0 = 1$, and for $n\ge 1$, define
    \[
    E_n:=
    \begin{cases}
    1, & 1\le n\le D,\\[2pt] E^{\mathrm{base}}_{\,n-D}(X_{D+1},\ldots,X_n), & n>D,
    \end{cases}
    \]
    where $E^{\mathrm{base}}_{k}(X_{D+1},\ldots,X_{D+k})$ denotes the base $\calP$-e-process evaluated on $X_{D+1},\ldots,X_{D+k}$. 
    Clearly, $(E_n)_{n\ge 1}$ is a $\calP$-e-process that also satisfies Assumption~\ref{assmp:E}. In addition, for every $Q \in \mathcal Q$ with $\KLinf(Q,\calP) > 0$ and $\eta \in (0,1)$,  the index $N_\eta(Q)$ for $E_n$ (defined in~\eqref{eq:warmstartindex}) satisfies $N_\eta(Q) \ge D+1$, $Q$-a.s. Further, for $\alpha\in (0,1)$, the stopping time $\tau_\alpha := \min\{n \ge 1: E_n \ge 1/\alpha \}$ satisfies $\tau_\alpha \ge D + 1$, since for $n\le D$, $E_n = 1$ and $1/\alpha > 1$. Finally, by setting $D:= \lceil{ \exp(X^2_1/2) \rceil}$ to have an infinite expectation under every $Q\in\mathcal Q$,  $(E_n)_{n\ge 1}$ serves as the required witness for the necessity of the integrable convergence condition.
    \end{example}

Theorem~\ref{th:generalasupperboundcond} provides a sufficient condition for a test to be optimal in $\alpha\to 0$ setting, matching the lower bounds in Theorem~\ref{th:lbsmallalpha}. Recall that Theorem~\ref{th:lowerbound2}  only applies to the composite alternative case; it is not relevant in the point-vs-point case discussed above (or the composite-vs-point case discussed later in Section~\ref{sec:ub_nonpcompnullparampointalt}). In the following theorem, we present sufficient conditions for a sequential test to achieve the optimal sample complexity bound in the setting of Theorem~\ref{th:lowerbound2}. These sufficient conditions correspond to fast concentration of empirical $\KLinf$ under the null and alternative, which are explained after the theorem, and which we prove hold in the non-parametric examples.

\begin{theorem}\label{generalupperboundcond2}
    Let $\calP$ and $\mathcal Q$ be arbitrary subsets of probability measures, and $\{Q_k\}_k$ be a sequence in $\mathcal Q$ such that $\KLinf(Q_k,\calP) \to 0$ as $k\to \infty$. For $n\in \N$, let ${F}_n$ denote the empirical distribution of $n$ i.i.d.\ observations. Suppose the following hold.
    \begin{enumerate}
    \item\label{cond1} There exists a constant  $C> 0$ such that for every $\alpha \in (0,1/2)$, 
    \[ \sup_{P\in \calP}~ P\lrs{\exists n \geq \log\tfrac{1}{\alpha} : n \KLinf({F}_n, \calP) \ge C(\log\tfrac{1}{\alpha} + \ln\ln(e+n))} \le \alpha.\]

    \item\label{cond2}There exist positive constants $k_0, C'$, and a function $H: M_1(\R) \rightarrow \R_+$, such that  $\forall n\ge 1$, and for every $\alpha\in (0,1)$,
        \begin{align*} \sup\limits_{k\ge k_0}~ &Q_k\left[\sqrt{\KLinf({F}_{n}, \calP)} - \sqrt{\KLinf(Q_k, \calP)}   \le - \sqrt{H(Q_k)}\sqrt{\frac{C' \log(1/\alpha) }{n}} \right] \le \alpha.  \end{align*}
    \end{enumerate}
    Then, for every fixed $\alpha\in (0,1/2)$,  there exists  $\tau_\alpha$ that is $\alpha$-correct under $\calP$. It has power one against $\mathcal Q$  if $\liminf_n \KLinf(F_n, \calP) > 0$, $Q$-a.s. for all $Q\in\mathcal Q$. Further, if $\limsup_{k\rightarrow \infty} H(Q_k) < \infty$, then 
    \[ \limsup\limits_{k\rightarrow \infty} ~\frac{\Exp{Q_k}{\tau_\alpha}}{\operatorname{KL^{-1}_{inf}}(Q_k,\calP) \log\log \operatorname{KL^{-1}_{inf}}(Q_k, \calP)} < \infty. \numberthis\label{eq:loglogub}\]
\end{theorem}
\begin{remark}\label{rem:parametericub2}
    The above theorem can handle non-parametric $\mathcal Q$, but if $\mathcal Q = \{Q_\theta\}_{\theta \in \Theta}$ is parametric for a parameter family $\Theta$, we can even replace the empirical distribution $F_n$ in the above theorem with $Q_{\theta_n}$, for some natural estimator $\theta_n$ like the MLE, and the theorem still holds unchanged. More generally, sometimes it may not be straightforward to prove a concentration result for $\KLinf(F_n,\calP)$. Instead, one may have similar guarantees for a different statistic (or process, $M_n$), whose asymptotic growth rate remains close to that for $\KLinf$. In Appendix~\ref{app:generalUB2} (Theorem~\ref{app_generalupperboundcond2}), we present a generalization of Theorem~\ref{generalupperboundcond2} that accounts for this. 

\end{remark}
To build intuition for the conditions and quantities in Theorem~\ref{generalupperboundcond2} above, first note that Condition~\ref{cond1} requires an $O(\log\log n)$ time-uniform concentration bound for the empirical $\KLinf$ process around its null-limit value, namely $0$, uniformly over all elements in the null class. Next, Condition~\ref{cond2} corresponds to a sub-Gaussian concentration of the square-root of the empirical $\KLinf$ around its limiting value under all but finitely many elements of the alternative sequence $\{Q_k\}$. These conditions are stronger than those in Theorem~\ref{th:generalasupperboundcond}, which suffice for achieving the lower bound in Theorem~\ref{th:lbsmallalpha}. In that setting, it is enough to construct a $\calP$-e-process whose limiting growth rate is $\KLinf(\cdot, \calP)$, without any additional concentration requirement. The proof of Theorem~\ref{generalupperboundcond2} is straightforward and is deferred to  Appendix~\ref{app:generalupperboundcond2proof}. 

To establish the tightness of the lower bounds in Theorems~\ref{th:lbsmallalpha} and~\ref{th:lowerbound2}, we now construct tests that satisfy the sufficient conditions of Theorems~\ref{th:generalasupperboundcond} and~\ref{generalupperboundcond2}, respectively, across a range of representative testing problems. We note that the universal inference e-process \citep{wasserman2020universal} often achieves the condition in Theorem~\ref{th:generalasupperboundcond} for quite general composite nulls and alternatives (see \citet{dixit2025anytime}), demonstrating the tightness of our lower bound in Theorem~\ref{th:lbsmallalpha} for an extremely wide range of problems. Since universal inference requires likelihoods (and usually also reference measures), it is mostly relevant for parametric problems, so our examples below focus on non-parametric settings.

\subsection{Non-parametric composite null, parametric point alternative} \label{sec:ub_nonpcompnullparampointalt}
Suppose $\mathcal P$ is the class of one-sided 1-sub-Gaussian distributions, and for the alternative, we have a Gaussian distribution with a positive mean. Formally, for $m > 0$,
\[ \mathcal P = \lrset{P: \Exp{P}{e^{\theta X - \frac{\theta^2}{2}}}\le 1, \quad \forall \theta\ge 0} \quad \text{and} \quad \mathcal Q = \{N(m,1)\}. \numberthis\label{eq:PQ1}  \]
One can check that the distributions in $\mathcal P$ have a nonpositive mean. \citet[\S 5.3]{larsson2024numeraire} argue that the numeraire e-variable for this problem is $E^*_Q(X) := e^{m X - m^2/2}$. Note that 
\[ \Exp{Q}{\log E^*_Q(X)} = \frac{m^2}{2} = \KL(Q,N(0,1)) \ge \KLinf(Q,\calP). \numberthis\label{eq:klinfsg} \] 
We note that~\eqref{eq:klinfsg}, in fact, holds as an equality. We will show this later. In the following corollary to Theorem~\ref{th:generalasupperboundcond}, we construct a $\calP$-e-process using $E^*_Q$, and show that it is optimal.
\begin{corollary}\label{cor:1}
    There exists an $\alpha$-correct, power-one stopping rule $\tau_\alpha$ for testing $\calP$ vs. $\mathcal Q$, as defined in~\eqref{eq:PQ1} such that for every $Q\in\mathcal Q$,
    \[ Q\lrs{\limsup\limits_{\alpha \rightarrow 0}~ \frac{\tau_\alpha}{\log({1}/{\alpha})} \le \frac{1}{\KLinf(Q,\mathcal P)} } = 1\quad\text{and}\quad \limsup\limits_{\alpha\to 0}~ \frac{\Exp{Q}{\tau_\alpha}}{\log(1/\alpha)} \le \frac{1}{\KLinf(Q,\calP)}. \]
\end{corollary}
\begin{proof}
The process $E_n:=E_{n-1}E^*_Q(X_n)$ with $E_0=1$ is a $\calP$-e-process that satisfies the following Q-a.s.:
\[ \liminf\limits_{n\rightarrow\infty}\lrp{\frac{1}{n} \log E_n} = \liminf\limits_{n\rightarrow\infty}\lrp{\frac{1}{n}\sum\limits_{i=1}^n \log E^*_Q(X_i)} \overset{(a)}{=} \Exp{Q}{\log E^*_Q(X)} \overset{(b)}{\ge} \KLinf(Q,\mathcal P),  \]
where $(a)$ follows from the SLLN and $(b)$ from~\eqref{eq:klinfsg}. Furthermore, in Lemma~\ref{lem:warmstart_cor:1}, we show that $E_n$ satisfies the integrable convergence condition (Assumption~\ref{ass:warm-start}). The two bounds in the corollary then follow for $\tau_\alpha$ defined in~\eqref{eq:generalstoppingtime} from Theorem~\ref{th:generalasupperboundcond}. 
\end{proof}
Recall that the lower bound of Theorem~\ref{th:lowerbound2} holds only when the alternative is composite, and in particular, when there exist several alternative distributions that are close to a common boundary null element. In the following section, we consider a composite vs. composite extension of the previous example and demonstrate the tightness of both lower bounds.

\subsection{Non-parametric composite null, parametric composite alternative} \label{sec:ub_nonpcompnullpcompalt}
In this example, our null is the collection of all 1-sub-Gaussian distributions. We consider a composite alternative consisting of all Gaussian distributions with a positive mean and unit variance. Formally, 
\[ \mathcal P = \lrset{P: \Exp{P}{e^{\theta X - \frac{\theta^2}{2}}}\le 1, \quad \forall \theta \in  \R} \quad\text{and}\quad \mathcal Q = \{ N(m,1): m > 0 \}. \numberthis\label{eq:PQ2} \]
Recall that $m_P = 0$ for all $P\in\cal P$. In particular, $N(0,1)\in \calP$. Hence, for any $Q\in\mathcal Q$ with mean $m_Q > 0$, we have $\KLinf(Q,\mathcal P) \le \KL(Q,N(0,1)) = m^2_Q/2$. In the following corollary of Theorem~\ref{th:generalasupperboundcond}, we present a test that is optimal for every (unknown) $Q\in\mathcal Q$, as $\alpha\rightarrow 0$. 
\begin{corollary}\label{cor:2}
    There exists an $\alpha$-correct, power-one stopping rule $\tau_\alpha$ for testing $\calP$ vs. $\mathcal Q$, as defined in~\eqref{eq:PQ2} such that for every $Q\in\mathcal Q$,
    \[ Q\lrs{\limsup\limits_{\alpha \rightarrow 0}~ \frac{\tau_\alpha}{\log({1}/{\alpha})} \le \frac{1}{\KLinf(Q,\mathcal P)} } = 1 \quad \text{and} \quad \limsup\limits_{\alpha\to 0}~ \frac{\Exp{Q}{\tau_\alpha}}{\log(1/\alpha)} \le \frac{1}{\KLinf(Q,\calP)}.\]
\end{corollary}
\begin{proof} For $\theta \in \R$, define
$E'_n(\theta) := \exp\{\theta \sum\nolimits_{i=1}^n X_i - \frac{n\theta^2}{2}\}$. 
From the condition on the elements of $\calP$, it is easy to verify that for every $\theta \in \R$, $E'_n(\theta)$ is a $\calP$-e-process. Observe that the e-process from Corollary~\ref{cor:1} corresponds to choosing $\theta = m$, the mean of the alternative distribution in~\eqref{eq:PQ1}. However, unlike in the point alternative setting of Section~\ref{sec:ub_nonpcompnullparampointalt}, here the test does not know $Q$. To get around this, we use a mixture of $E'_n(\cdot)$ instead, mixed with $N(0,1)$ prior on $\theta \in \R$, that is, 
\[ E_n := \frac{1}{\sqrt{2\pi}}\int\limits_{\R} E'_n(\theta) e^{-\frac{\theta^2}{2}} d\theta = \sqrt{\frac{1}{n+1}}~ e^{\frac{n^2 (\hat{m}_n)^2}{2(n+1)}}, \]
where $\hat{m}_n \!:= \!\frac{1}{n}\sum\nolimits_{i=1}^n X_i$. $E_n$ is a $\calP$-e-process that for every $Q\in\mathcal Q$ satisfies $Q$-a.s.:
\begin{align*}
    \liminf\limits_{n\rightarrow \infty}~\lrp{\frac{1}{n} \log E_n} = \liminf\limits_{n\rightarrow \infty}~\lrp{ \frac{1}{2n}\log\frac{1}{n+1} + \frac{n (\hat{m}_n)^2}{2(n+1)}} = \frac{m^2_Q}{2} \ge \KLinf(Q,\mathcal P).
\end{align*}
Furthermore, in Lemma~\ref{lem:warmstart_cor:2}, we show that $E_n$ satisfies the integrable convergence condition (Assumption~\ref{ass:warm-start}). The two bounds in the corollary then follow for $\tau_\alpha$ defined in~\eqref{eq:generalstoppingtime} from Theorem~\ref{th:generalasupperboundcond}. 
\end{proof}

Since the alternative in~\eqref{eq:PQ2} is composite with $\inf\nolimits_{Q\in\mathcal Q}\KLinf(Q,\calP) \le \inf\nolimits_{Q\in\mathcal Q} m^2_Q/2 = 0$, the lower bound in Theorem~\ref{th:lowerbound2} also holds for the testing of $\calP$ vs. $\mathcal Q$. In the following corollary of Theorem~\ref{generalupperboundcond2}, we present a sequential test that incurs an iterated logarithmic complexity along every sequence $Q_k \in\mathcal Q$ satisfying $\KLinf(Q_k,\calP) \rightarrow 0$ as $k\rightarrow \infty$, thus achieving the lower bound in Theorem~\ref{th:lowerbound2} along every direct sequence of alternative distributions.

\begin{corollary}
    There exists a power-one, $\alpha$-correct stopping rule $\tau_\alpha$ to test $\calP$ vs. $\mathcal Q$, as defined in~\eqref{eq:PQ2}, which satisfies the following for every sequence $Q_k\in\mathcal Q$ such that $\KLinf(Q_k,\calP) \rightarrow 0$ as $k\rightarrow \infty$: 
    \[ \limsup\limits_{k\rightarrow \infty}~\frac{\Exp{Q_k}{\tau_\alpha}}{\operatorname{KL^{-1}_{inf}}(Q_k,\mathcal P) \log\log\operatorname{KL^{-1}_{inf}}(Q_k,\mathcal P) } <  \infty.  \]
\end{corollary}
\begin{proof}
We first argue that for any $Q\in\mathcal Q$, $\KLinf(Q,\calP) = m^2_Q/2$. An upper bound follows trivially since $\KLinf(Q,\calP) \le \KL(Q,N(0,1)) = m^2_Q/2$. Thus, it only remains to show that $\KLinf(Q,\calP) \ge m^2_Q/2$. To this end, observe that for the one-sided sub-Gaussian null from the previous example (Section~\ref{sec:ub_nonpcompnullparampointalt}), denoted by ${\calP_1}$ in this section, \citet[\S 5.3]{larsson2024numeraire} argue that $\KLinf(Q,{\calP}_1) = m^2_Q/2$. Since $\calP \subset {\calP}_1$, we have $\KLinf(Q,\calP_1) \le \KLinf(Q,\calP)$, proving the lower bound. 

Now, fix an arbitrary sequence $\{Q_k\} \subset \mathcal Q$ such that $\KLinf(Q_k, \calP) \to 0$ as $k\to\infty$. Let $m_n := \sum_{i\le n} X_i /n$ denote the empirical average of $n$ i.i.d.\ samples. Since $\mathcal Q$ is parametric, where each distribution is uniquely identified with its mean, we will use the observation in  Remark~\ref{rem:parametericub2} and show that $\KLinf(N(m_n, 1),\calP)$, which equals $(\sum_i X_i)^2/2n^2$, satisfies the conditions of Theorem~\ref{generalupperboundcond2}, along with explicit constants $C$, $C'$, and $k_0$, as well as the function $H$. Then, Theorem~\ref{generalupperboundcond2} proves existence of a stopping time $\tau_\alpha$ that satisfies the bound in the corollary.

First, to see that Condition~\ref{cond1} holds, for $\theta \in \R$, define $E'_n(\theta) := \exp\{\theta \sum\nolimits_{i=1}^n X_i - \frac{n\theta^2}{2}\}$. From the definition of $\calP$, it is easy to verify that $E'_n(\theta)$ is a $\calP$-e-process. For $c \ge 6.6 e$, consider the heavy-near-zero mixture distribution supported on $\R$ \citep{robbins1968iterated} :
\[ \pi(\theta) := \frac{1}{Z_0} \frac{{\bf 1}_{|\theta| \le 1}}{|\theta| (\ln\frac{c}{|\theta|})(\ln\ln\frac{c}{|\theta|})^2}, \quad \text{ for } \theta \in \R,  \]
where $Z_0$ represents the normalization constant. $E_n := \int_{\R} E'_n(\theta) \pi(\theta) d\theta $ is also a $\calP$-e-process. 

In their Theorem 6 (with $S_n = \sum\nolimits_{i\le n} X_i$ and $V_n = n$), \citet{agrawal2025regret} bound the regret of the process $\log{E_n}$ with respect to $n \KLinf(N(m_n,1),\calP) = (\sum_{i\le n} X_i)^2/2n$. In particular, from Equation (8) in their paper, we see that Condition~\ref{cond1} holds with $C=4$. See also the discussion around Equation (9) in their work.

We now show that Condition~\ref{cond2} in Theorem~\ref{generalupperboundcond2} holds with $k_0 = 1$, $C' = 1$, and $H \equiv 1$. Observe that
\begin{align*}
    Q_k\!\lrs{ \sqrt{ \frac{(\sum_{i\le n} X_i)^2}{2n^2} } \!-\! \sqrt{\frac{m^2_{Q_k}}{2}} \le -\sqrt{\frac{\log(1/\alpha)}{n}} } \!\le\! Q_k\!\lrs{ \abs{ m_{Q_k} \!-\! \sum\limits_{i=1}^n \frac{X_i}{n} } \ge \sqrt{\frac{2\log(1/\alpha)}{n}} }.
\end{align*}
Now, using the Gaussian tail bound, the right-hand side above is at most $\alpha$.

Finally, observe that $\liminf_n \KLinf(N(m_n,1), \calP) = \liminf_n (\sum_i X_i)^2/2n^2 = m^2_Q/2$, $Q$-a.s., which is positive for each $Q\in\mathcal Q$, completing the proof of the corollary.
\end{proof}

\begin{remark}
     Note that the null considered in this section is smaller than the one-sided sub-Gaussian null discussed in Section~\ref{sec:ub_nonpcompnullparampointalt}. In Appendix~\ref{app:onesidedcompnullpcompalt}, we present tests that achieve the two lower bounds corresponding to this broader one-sided sub-Gaussian null and the composite alternative $\cal Q$ considered here. 
\end{remark}

\subsection{Non-parametric composite null, non-parametric composite alternative}\label{sec:testingmeanofbdd}
For $0<m_0<1$, in this example, we consider 
\[\mathcal P = \{P \in M_1[0,1]: m_P = m_0\} \quad\text{and}\quad \mathcal Q = \{Q \in M_1[0,1]: m_Q <  m_0\}. \numberthis\label{eq:PQ3} \] 
As in the previous section, we present sequential tests that demonstrate the tightness of the bounds in Theorems~\ref{th:lbsmallalpha} and~\ref{th:lowerbound2} for testing $\calP$ vs. $\mathcal Q$, as defined in~\eqref{eq:PQ3}. We begin by introducing some notation that will be used throughout this section.

For a probability measure $Q \in M_1[0,1]$ and $x\in [0,1]$, we write
\[ \KLinf(Q,x) := \inf\lrset{\KL(Q,P):~ P\in M_1[0,1], ~m_P = x}. \numberthis\label{eq:klinfbdd} \]
Observe that for $Q\in\mathcal Q$, $\KLinf(Q,m_0) = \KLinf(Q,\mathcal P)$. Moreover, $\KLinf(Q, m_Q) = 0$. Let $F_n$ denote the empirical distribution for $n$ observations $X_1, \dots, X_n$. With these, in the following corollary to Theorem~\ref{th:generalasupperboundcond}, we present a sequential test that is optimal for every unknown $Q\in\mathcal Q$ in the $\alpha\rightarrow 0$ limit. 

\begin{corollary}\label{cor:3}
    There exists an $\alpha$-correct, power-one stopping rule $\tau_\alpha$ for testing $\calP$ vs $\mathcal Q$, as defined in~\eqref{eq:PQ3}, such that for every $Q\in\mathcal Q$,
    \[ Q\lrs{\limsup\limits_{\alpha\rightarrow 0} \frac{\tau_\alpha}{\log({1}/{\alpha})} \le \frac{1}{\KLinf(Q,\mathcal P)}  } = 1 \quad\text{and}\quad \limsup\limits_{\alpha\to 0}~\frac{\Exp{Q}{\tau_\alpha}}{\log(1/\alpha)} \le \frac{1}{\KLinf(Q,\calP)}. \]
\end{corollary}
\begin{proof}
For an appropriate nonnegative $R_n$, define
\[ C_n := \lrset{x: n\KLinf(F_n,x) - R_n < \log\frac{1}{\alpha}}. \]
Consider a test that stops at time $\tau_\alpha$ given by
\begin{align*}
    \tau_\alpha 
    &= \min\lrset{ n : m_0 \notin C_n } = \min\lrset{n: n\KLinf(F_n, m_0) - R_n \ge \log\frac{1}{\alpha} }. \numberthis\label{eq:stoppingtime}
\end{align*}
Using the martingale construction in \citet[Lemma F.1]{agrawal2021optimal}, it follows that with $R_n = \log(1+n) + 1$, the process 
$$E_n := e^{n\KLinf(F_n, m_0) - R_n} = e^{n\KLinf(F_n, \mathcal P) - R_n}$$ 
is a $\calP$-e-process. Moreover, since $\KLinf(\cdot, m_0)$ is a continuous function \citep[Theorem 7]{honda2010asymptotically} and under $Q\in\mathcal Q$, $F_n \Rightarrow Q$ a.s., $E_n$ satisfies the growth condition (Assumption~\ref{assmp:E}) under every $Q\in\mathcal Q$. Furthermore, in Lemma~\ref{lem:warmstart_cor:3}, we show that $E_n$ also satisfies the integrable convergence condition (Assumption~\ref{ass:warm-start}). The two bounds in the corollary then follow from Theorem~\ref{th:generalasupperboundcond}. 
\end{proof}
Stopping times of the form~\eqref{eq:stoppingtime} have previously been used in asymptotically optimal pure exploration bandit algorithms (whose sample complexity exactly matches the corresponding lower bound in the $\alpha\rightarrow 0$ regime). \cite{agrawal2020optimal} were the first to propose such bandit algorithms in non-parametric settings, including that considered in this example. However, the test described in Corollary~\ref{cor:3} has an optimal dependence on $\alpha$ but a suboptimal dependence on $\KLinf(\cdot,m_0)$ (it is loose in terms of Theorem~\ref{th:lowerbound2}). In Corollary~\ref{th:loglogbddsupp} below, we demonstrate another sequential test that achieves the lower bound in Theorem~\ref{th:lowerbound2}. %


\begin{corollary}\label{th:loglogbddsupp} There exists an $\alpha$-correct, power-one stopping rule $\tau_\alpha$ to test $\calP$ vs. $\mathcal Q$, as defined in~\eqref{eq:PQ3}, such that for any sequence $Q_k \in \mathcal Q$, $Q_k \Rightarrow P \in \mathcal P$ with $\V[P] > 0$ as $k\rightarrow\infty$, we have
    \[ \limsup\limits_{k\rightarrow \infty} ~\frac{\Exp{Q_k}{\tau_\alpha}}{\operatorname{KL^{-1}_{inf}}(Q_k,\mathcal P)\log\log\operatorname{KL^{-1}_{inf}}(Q_k,\mathcal P)} <  \infty. \]
\end{corollary}
The proof for Corollary~\ref{th:loglogbddsupp} relies on some known and new properties of $\KLinf$, which we present in Appendix~\ref{app:kltildeinf}. We equip the space $M_1[0,1]$ with L\'evy metric, and prove the convergence of the optimizer $\lambda_{m_0}(\cdot)$ for $\KLinf(\cdot, m_0)$ along the specific sequence in $M_1[0,1]$ under consideration. We further prove an asymptotic expansion and a concentration result for $\KLinf$ in Proposition~\ref{prop:constconv} and Lemma~\ref{lem:tumultconcklinfDH}, respectively. We note that these results may be of independent interest. 
\begin{remark}
    In Corollary~\ref{th:loglogbddsupp}, we prove the stated dependence on $\KLinf(\cdot, \mathcal P)$ along every sequence $\{Q_k\}$ that converges to some $P \ne \delta_{m_0}$. Notice that this covers every direct sequence appearing in Theorem~\ref{th:lowerbound2} for the bounded-mean testing problem. Indeed, if ${\bf Q} = \{Q_k\}$ is a sequence directly approaching the null with anchor $P_{\bf Q} \in \calP$,  then there is a $C \ge 1$ such that ${\bf Q} \in \mathcal S_C$. Hence, $\KL(Q_k, P_{\bf Q}) \le C \KLinf(Q_k, \calP) \to 0$ as $k\to \infty$, and hence, $Q_k \Rightarrow P_{\bf Q}$. Further, if $P_{\bf Q} = \delta_{m_0}$, then  $\KL(Q_k, P_{\bf Q})=\infty$ for every $Q_k$, leading to a contradiction.
\end{remark}

We conclude this section with a proof for Corollary~\ref{th:loglogbddsupp}. All supporting results with proofs are presented in Appendix~\ref{app:kltildeinf}.

\begin{proof}[{Proof of Corollary~\ref{th:loglogbddsupp}}]
Fix an arbitrary sequence $\{Q_k\}\subset \mathcal Q$ satisfying the condition in the corollary statement. To prove the corollary, we show that the conditions in Theorem~\ref{generalupperboundcond2} are satisfied for $\KLinf(F_n,m_0)$ with constants $C>0$, $C'> 0$, and $k_0 \ge 1$, and a function $H$. 

First, $\liminf_n \KLinf(F_n, m_0) = \KLinf(Q,m_0)$, $Q$-a.s., which is strictly positive for each $Q\in\mathcal Q$. This follows from continuity of $\KLinf(\cdot, m_0)$. 

Next, to see that Condition~\ref{cond1} holds, 
first observe that for $\lambda \in I_{m_0} := [-1/m_0,1/(1-m_0)]$, $W_n(\lambda):= \prod_{i=1}^n (1-\lambda(X_i - m_0))$ is a $\calP$-e-process. Define the mixture process  
\[ E_n \!:=\! \int\limits_{I_{m_0}} \!W_n(\lambda) \pi(\lambda) d\lambda, ~ ~\text{ where }~~ \pi(\lambda) \!=\! \begin{cases}
    \frac{\ln\ln(6.6e)}{2\lambda \ln\lrp{\frac{6.6e}{(1-m_0) \lambda}}\lrs{ \ln\ln\lrp{ \frac{6.6e}{(1-m_0)\lambda} } }^2}, & 0 < \lambda \le \frac{1}{1-m_0},\\
    \frac{\ln\ln(6.6e)}{2|\lambda| \ln\lrp{\frac{6.6e}{m_0 |\lambda|}} \lrs{ \ln\ln \lrp{\frac{6.6e}{m_0 |\lambda|}} }^2 }, & -\frac{1}{m_0} \le \lambda < 0
\end{cases}  \]
is a modified Robbins' heavy-near-zero prior supported on $I_{m_0}$. In particular, one can verify that $\pi$ is a density over this interval. Then, 
\[ \KLinf(F_n, m_0) = \max\limits_{\lambda \in I_{m_0}}~ \lrp{\frac{1}{n}\log W_n(\lambda)}, \]
where the equality follows from the dual formulation for $\KLinf(\cdot, m_0)$, reviewed in Appendix~\ref{app:kltildeinf}. 

Define $R_n := n \KLinf(F_n,m_0) - \log E_n$. To get some intuition, $\log E_n$ can be interpreted as the $\log$ wealth of a mixture strategy and $n\KLinf(F_n, m_0)$ as the maximum $\log$ wealth in hindsight that can be obtained using any fixed strategy. Then, $R_n$ corresponds to the regret of the mixture strategy with respect to the best in hindsight. 

It is easy to achieve a worst-case regret of $O(\log n)$ for bounded observations, for example, using the uniform mixture from Corollary~\ref{cor:3}. In fact, one cannot have a smaller worst-case regret, meaning that there exists a sequence of bounded observations where that mixture will suffer $O(\log n)$ regret. However, following \cite{orabona2023tight}, recently, \citet[Theorem 4.1]{agrawal2026cover} established pathwise regret for the aforementioned Robbins' mixture wealth. They show that the regret is $O(\log\log n)$ on a set that contains most but not all sequences of bounded observations, under the null. On the complement set, it is linear. In particular, in Equation (7) of their work, the authors show that Condition~\ref{cond1} in Theorem~\ref{generalupperboundcond2} is satisfied for some $C \ge 3$. 

We now show that Condition~\ref{cond2} also holds. Since $\lim_k\KLinf(Q_k, \calP) = 0$, and $Q_k\Rightarrow P \in \calP$ with $\V[P] > 0$, there exists $k_0 \in \N$ such that for all $k \ge k_0$, the dual optimizer $\lambda_{m_0}(Q_k) \in [-1,1]$ (Lemma~\ref{lem:KLinflambdaconvg0}). In the rest of the proof, we fix such a $k \ge k_0$. Lemma~\ref{lem:tumultconcklinfDH} shows that $\KLinf(F_n,m_0)$ satisfies Condition~\ref{cond2} with $C'= 1$ and $H(Q_k) = D(Q_k,m_0)/C(Q_k,m_0)$, where $C(Q_k,m_0)$ and $D(Q_k,m_0)$ are constants defined in Lemma~\ref{lem:sqrtklinftildemultdev}. 

Finally, in Proposition~\ref{prop:constconv}, we show that
    \[ \lim\limits_{k\rightarrow \infty} ~\frac{D(Q_k,m_0)}{C(Q_k,m_0)} = \frac{1}{\V[P]} < \infty,\]
completing the proof. 
\end{proof}
We note that a higher-order analysis of $\KLinf$-based stopping rule used in Corollary~\ref{cor:3} was recently obtained by \citet{deep2026beyond}, who establish a central limit theorem for the stopping time in Wald's regime for the bounded-mean testing problem. Developing such a higher-order theory for general composite hypotheses, and in the Farrell regime, remains open.

\subsection{Hypotheses generated by finitely many constraints}\label{sec:ub_finitelygenhypo}
Finally, we demonstrate the tightness of the bound in Theorem~\ref{th:lbsmallalpha} for a much broader class of testing problems.
For $\calX \subset \R$ and a set of \emph{constraint functions} $\{\phi_1, \dots, \phi_K\}$, where for  all $i\in [K]$, $\phi_i: \calX \to \R$ are measurable, and $\phi_K \not\equiv 0$ on $\calX$, define
\begin{align}
    &\mathcal P \!=\! \lrset{\!P\!\in\! M_1(\calX)\!:\max\limits_{i\in[K]}~ \Exp{P}{|\phi_i(X)|} < \infty,  ~\max\limits_{i\in[K-1]} \Exp{P}{\phi_i(X)} \!\le\! 0 ~\text{ and }~\Exp{P}{\phi_K(X)} \!=\! 0}, \label{eq:generalnull}\\ 
    &\mathcal Q \!=\! \lrset{\!Q \!\in\!  M_1(\calX)\!:\max\limits_{i\in[K]} ~\Exp{Q}{|\phi_i(X)|} < \infty,    ~\max\limits_{i\in [K-1]}\Exp{Q}{\phi_i(X)} \!\le\! 0 ~\text{ and }~ \Exp{Q}{\phi_K(X)} \!<\! 0 }. \label{eq:generalalt}
\end{align}
By convention, when $K=1$, we omit the constraints 
$$\max_{i\in [K-1]} \Exp{P}{\phi_i(X)} \le 0 \quad \text{ and } \quad \max_{i\in[K-1]} \Exp{Q}{\phi_i(X)} \le 0$$
from $\calP$ and $\mathcal Q$, and in this case, the null and alternative become 
\begin{align*}
&\calP = \{P \in M_1(\calX): \Exp{P}{|\phi_1(X)|} < \infty, \text{ and } \Exp{P}{\phi_1(X)} = 0\} \\
&\mathcal Q = \{Q \in M_1(\calX): \Exp{Q}{|\phi_1(X)|} < \infty, \text{ and } \Exp{Q}{\phi_1(X)} < 0 \},
\end{align*}
respectively.  Clearly, the null and the alternative considered above are non-intersecting. This is a very general framework that encompasses many common testing problems. Some running examples to keep in mind that fit naturally into this framework are discussed below.  
\begin{itemize}
    \item {\bf Bounded mean testing.}  This setup was studied in detail in Section~\ref{sec:testingmeanofbdd}. Here, $\mathcal{X} = [0,1]$, $K = 1$, and $\phi_1(x) = x - m_0$. 
    \item {\bf Heavy-tailed mean testing.} Consider the problem of testing the mean of the distributions supported on $\R$ (that is, $\mathcal{X} = \R$) with a bounded second moment. For fixed and known constants $c > 0$ and $B > c^2$, this corresponds to the null and alternative given by $\{P \in  M_1(\R) : \Exp{P}{X^2} \le B,~ m_P = c\}$ and $\{Q \in  M_1(\R) : \Exp{Q}{X^2} \le B,~ m_Q < c\}$.
    \item {\bf Testing quantile for distributions supported on $\R$.} Here, $\calX = \R$ and $K= 1$. For $\beta \in (0,1)$ and $q \in \R$, define $\phi_1(x) :=  {\bf 1}_{(-\infty, q]}(x) - \beta$. Then, the null and alternative sets are given by $\calP = \lrset{  P \in  M_1(\R): \Exp{P}{\phi_1(X)} = 0}$ and $\mathcal Q = \lrset{Q \in  M_1(\R): \Exp{Q}{\phi_1(X)} < 0}$. 
\end{itemize}


\citet[Theorem 1]{clerico2024optimal} and \citet[Corollary 3.3]{larsson2025variables} show that every admissible e-variable for $\mathcal{P}$, $\mathcal P$-quasi surely equals 
\[ S^{\lambda,\kappa}(x) = 1 + \sum\limits_{i=1}^{K-1} \lambda_i \phi_i(x) + \kappa_1 \phi_K(x) + \kappa_2(- \phi_K(x)),\numberthis\label{eq:SandPi}\] 
where $(\lambda,\kappa) \in \Pi := \{(\lambda,\kappa) \in \R^{K+1}_+: ~ S^{\lambda,\kappa} \ge 0~~ \calP\text{-q.s.}\}$. 


In general, for $\mathcal{X} = \R$, without restrictions on the form of the constraint functions $\phi_i$, $\KLinf(Q, \mathcal{P})$ can be zero for all $Q \in  M_1(\R)$, implying that testing $\mathcal{P}$ against any alternative would require an unbounded number of samples on average. For example, this occurs if $K=1$ and $\phi_1(x) = x$ is linear. Given any $Q \in  M_1(\R)$, one can always construct another distribution with a prescribed mean that is arbitrarily close to $Q$ in $\KL$. We refer the reader to \citet[\S 2]{agrawal2020optimal} for a detailed discussion and proof.
In the aforementioned representation, note that if $\mathcal X = \R$, $K=1$, and $\phi_1(x) = x$ is a linear function, then $\kappa_1 = \kappa_2$~ to ensure the non-negativity of $S(x)$ 
and thus no nontrivial e-variables exist. This phenomenon is related to the previous observation that for this setting, $\KLinf(Q, \mathcal{P}) = 0$ for any $Q \in  M_1(\R)$. This connection will become apparent when we present the dual formulation of $\KLinf$. 

To exclude such degenerate cases and to facilitate the analysis below, we impose compactness of the null class and a point-charging condition that allows $\calP$-q.s. inequalities to be interpreted pointwise on $\calX$.

\begin{assumption}\label{assmp:compactnull}
    Assume that $\calP$ in~\eqref{eq:generalnull} is compact in the L\'evy metric and for each $x\in\calX$, there is a distribution $P_x \in \calP$ such that $P_x(\{x\}) > 0$. 
\end{assumption}
 These assumptions are readily verifiable in several important examples. Note that under Assumption~\ref{assmp:compactnull}, the $\calP$-q.s. inequality in the definition of $\Pi$ in~\eqref{eq:SandPi} reduces to a pointwise inequality over $\calX$. Following \cite{larsson2025variables}, we now introduce some terminology. For a vector $\rho \in \R^{K+1}_+$, we denote its support by $\supp(\rho) := \lrset{i : \rho_i > 0}$. An index set $I\subset [K+1]$ is called redundant if there exists a non-zero $\rho \in \R^{K+1}_+$ with $\supp(\rho) \subseteq I$ such that $\sum_{i=1}^{K-1} \rho_i \phi_i + (\rho_K-\rho_{K+1}) \phi_K  = 0$, $\calP$-quasi surely. Define
\[ \calK := \lrset{(\lambda, \kappa) \in \R^{K+1}_+: \supp(\lambda,\kappa) \text{ is not redundant}}. \numberthis\label{eq:setK}\]
Note that in the above definition, $\lambda \in \R^{K-1}_+$ and $\kappa \in \R^2_+$. 

We now turn to studying $\KLinf(\cdot,\calP)$ for the setting at hand. The dual representations for $\KLinf$ developed in specific settings by \cite{honda2010asymptotically,agrawal2020optimal,shekhar2026dual} can be extended to the more general setup of this section. However, in the following theorem, we use a different approach to arrive at an alternative representation for $\KLinf$.
\begin{theorem}\label{th:KLinfdual}
Under Assumption~\ref{assmp:compactnull}, for $\calP$ in~\eqref{eq:generalnull} and $Q\in M_1(\calX)$ such that $\KLinf(Q,\calP) < \infty$ and $\max\nolimits_{i\in [K]} \mathbb{E}_Q[|\phi_i(X)|] < \infty$, 
\begin{align*} 
    \KLinf(Q,\mathcal P) = \max\limits_{(\lambda,\kappa) \in \Pi \cap \calK}~\Exp{Q}{\log S^{\lambda,\kappa}(X)}.
\end{align*}
Further, if $(\lambda_\star, \kappa_\star) $ is any maximizer in $\Pi \cap \calK$, then $S^{\lambda_\star, \kappa_\star} > 0$, Q-a.s., and 
\[ \Exp{Q}{\frac{1}{S^{\lambda_\star, \kappa_\star}}} \le 1. \]
\end{theorem}
The alternative representation above connects $\KLinf$ to the finite-dimensional class of e-variables generated by the constraints. Compactness of $\calP$ (Assumption~\ref{assmp:compactnull}) ensures existence of the reverse $\KL$ projection of $Q$ onto $\calP$, which, together with the complete-class representation of e-variables for $\calP$ in \citet{larsson2025variables}, yields Theorem~\ref{th:KLinfdual}. We refer the reader to Appendix~\ref{app:finitelyconstraineddual} for a proof.


\begin{remark}
    The notion of separation between the unknown alternative $Q$ and the null set $\calP$ that naturally appears in the lower bounds and tests in this paper ($\KLinf(Q,\calP)$), is closely related, yet different from the infimum of $\KL$ divergences studied in \cite{larsson2024numeraire}. For brevity, we refer to the latter as $\infKL(Q,\calP)$. While in $\KLinf(Q,\calP)$, the infimum is over  $\calP$, in $\infKL(Q,\calP)$ it is over bipolar of $\calP$ (denoted as $\calP^{\circ\circ}$). Since $ \calP \subset \calP^{\circ\circ}$, we have $\infKL(Q,\calP) \le \KLinf(Q,\calP)$. While the exact relation between $\KLinf$ and $\infKL$ is not yet completely understood, Theorem~\ref{th:KLinfdual} shows that the two are the same for the rather general setup considered in this section (but they are known to be different in general).
   Further, the aforementioned paper shows that $\infKL$ is the correct quantity determining the hardness of testing using a single sample. However, our tight lower and upper bounds show that $\KLinf$ is the right quantity that captures the asymptotic difficulty of sequential testing.
\end{remark}


With the set of all admissible e-variables for $\calP$, and their relation to $\KLinf(\cdot,\calP)$ from Theorem~\ref{th:KLinfdual},  we are now ready to construct a $\calP$-e-process for testing $\calP$ vs. $\mathcal Q$, as defined in~\eqref{eq:generalnull} and~\eqref{eq:generalalt}, respectively. We present a stopping rule that is optimal as $\alpha\rightarrow 0$ in the following corollary to Theorem~\ref{th:generalasupperboundcond}. 
\begin{corollary}\label{cor:4}
    There exists an $\alpha$-correct, power-one stopping rule $\tau_\alpha$ for testing $\calP$ against $\mathcal Q$ from~\eqref{eq:generalnull} and~\eqref{eq:generalalt} such that for every $Q\in\mathcal Q$ with $\KLinf(Q,\calP) < \infty$, 
    \[ Q\lrs{ \limsup\limits_{\alpha\rightarrow 0}~\frac{\tau_\alpha}{\log({1}/{\alpha})} \le \frac{1}{\KLinf(Q,\mathcal P)}  } = 1 \quad\text{and}\quad \limsup\limits_{\alpha\to 0}~\frac{\Exp{Q}{\tau_\alpha}}{\log(1/\alpha)} \le \frac{1}{\KLinf(Q,\calP)}. \]
\end{corollary}
\begin{proof}
Let $\mathcal C := \operatorname{conv(\Pi\cap \calK)}$ denote the convex hull of $\Pi \cap \calK$, and $\mu_{\mathcal C}$ be a uniform probability measure on $\mathcal C$ relative to its affine hull. Note that a uniform prior is feasible since $\Pi \cap \calK \subseteq \R^{K+1}$ is compact \citep[Lemma 7.2]{larsson2025variables}, and hence, its convex hull, $\mathcal C$, is compact in finite dimensions. Moreover, $\Pi$ is convex, and hence, $\mathcal C \subseteq \Pi$. Define 
\[  E_n :=  \int\limits_{\mathcal C}\prod\limits_{i=1}^n  S^{\lambda,\kappa}(X_i) ~\mu_{\mathcal C}(d\lambda,d\kappa) = \int\limits_{\mathcal C}~e^{\sum\limits_{i=1}^n \log S^{\lambda,\kappa}(X_i)}~\mu_{\mathcal C} (d \lambda,d \kappa). \]
From~\eqref{eq:SandPi}, since for each fixed $(\lambda,\kappa) \in \Pi$, $S^{\lambda,\kappa}(X)$ is a $\calP$-e-variable, $\lrset{E_n}$, which is a mixture of $S^{\lambda,\kappa}(\cdot)$ for $(\lambda,\kappa) \in \mathcal C \subseteq \Pi$ is a $\calP$-e-process, and $\tau_\alpha$ from~\eqref{eq:generalstoppingtime} is an $\alpha$-correct stopping rule. 

We now relate the $\calP$-e-process $E_n$ to $\KLinf$. First, recall that the dimension of the affine hull of $\mathcal C$ is at most $K+1$, and that the functions $g_t(\lambda,\kappa) := \log S^{\lambda,\kappa}(X_t)$ are exp-concave in $(\lambda,\kappa)$. Then, \citet[Lemma F.1]{agrawal2021optimal}, along with the fact that $\Pi \cap \calK \subseteq \mathcal C$, yields that 
\[  \sup\limits_{(\lambda,\kappa) \in \Pi \cap \calK}  \sum\limits_{i=1}^n  \log S^{\lambda,\kappa}(X_i) \!\le\!  \sup\limits_{(\lambda,\kappa) \in \mathcal C} \sum\limits_{i=1}^n \log S^{\lambda,\kappa}(X_i) \!\le\! \log E_n \!+\! (K\!+\!1) \log(n+1) \!+\! 1. \numberthis\label{eq:seqineq}\]
Now, under $Q\in\cal Q$, let $F_n$ denote the empirical distribution of $n$ i.i.d.\ observations. We first  note that $\KLinf(F_n, \calP) < \infty$. This is because for each observation $X_i$, Assumption~\ref{assmp:compactnull} gives $P_{X_i} \in \calP$ such that $P_{X_i}(X_i) > 0$. Then, $\bar{P}_n := \tfrac{1}{n}\sum_i P_{X_i}$ also lies in  $\calP$ (convexity of $\calP$), and $F_n \ll \bar{P}_n$.  Thus, $\KLinf(F_n, \calP)$ exhibits the dual representation given in  Theorem~\ref{th:KLinfdual}.

Finally, dividing the left- and right-most sides of the sequence of inequalities in~\eqref{eq:seqineq} by $n$, and taking the limit as $n\rightarrow \infty$, we get
\begin{align*} \liminf\limits_{n\rightarrow \infty}~ \frac{1}{n} \log E_n 
&\ge \liminf\limits_{n\rightarrow \infty}~\lrp{\sup\limits_{(\lambda,\kappa) \in \Pi \cap \calK}~\frac{1}{n}\sum_i \log S^{\lambda,\kappa}(X_i)} \\
&= \liminf\limits_{n\rightarrow \infty}~\KLinf(F_n, \mathcal P) \tag{Theorem~\ref{th:KLinfdual}}\\
&\ge \KLinf(Q,\mathcal P),
\end{align*}
where the last inequality follows from lower-semicontinuity of $\KLinf(\cdot, \mathcal P)$ on $\cal Q$, which in turn follows from joint lower semicontinuity of $\KL(\cdot,\cdot)$ in the L\'evy metric, and compactness of $\cal P$ under the topology generated by L\'evy metric \citep[Theorem 2, pg. 116]{berge1997topological}. Thus, $E_n$ satisfies the growth condition in Assumption~\ref{assmp:E}. Further, since $\calP$ is compact and for every $Q\in\mathcal Q$,  $\KLinf(Q,\calP)$ is lower semicontinuous, $\KLinf(Q,\calP) > 0$ for every $Q\in \mathcal Q$. Finally, in Lemma~\ref{lem:warmstart_cor:4}, we show that it also satisfies the integrable convergence condition (Assumption~\ref{ass:warm-start}). The bounds in the corollary then follow from Theorem~\ref{th:generalasupperboundcond}, completing the proof.
\end{proof}

\section{Joint optimality}\label{sec:jointoptimality}
In Section~\ref{sec:lb}, we presented two lower bounds, which informally capture two distinct costs: the small-error regime ($\alpha\to 0$) and the small-gap regime ($\KLinf \to 0$). In Section~\ref{sec:ub}, we then provided sufficient conditions under which tests attain the corresponding lower bounds in each regime, while treating the remaining parameters as fixed. In this section, we show that, given tests that achieve these individual optimality guarantees, it is straightforward to construct a test that is simultaneously optimal in both regimes, that is, one that attains both the lower bounds on expected stopping times in their respective asymptotic limits.
\begin{proposition}
    Let $\calP$ and $\mathcal Q$ be arbitrary collection of probability measures that satisfy $$\inf_{Q\in\mathcal Q} \KLinf(Q, \calP) = 0.$$ For every $\beta \in (0,0.5)$, suppose $\tauone_\beta\in\mathcal T_\beta$ and $\tautwo_\beta \in \mathcal T_\beta$ achieve the small-error and small-gap optimality guarantees as in Theorem~\ref{th:generalasupperboundcond} and Theorem~\ref{generalupperboundcond2}, respectively. Then, for every $\alpha \in (0,1)$, there exists a test $\tau_\alpha \in \mathcal T_\alpha$, that simultaneously achieves both the expected optimality guarantees, that is, it matches the small-error and small-gap lower bounds in their respective asymptotic regimes.
\end{proposition}
\begin{proof}
    The idea is to run the two level-$\alpha/2$ tests in parallel and reject the null as soon as either test rejects. Accordingly, we define $\tau_\alpha := \tauone_{\alpha/2} \wedge \tautwo_{\alpha/2}$, and hence, 
    $$\{\tau_\alpha < \infty\} = \{\tauone_{\alpha/2} < \infty\}  \cup \{ \tautwo_{\alpha/2} < \infty\}.$$ 
    It then clearly follows that $\tau_\alpha \in \mathcal T_\alpha$ is $\alpha$-correct and power one. Further, since $$\tau_\alpha \le \tauone_{\alpha/2}\quad \text{ and } \quad \tau_\alpha \le \tautwo_{\alpha/2}, \quad \text{ pathwise},$$ 
    we also have 
    $$ \Exp{Q}{\tau_\alpha} \le \mathbb{E}_{Q}[\tauone_{\alpha/2}]\quad \text{ and } \quad \Exp{Q}{\tau_\alpha} \le \mathbb{E}_{Q}[\tautwo_{\alpha/2}], \qquad \forall Q\in \mathcal Q.$$ 
    Therefore, $\tau_\alpha$ inherits the expected optimality properties of $\tauone_{\alpha/2}$ and $\tautwo_{\alpha/2}$ in the small-error and small-gap regimes, respectively, completing the proof.
    \end{proof}
While the above proposition achieves joint optimality by combining stopping times, an alternative approach is to combine the underlying wealth processes themselves. In particular, one may form a convex combination of two wealth processes, each attaining a single-regime optimality guarantee, and threshold the resulting ``hedged'' wealth process. We refer the reader to \cite{agrawal2025regret,agrawal2026cover} for related discussions on hedged e-processes.

\section{Conclusion}\label{sec:conc}
We establish tight lower bounds on the sample complexity of $\alpha$-correct, power-one sequential tests under two distinct asymptotic regimes. In the Wald's regime, where $\alpha\to0$ for a fixed alternative, our lower bound holds without structural assumptions on the null or alternative classes. In the Farrell's regime, where $\alpha$ is fixed and $\KLinf\to0$, we establish an iterated-logarithm lower bound along sequences that have several alternatives approaching a common null boundary component; along every such sequence, all but an asymptotically sparse set of alternatives incur the Farrell-rate sample complexity. These results require no common dominating reference measure and substantially extend classical parametric lower bounds to composite and non-parametric settings. We further provide general sufficient conditions for matching upper bounds and verify them in several parametric and non-parametric examples. Future work could consider non-i.i.d.\ settings, or investigate higher-order notions of optimality in either regime.


\section*{Acknowledgements}
{The authors thank Hongjian Wang for helpful conversations. SA acknowledges the generous support from the Pratiksha Trust, Bangalore, through the Young Investigator Award, and the DST INSPIRE Faculty Grant IFA24-ENG-389. AR acknowledges funding from NSF grant 2310718 and a Sloan fellowship.}

\bibliography{BibTex}       

\begin{appendix}

\section{On $\KLinf$ for bounded distributions}\label{app:kltildeinf}
Recall that for $0 < m_0 < 1$ and $Q\in M_1[0,1]$, 
\begin{align*} 
    \KLinf(Q,m_0) &:= \inf\limits_{\substack{P' \in M_1[0,1] \\ m_{P'} = m_0}}~\KL(Q,P') \\
    &= \max\limits_{\lambda \in [-\frac{1}{m_0}, \frac{1}{1-m_0}]}~\Exp{Q}{\log(1-\lambda (X-m_0) )} \tag{\cite{honda2010asymptotically}},
\end{align*}
and let 
\[ \lambda_{m_0}(Q) := \argmax\limits_{\lambda\in [-\frac{1}{m_0}, \frac{1}{1-m_0}]}~\Exp{Q}{\log(1-\lambda(X-m_0))}. \]
(Pick any one element from the set of maximizers, if there are multiple.) In this Appendix, we state and prove certain properties and concentration results for $\KLinf(\cdot,m_0)$ and $\lambda_{m_0}(\cdot)$. 

\subsection{Properties} For fixed $m_0\in (0,1)$ and $Q\in M_1[0,1]$ such that $Q \ne \delta_{m_0}$, it follows from the definition and concavity of $\log(\cdot)$ that $\Exp{Q}{\log(1-\lambda(X-m_0))}$ is a strictly concave function. Hence, there is a unique maximizer $\lambda$ for $\KLinf(Q,m_0)$; denote it by $\lambda_{m_0}(Q)$. 

Continuity of $\KLinf(\cdot, m_0)$ is well known; see \citet[Theorem 7]{honda2010asymptotically} who state the result without giving a complete proof due to space constraints and instead note that it is somewhat complicated. We refer the reader to \citet[\S 4, Lemma 4.11]{thesis} for a complete and self-contained proof of the result for $\KLinf(\cdot,m_0)$ in a much more general setting. Note that the continuity is with respect to the weak topology, or equivalently, the topology induced by the L\'evy metric, in $ M_1[0,1]$. 
    
In the following lemma, we show the convergence of dual optimizers along a sequence of alternative distributions satisfying certain conditions. Recall, for a sequence of probability measures $\{Q_n\}$, and another probability measure $P$, we use $Q_n \Rightarrow P$ to mean that $Q_n$ converges weakly to $P$.
    \begin{lemma}\label{lem:KLinflambdaconvg0}
    For $m_0\in (0,1)$ and a sequence $Q_n \in  M_1[0,1]$ such that $Q_n \Rightarrow P \ne \delta_{m_0}$ with $m_{Q_n} \le m_0$ and $m_P = m_0$, we have $\lambda_{m_0}(Q_n) \rightarrow 0$ as $n \rightarrow \infty$. 
\end{lemma}

\begin{proof}
    Recall that $I_{m_0}:=
    [-\tfrac{1}{m_0},\tfrac{1}{1-m_0}]$. For this proof, let the sequence of optimizers be denoted by 
    $\lambda_n:=\lambda_{m_0}(Q_n)$. 
    By continuity of $\mathrm{KL}_{\inf}(\cdot,m_0)$ under weak
    convergence,
    $$
    \mathrm{KL}_{\inf}(Q_n,m_0)
    \longrightarrow
    \mathrm{KL}_{\inf}(P,m_0)=0,
    $$
    where the last equality follows since $P\in M_1[0,1]$ with $m_P=m_0$. 
    
    Next, consider an arbitrary convergent subsequence (existence is guaranteed since it is a sequence in a compact set), still denoted by $\lambda_n$, with $\lambda_n\to\bar\lambda\in I_{m_0}$. For $(\lambda,x)\in I_{m_0}\times[0,1]$, define
    $$
    \ell(\lambda,x)
    :=
    \log\bigl(1-\lambda(x-m_0)\bigr)
    $$
(with the convention $\log 0=-\infty$), and, for $M>0$, define
$$
\ell_M(\lambda,x)
:=
\log\left(
\bigl(1-\lambda(x-m_0)\bigr)\vee e^{-M}
\right).
$$
The function $\ell_M$ is bounded and continuous on the compact set
$I_{m_0}\times[0,1]$. Therefore, since $\lambda_n\to\bar\lambda$
and $Q_n\Rightarrow P$,
$$
\int \ell_M(\lambda_n,x)\,dQ_n(x)
\longrightarrow
\int \ell_M(\bar\lambda,x)\,dP(x).
$$
Since $\lambda_n$ maximizes the dual
representation of $\mathrm{KL}_{\inf}(Q_n,m_0)$, it follows that
$$
0
=
\lim_{n\to\infty}\mathrm{KL}_{\inf}(Q_n,m_0)
=
\lim_{n\to\infty}
\int \ell(\lambda_n,x)\,dQ_n(x)
\leq
\int \ell_M(\bar\lambda,x)\,dP(x),
$$
where in the last inequality, we used that $\ell(\lambda_n, \cdot) \le \ell_M(\lambda_n, \cdot)$ for each $\lambda_n$, and bounded convergence theorem. 
Letting $M\to\infty$, using
$\ell_M\to\ell$, and monotone convergence theorem,
$$
0
\leq
\int \ell(\bar\lambda,x)\,dP(x).
$$
Using the above and the inequality $\log u\leq u-1$ gives
$$
0 \le \int \ell(\bar\lambda,x)\,dP(x)
\leq
\int -\bar\lambda(x-m_0)\,dP(x)
=
-\bar\lambda(m_P-m_0)
=
0,
$$
since $m_P = m_0$. Thus equality holds, that is,
\[ \Exp{P}{\ell(\bar{\lambda},X)} = 0. \]
Since the equality in Jensen's inequality  is possible
only when $\bar{\lambda} = 0$ (since $P\ne \delta_{m_0}$), we get that every convergent subsequence of $\{\lambda_n\}$ has limit zero. Since $I_{m_0}$ is compact, we conclude that $\lambda_{m_0}(Q_n)\to0$.
\end{proof}

\begin{proposition}\label{prop:constconv}
Let \( m_0 \in (0,1) \), and let \( \{Q_k\}_{k \geq 1} \subset  M_1[0,1] \) be a sequence of probability measures with $m_{Q_k} < m_0$ that weakly converges to $P \in  M_1[0,1]\setminus\lrset{\delta_{m_0}}$ with $m_{P} = m_0$. Then, the following asymptotic expansion (in $k$) holds: 
    \[ \KLinf(Q_k, m_0) = \frac{1}{2} \frac{\lrp{\Exp{Q_k}{X-m_0}}^2}{\Exp{Q_k}{(X-m_0)^2}} + o\lrp{(\lambda_{Q_k})^2}, \] 
where \( \lambda_{Q_k} = \lambda_{m_0}(Q_k) \) is the unique dual maximizer for $\KLinf(Q_k,m_0)$. In addition, 
\[
\lim_{k \to \infty} \frac{\lrp{\log\lrp{\frac{1+\lambda_{Q_k} m_0}{1-\lambda_{Q_k}(1-m_0)}}}^2}{2\KLinf(Q_k,m_0)} = \frac{1}{\V[P]}.\] 
\end{proposition}
\begin{proof} First, observe that $\KLinf(Q_k, m_0) \rightarrow  0$ and $\lambda_{Q_k} \rightarrow 0$ as $k\rightarrow\infty$. These follow from continuity of $\KLinf(\cdot, m_0)$ in the topology of weak convergence and Lemma~\ref{lem:KLinflambdaconvg0}, respectively. Thus, there exists $k_0$ such that for all $k\ge k_0$, $\lambda_{Q_k} \in (-1, 1)$ and hence, satisfies the first order conditions for optimality. Going forward, we only consider $k\ge k_0$.  

To see the asymptotic expansion for $\KLinf(Q_k,m_0)$, define 
\begin{align*}
&f_k(\lambda) := \mathbb{E}_{Q_k}[\log(1 - \lambda(X - m_0))] \quad\text{and}\quad g_k(\lambda) := -f'_k(\lambda) = \mathbb{E}_{Q_k}\left[\frac{X - m_0}{1 - \lambda(X - m_0)}\right]. 
\end{align*}
Further, let
\[ a_k := \Exp{Q_k}{X-m_0}, \quad \text{and}\quad b_k := \Exp{Q_k}{(X-m_0)^2}. \]
From the dominated convergence theorem, we see that $g_k(\cdot)$ is infinitely differentiable for $\lambda\in (-1,1)$, so that \( g_k(0) = a_k \), and \( g_k'(0) = b_k \). 

Since $Q_k\Rightarrow P \ne \delta_{m_0}$, we have $b_k \to \operatorname{Var}[{P}] > 0$ and $\lambda_{Q_k} \to 0$. For $|\lambda| < 1$, all the relevant third derivatives are uniformly bounded, and hence, Taylor's expansion gives that uniformly in $k$, 
\begin{align*}
&f_k(\lambda) = -a_k  \lambda - 0.5 b_k  \lambda^2 + O(\lambda^3), \quad \text{ and }\quad g_k(\lambda) = a_k + \lambda b_k + O(\lambda^2).
\end{align*}
Further, the dual optimizer \( \lambda_{Q_k} \) satisfies $ f_k'(\lambda_{Q_k}) = -g_k(\lambda_{Q_k}) = 0.$ Solving for $\lambda_{Q_k}$ from $g_k(\lambda_{Q_k}) = 0$ gives
\[ a_k = -\lambda_{Q_k} b_k + O(\lambda^2_{Q_k}). \]
Using the form of $a_k$ from above, we get
\begin{align*}
f_k(\lambda_{Q_k})  = 0.5 b_k \lambda^2_{Q_k} + O(\lambda^3_{Q_k}).
\end{align*}
Also, 
\[ \frac{a^2_k}{b_k} = b_k \lambda^2_{Q_k} + O(\lambda^3_{Q_k}). \]
Since $\KLinf(Q_k,m_0) = f_k(\lambda_{Q_k})$, we have the claimed expansion in the proposition. 

To prove the limit in the lemma statement, we now consider the function in the numerator. Let \( \phi(\lambda) := \log(1 + \lambda m_0) - \log(1 - \lambda(1 - m_0)) \), so that on expanding $\phi$ around $0$, we get
\[ \phi(\lambda) = \lambda + \frac{\lambda^2}{2}\lrp{- m^2_0+(1-m_0)^2} + o(\lambda^2),  \]
which gives
\[
(\phi(\lambda))^2 = \lambda^2 + \lambda^3(-m^2_0 + (1-m_0)^2) + o(\lambda^3).
\]

Combining,
\begin{align*}
    \frac{\phi^2(\lambda_{Q_k})}{2\KLinf(Q_k,m_0)} 
    &= \frac{(\lambda_{Q_k})^2 + (\lambda_{Q_k})^3(-m^2_0 + (1-m_0)^2) + o((\lambda_{Q_k})^3)}{\Exp{Q_k}{(X-m_0)^2} (\lambda_{Q_k})^2 + o((\lambda_{Q_k})^2)}.
\end{align*}
Taking limits on both sides as $k\rightarrow \infty$, we get
\begin{align*}
    \lim\limits_{k\rightarrow\infty}~\frac{\phi^2(\lambda_{Q_k})}{2\KLinf(Q_k,m_0)} 
    &= \lim\limits_{k\rightarrow\infty}~\frac{1 + \lambda_{Q_k}(-m^2_0+ (1-m_0)^2) + o(\lambda_{Q_k})}{ \Exp{Q_k}{(X-m_0)^2}  + o(1)}\\
    &= \frac{1}{\V[P]},
\end{align*}
proving the limit. 
\end{proof}

\subsection{Concentration}

\begin{lemma}\label{lem:sqrtklinftildemultdev}
    For $Q \in  M_1[0,1]$, with $\lambda_{m_0}(Q) \in [-1,1]$, let $\hat{Q}_n$ denote the empirical distribution for $n$ i.i.d. samples from $Q$. For $1 > m_0 > m_Q > 0$, and $\epsilon >0$, 
    \[ Q\lrs{\sqrt{\KLinf(\hat{Q}_n, m_0)} \le \sqrt{\KLinf(Q,m_0)} - \epsilon} \le e^{-n\epsilon^2\frac{ (C(Q,m_0))}{D(Q,m_0)}}, \]
    where for the unique $\lambda_{m_0}\!(Q) \ge 0$ such that  $\KLinf(Q,m_0)\!=\! \Exp{Q}{\log(1-\lambda_{m_0}(Q)(X-m_0))}$,   
    \[ C(Q,m_0) := 2\KLinf(Q,m_0)\] 
    and 
    \[D(Q,m_0) :=  \lrp{\log\lrp{1+m_0\lambda_{m_0}(Q)}-\log\lrp{1-\lambda_{m_0}(Q)(1-m_0)}}^2.\]
\end{lemma}
\begin{proof} 
First, observe that the condition that $\lambda_{m_0}(Q) \in [-1,1]$ implies that $Q$ is not a point mass, which further implies that there is a unique dual maximizer $\lambda_{m_0}(Q)$. Without loss of generality, we consider $\epsilon < \sqrt{\KLinf(Q,m_0)}$, as otherwise, the bound holds trivially. Consider the following inequalities:
    \begin{align*}
        &Q\lrs{ \sqrt{\KLinf(\hat{Q}_n, m_0)} \le \sqrt{\KLinf(Q,m_0)} - \epsilon }\\
        &=Q\lrs{\KLinf(\hat{Q}_n, m_0) \le \KLinf(Q,m_0) + \epsilon^2 - 2\epsilon \sqrt{\KLinf(Q,m_0)} } \\
        &\overset{(a)}{\le} Q\lrs{\KLinf(\hat{Q}_n, m_0) \le \KLinf(Q,m_0) - \epsilon \sqrt{\KLinf(Q,m_0)} }, 
    \end{align*}
    where $(a)$ follows from the assumption on $\epsilon$. The inequality now follows from Lemma~\ref{lem:multdevKLinftilde}.
\end{proof}

\begin{lemma}\label{lem:multdevKLinftilde} For $Q \in  M_1[0,1]$ with $\lambda_{m_0}(Q) \in [-1,1]$, $1 > m_0 > m_Q > 0$, and $\epsilon > 0$, 
    \begin{align*}
        \log Q&\lrs{ \KLinf(\hat{Q}_n, m_0) \le \KLinf(Q,m_0) - \epsilon \sqrt{\KLinf(Q,m_0)} } \\
        &\qquad\qquad \le \frac{-2n\epsilon^2 \KLinf(Q,m_0) }{(\log\lrp{1+m_0\lambda_{m_0}(Q)}- \log\lrp{1-\lambda_{m_0}(Q)(1-m_0)})^2},
    \end{align*}
where $\lambda_{m_0}(Q) \ge 0$ uniquely satisfies $\KLinf(Q,m_0) = \Exp{Q}{\log(1-\lambda_{m_0}(Q)(X-m_0))}$, and $\hat{Q}_n$ denotes the empirical distribution with $n$ i.i.d. samples from $Q$.
\end{lemma}

\begin{proof}
    Recall that for measure $Q$ in consideration, 
    $$\KLinf(Q,m_0) = \max\limits_{\lambda\in \lrs{-1, 1}}~ \Exp{Q}{\log\lrp{1-\lambda(X - m_0)}},$$
    and $\lambda_{m_0}(Q)$ denotes the unique optimizer in the above expression (since $Q$ is not a point mass, as otherwise, $\lambda_{m_0}(Q)$ would be outside $[-1,1]$, violating the assumption). Further, $\lambda_{m_0}(Q) \in [-1,1] \subset (-\frac{1}{m_0}, \frac{1}{1-m_0})$. For simplicity of notation, we will drop its dependence on $Q$ and $m_0$ from the notation, and instead call it $\lambda^*$ in the rest of this proof. 

    Next, without loss of generality, we assume $\epsilon < \sqrt{\KLinf(Q,m_0)}$ as otherwise, the bound trivially holds. Clearly, 
    \begin{align*}
        & Q\lrs{ \KLinf(\hat{Q}_n, m_0) \le \KLinf(Q,m_0) - \epsilon \sqrt{\KLinf(Q,m_0)} }\\
        &\qquad
        \le Q\lrs{ \sum\limits_{i=1}^n \log\lrp{1- \lambda^*(X_i - m_0)} \le n\KLinf(Q,m_0) - n\epsilon \sqrt{\KLinf(Q,m_0)} }, \numberthis\label{eq:KLinf_multform}
    \end{align*}
    where the inequality follows since we replaced the maximum on the left with a feasible (and possibly sub-optimal) choice of the variable being optimized. Now, let 
    \[ Y(x) := \log(1-\lambda^*(x - m_0)). \]
    Since 
    $$\lambda^* \in [-1,1] \subset \left(-\frac{1}{m_0}, \frac{1}{1-m_0}\right),$$ 
    we have  $(1-\lambda^*(x-m_0)) > 0$ for all $x\in [0,1]$. Moreover, since $m_0 > m_Q$, we also have $\lambda^* > 0$. Thus,    
    $$Y(x) \in [ \log(1-\lambda^*(1-m_0)), ~\log(1+m_0\lambda^*)], \quad \text{ for all } x\in[0,1].$$ 
    With this, observe that $Y(X_i)$ are i.i.d., and bounded in  $ [ \log(1-\lambda^*(1-m_0)), ~\log(1+m_0\lambda^*)]$. Moreover, $\Exp{X\sim Q}{Y(X)} = \KLinf(Q,m_0)$. Hence, using Hoeffding's inequality to further bound the right-hand side in~\eqref{eq:KLinf_multform}, we get 
    \begin{align*}
        \log Q&\lrs{ \KLinf(\hat{Q}_n, m_0) \le \KLinf(Q,m_0) - \epsilon\sqrt{\KLinf(Q,m_0)} } \\
        &\qquad\qquad \le  \frac{-2n\epsilon^2 \KLinf(Q,m_0)}{\lrp{\log\lrp{1+m_0\lambda^*}-\log\lrp{1-\lambda^*(1-m_0)}}^2},
    \end{align*}
    proving the desired bound. 
\end{proof}

\begin{lemma}\label{lem:tumultconcklinfDH}
    Let $m_0 \in (0,1)$, $Q\in  M_1[0,1]$ be such that $1 > m_0 > m_Q > 0$, $\lambda_{m_0}(Q) \in [-1,1]$, and $\hat{Q}_n$ denote the empirical distribution with $n$ i.i.d. samples from $Q$. Then, for $C(Q,m_0)$ and $D(Q,m_0)$ as in Lemma~\ref{lem:sqrtklinftildemultdev}, for all $n\in \N$ and $\alpha \in (0,1)$, we have 
    \[ Q\lrs{\sqrt{\KLinf(\hat{Q}_n,m_0)} \le \sqrt{\KLinf(Q,m_0)} - \sqrt{\frac{\log(1/\alpha)}{n}} \sqrt{\frac{D(Q,m_0)}{C(Q,m_0)}}  } \le \alpha.\]
\end{lemma}
\begin{proof}
The inequality in the lemma follows by setting the right-hand side in the bound of Lemma~\ref{lem:sqrtklinftildemultdev} to at most $\alpha$. 
\end{proof}

\begin{lemma}\label{lem:subexpconcKLinf} 
    Let $m_0 \in (0,1)$, $Q \in M_1[0,1]$ be such that $1 > m_0 > m_Q > 0$, and $F_n$ denote the empirical distribution with $n$ i.i.d.\ samples from $Q$. Then for $\epsilon > 0$,
    \[ Q\lrs{ \KLinf(F_n, m_0) - \KLinf(Q,m_0) \le -\epsilon } \le e^{-n\epsilon^2/\lrp{4 \lrp{\epsilon + e^{-\KLinf(Q,m_0)}\frac{1-m_Q}{1-m_0} + e^{\KLinf(Q,m_0)}}}}.\]
    Further, for $\epsilon \in (0,e^{\KLinf(Q,m_0)}+e^{-\KLinf(Q,m_0)}\tfrac{1-m_Q}{1-m_0})$,
    \[ Q\lrs{ \KLinf(F_n, m_0) - \KLinf(Q,m_0) \le -\epsilon } \le e^{-\frac{n\epsilon^2}{8} \lrp{\frac{e^{-\KLinf(Q,m_0)}}{1+e^{-2\KLinf(Q,m_0)}\frac{1-m_Q}{1-m_0}}}}. \]
\end{lemma}
\begin{proof}
    Consider the following inequalities:
    \begin{align*}
        &Q\lrs{ \KLinf(F_n, m_0) - \KLinf(Q, m_0) \le -\epsilon }\\
        &= Q\lrs{ \max\limits_{\lambda\in I_{m_0}} \frac{1}{n} \sum\limits_{i=1}^n \log(1-\lambda(X_i - m_0)) - \Exp{Q}{\log(1-\lambda_{m_0}(Q)(X-m_0))} \le -\epsilon  }\\
        &\le Q\lrs{  \frac{1}{n} \sum\limits_{i=1}^n \log(1-\lambda_{m_0}(Q)(X_i - m_0)) - \Exp{Q}{\log(1-\lambda_{m_0}(Q)(X-m_0))} \le -\epsilon  }.
    \end{align*}
    Let 
    \[ Y_i := \log(1-\lambda_{m_0}(Q)(X_i - m_0)) - \Exp{Q}{\log(1-\lambda_{m_0}(Q)(X_i - m_0))}. \]
    Clearly, $\Exp{Q}{Y_i} = 0$. Further, 
    \[ \Exp{Q}{e^{Y_i}}  = e^{-\KLinf(Q,m_0)}\lrp{1+\lambda_{m_0}(Q)(m_0 - m_Q)} \le \frac{1-m_Q}{1-m_0}e^{-\KLinf(Q,m_0)}, \]
    and
    \[\Exp{Q}{e^{-Y_i}} = e^{\KLinf(Q,m_0)} \Exp{Q}{\frac{1}{1-\lambda_{m_0}(Q)(X - m_0)}} \le e^{\KLinf(Q,m_0)}, \]
    where the last inequality follows from optimality of $\lambda_{m_0}(Q)$. Using Lemma~\ref{lem:bernstein}, for $\epsilon > 0$, we get the first bound. The inequality in the other case follows from the above inequality by using the condition on $\epsilon$. 
\end{proof}

\section{Technical details for Section~\ref{sec:lb}}\label{app:proofs_seclb}

Our lower bound proofs use a standard change-of-measure inequality, corresponding to a data processing inequality for KL divergence at stopping times. We begin by recalling this result from \citet[Lemma 1]{kaufmann2016complexity}.

\begin{lemma}\label{lem:DP}
    Consider any two probability measures $P$ and $Q$ such that $Q\ll P$. Let $\tau$ denote a stopping time such that $\Exp{Q}{\tau} < \infty$, and let $\mathcal E$ be an $\mathcal F_\tau$-measurable event. Then,
    \[ \Exp{Q}{\tau} \KL(Q, P) \ge  \kl(Q[\calE], P[\calE]), \]
    where for $p\in [0,1]$ and $q\in [0,1]$, $\kl(p,q)$ denotes the $\KL$ divergence between Bernoulli distributions with means $p$ and $q$, respectively. 
\end{lemma}
We will apply this inequality with carefully chosen events $\calE$ on the right hand side and the measure $P$ on the left hand side, to derive a lower bound for $\Exp{Q}{\tau}$.

\subsection{Proof for Theorem~\ref{th:lbsmallalpha}}\label{app:th:lbsmallalpha_proof} In this section, we prove the lower bound in Theorem~\ref{th:lbsmallalpha}. The proof of~\eqref{eq:lbexpectedineq} generalizes that of \cite{kaufmann2016complexity}, who prove lower bounds for the sample complexity of the best-arm identification problem in the stochastic multi-armed bandit setting with parametric arm distributions. A similar but asymptotic lower bound was proved in \citet[\S 1.4]{garivier2019explore} in a great generality, but for the regret-minimization problem in the stochastic multiarmed bandit setting. 

\begin{proof} We first note that if there does not exist any $P \in \mathcal P$ such that $\KL(Q,P) < \infty$, we have $\KLinf(Q,\mathcal P) = \infty$, and the three bounds in~\eqref{eq:lbexpectedineq},~\eqref{eq:aslbongrid}, and~\eqref{eq:aslbineq} hold immediately. So we consider $\KLinf(Q,\mathcal P) < \infty$ without loss of generality in the rest of the proof.

Next, the bound in~\eqref{eq:lbexpectedineq} immediately follows whenever $\Exp{Q}{\tau_\alpha} = \infty$. So, let us consider the case when this does not hold. Since $\KLinf(Q,\calP) < \infty$, there exists $P \in \mathcal P$ such $0 < \KL(Q,P) < \infty$, which implies that the Radon-Nikodym derivative ${dQ}/{dP}$ exists $Q$-a.s. For $n\in\N$, let $X_1, X_2,  \dots $ be i.i.d. samples from $Q \in \mathcal Q$,  $\mathcal F_n := \sigma\lrp{X_1, \dots, X_n}$, and $$\mathcal L^{(n)}_{Q,P}:= \sum\nolimits_{i=1}^n \log \frac{dQ}{dP}(X_i),$$
denote the log-likelihood ratio of observing $n$ i.i.d. samples from $Q$ and that under $P$. Since $\Exp{Q}{\tau_\alpha} < \infty$ and $\KL(Q,P) < \infty$, using Lemma~\ref{lem:DP}, we have
\[ \Exp{Q}{\mathcal L^{\tau_\alpha}_{Q,P}} = \Exp{Q}{\tau_\alpha} \KL(Q,P) \ge \kl(Q(\calE), P(\calE)), \quad \text{ for all } \calE \in \mathcal F_{\tau_\alpha}. \]
Now, choose $\mathcal E = \lrset{ \tau_\alpha < \infty }$ so that $Q(\calE) = 1$ and $P(\calE) \le \alpha$, and we get 
\[ \Exp{Q}{\tau_\alpha}\KL(Q,P) \ge \kl(1,\alpha) = \log(1/\alpha). \]
We get~\eqref{eq:lbexpectedineq} by optimizing over $P\in\cal P$.

We now prove~\eqref{eq:aslbongrid}. Fix a rational $c \in (0,1/\KLinf(Q,\calP))$. Let $\epsilon > 0$ and $\delta > 0$ be small enough to ensure $c (\KLinf(Q,\calP) + \epsilon + \delta) < 1$. Since $\KLinf(Q,\calP) < \infty$, there exists $P_\epsilon \in \calP$ such that 
\[ \KL(Q,P_\epsilon) \le \KLinf(Q,\calP) + \epsilon < \infty. \]
Clearly, $Q \ll P_\epsilon$, and hence, the Radon-Nikodym derivative ${dQ}/{dP_\epsilon}$ exists ($Q$-a.s.), and we define the log-likelihood ratio process 
\[  \mathcal L^{(n)}_{Q,\epsilon} := \sum\limits_{i=1}^n \log\frac{dQ}{dP_\epsilon}(X_i), \qquad \forall n \ge 1, \]
with $\mathcal L^{(0)}_{Q,\epsilon} = 0$. 
Consider an arbitrary $\alpha$-correct power-one sequential test $\tau_\alpha$, and let $f_{nc} = \floor{nc}$. We will show that the event $A_{n}(c)$, defined below, cannot happen infinitely often under any $Q\in\mathcal Q$: 
$$A_n(c) := \{ \tau_{e^{-n}} \le f_{nc} \} = \{ \tau_{e^{-n}} \le nc \},$$ 
where the last equality follows since $\tau_{e^{-n}} \in \N$. Next, define 
$$G_n := \{ \mathcal L^{(f_{nc})}_{Q,\epsilon} \le f_{nc} (\KL(Q,P_\epsilon) + \delta) \}.$$ 
Then, $A_n(c) \in \mathcal F_{f_{nc}}$ and $G_n \in \mathcal F_{f_{nc}}$, and hence, 
\begin{align*}
    e^{-n} \ge P_\epsilon\lrs{ \tau_{e^{-n}} < \infty } 
    &\ge 
    P_\epsilon\lrs{ A_n(c), ~ G_n }\\
    &\ge \Exp{Q}{ e^{-\mathcal L^{(f_{nc})}_{Q,\epsilon}} {\bf 1}\lrset{A_n(c), ~  G_n } } \\
    &\ge e^{-{f_{nc}} \lrp{\KL(Q,P_\epsilon) + \delta}}  Q\lrs{ A_n(c), ~  G_n   } \\
    & \ge e^{-{nc} \lrp{\KLinf(Q,\calP) + \epsilon + \delta}}  Q\lrs{ A_n(c), ~  G_n }.
\end{align*}
On rearranging the above inequality,  we get
\[ Q\lrs{A_n(c), ~ G_n} \le e^{-n \lrp{ 1 - c(\KLinf(Q, \calP) +  \epsilon + \delta)}}. \]
Since the exponent on the right hand side above is strictly negative, we get
\[ \sum\limits_{n=1}^\infty Q[ A_n(c), ~ G_n ] \le \sum\limits_{n=1}^\infty e^{-n \lrp{1-c\lrp{\KLinf(Q,\calP) + \epsilon + \delta}}} < \infty,\]
which implies $Q[ \{A_n(c), ~ G_n \} \text{  i.o.} ] = 0$ \citep[\S 2.1]{williams1991probability}. 

Next, since $\mathcal L^{(f_{nc})}_{Q,\epsilon}/f_{nc} \to \KL(Q,P_\epsilon)$ almost surely as $n\to \infty$ (SLLN), $Q[ G^c_n \text{ i.o.}] = 0$. Together, these imply $Q[A_n(c) \text{ i.o.}] = 0$, which gives
\[ Q\lrs{\liminf\limits_{n\to\infty}\frac{\tau_{e^{-n}}}{n} < c} = 0 \quad  \equiv \quad Q\lrs{ \liminf\limits_{n\to\infty} \frac{\tau_{e^{-n}}}{n} \ge c  } = 1. \]
Since the above is true for every rational $c\in (0,1/\KLinf(Q,\calP))$, intersecting the event over all rational numbers in this interval gives~\eqref{eq:aslbongrid}.

We now use~\eqref{eq:aslbongrid} to prove~\eqref{eq:aslbineq}. Fix $\omega\in\Omega_0$, where 
$$\Omega_0 := \{ \liminf\nolimits_{n} \tau_{e^{-n}}/n \ge 1/\KLinf(Q,\calP) \}.$$ 
From~\eqref{eq:aslbongrid}, $Q[\Omega_0] = 1$. Next, fix a rational $c\in(0,1/\KLinf(Q,\calP))$. There exists integer $N(\omega,c)$ such that $\forall n\ge N(\omega,c)$, $\tau_{e^{-n}}(\omega)>cn$. Let $\alpha\in(0,e^{-N(\omega,c)}]$, and define $n:=\lfloor\log(1/\alpha)\rfloor$. Then,  
$$n\ge N(\omega,c)\quad\text{ and }\quad n\le \log(1/\alpha)<n+1.$$ 
Equivalently, $e^{-(n+1)}<\alpha\le e^{-n}$. Now, since $\alpha\le e^{-n}$, by the nesting assumption, $\tau_{\alpha}(\omega)\ge \tau_{e^{-n}}(\omega)$. Therefore, $\tau_{\alpha}(\omega)\ge \tau_{e^{-n}}(\omega)>cn$. Also since $\log(1/\alpha)<n+1$, we get that
\[
\frac{\tau_{\alpha}(\omega)}{\log(1/\alpha)}>\frac{cn}{n+1}.
\]
Taking limit as $\alpha\to 0$, we have $n\to\infty$ and $\frac{cn}{n+1}\to c$, and hence,
\[
\liminf_{\alpha\to 0}\frac{\tau_{\alpha}(\omega)}{\log(1/\alpha)}\ge c.
\]
Since the above bound holds for all rational $c < 1/\KLinf(Q,\calP)$, we also have
\[
\liminf_{\alpha\to 0}\frac{\tau_{\alpha}(\omega)}{\log(1/\alpha)}\ge \frac{1}{\KLinf(Q,\calP)}.
\]
The bound in~\eqref{eq:aslbineq} now follows by recalling that $\omega\in\Omega_0$, and $Q[\Omega_0] = 1$.
\end{proof}

\subsection{Proof of Corollary~\ref{cor:lilseq}}\label{app:proof_cor:lilseq}
Fix $\tau_\alpha\in\mathcal T_\alpha$. By assumption, there exists $C\ge 1$ such that $\mathcal S_C\neq\emptyset$. Fix a direct sequence $\mathbf Q=(Q_1,Q_2,\ldots)\in\mathcal S_C$. 
Applying Theorem~\ref{th:lowerbound2} with $\gamma=1$, it follows for the fixed pair $(\tau_\alpha,\mathbf Q)$,
\[
0 \le \limsup_{N\to\infty}
\frac{\nu\left(\tau_\alpha,\frac{1}{4C},N,\mathbf Q\right)}{N}
\le 
\limsup_{N\to\infty}
\sup_{\tau'_\alpha\in\mathcal T_\alpha}
\sup_{\mathbf Q'\in\mathcal S_C}
\frac{\nu\left(\tau'_\alpha,\frac{1}{4C},N,\mathbf Q'\right)}{N}
=0.
\]
By the definition of $\nu$, and the above sequence of inequalities, we have that the set
\[
A:=\lrset{ i\ge1: \Exp{Q_i}{\tau_\alpha} < \frac{1}{4C} F\left(\sqrt{\KLinf(Q_i,\calP)}\right)}
\]
has asymptotic density zero. Hence, its complement is infinite. Consequently, there exists a subsequence $i_1<i_2<\cdots$ such that $\KLinf(Q_{i_j},\mathcal P)\longrightarrow 0$ as $j\to\infty$, and 
\[
\frac{\Exp{Q_{i_j}}{\tau_\alpha}}
{F\left(\sqrt{\KLinf(Q_{i_j},\calP)}\right)}
\ge \frac{1}{4C},
\qquad  \forall j\ge1.
\]
From the definition of the limsup, it now follows that 
\[
\limsup_{\substack{Q\in\mathcal Q \\
\KLinf(Q,\mathcal P)\to0}}
\frac{\Exp{Q}{\tau_\alpha}}
{F\left(\sqrt{\KLinf(Q,\mathcal P)}\right)}
\ge
\limsup_{j\to\infty}
\frac{\Exp{Q_{i_j}}{\tau_\alpha}}
{F\left(\sqrt{\KLinf(Q_{i_j},\mathcal P)}\right)}
\ge \frac{1}{4C}>0,
\]
completing the proof.

\subsection{Proof of Lemma~\ref{lem:changeofmeasure}}\label{app:proof_changeofmeasure}
    If $\Exp{Q}{\tau} = \infty$ or $\KL(Q,P) = \infty$, then the stated bound is immediate. Hence, we suppose that both these quantities are finite. In particular, $Q\ll P$. Using Lemma~\ref{lem:DP}, we have
    \[ \Exp{Q}{\tau}\KL(Q,P) \ge \kl\lrp{ {Q}[\calE], {P}[\calE] } \overset{(a)}{\ge} \kl\lrp{\frac{1}{2},{P}[\calE]} = \frac{1}{2}\log\frac{1}{4{P}[\calE](1-{P}[\calE])}.\]
    The inequality $(a)$ follows since for $p\ge q$, $\kl(p,q)$ increases monotonically in $p$ and decreases in $q$. Rearranging the above inequality, we get the desired bound on ${P}[\calE]$.

\subsection{Proof of Lemma~\ref{lem:highprobevent}}\label{app:proof_highprobevent}
    Since $\tau_\alpha$ is the stopping time of an $\alpha$-correct, power-one sequential test, we have $Q[\calE_U] = 1$. We now show that 
    \[ Q\lrs{\tau_\alpha < \frac{d_l}{\KLinf(Q,\mathcal P)}} < 0.5 - \epsilon. \numberthis\label{eq:intineq1}  \]
    We will assume otherwise and derive a contradiction. To this end, consider an auxiliary algorithm $\A'$ that runs the sequential test for initial ${d_l}/{\KLinf(Q,\mathcal P)}$ steps. If the test halts by this time, then $\A'$ also halts and rejects the null. Otherwise, $\A'$ stops at time ${d_l}/{\KLinf(Q,\mathcal P)}$, and accepts the null. In the following, we denote by $Q_{\A'}[\cdot]$ the probability of an event, when sampling from $Q$ using algorithm $\A'$. Similarly, we use $\Exp{Q,\A'}{\cdot}$ to denote the expectation when sampling from $Q$ using the algorithm $\A'$. 
    
    Define the event $\calE_V := \lrset{ \A' \text{ rejects null} }$. Clearly, 
    \[ Q_{\A'}[\calE_V] = Q\lrs{\tau_\alpha \le \frac{d_l}{\KLinf(Q,\mathcal P)} } \ge 0.5 - \epsilon  \tag{by assumption}. \]

    On the other hand, for $P \in \mathcal P$, 
    \[P_{\A'}\lrs{ \calE_V } = P\left[\tau_\alpha < \frac{d_l}{\KLinf(Q,\mathcal P)}\right] \le \alpha < 0.5 - \epsilon \tag{$\tau_\alpha$ is $\alpha$-correct}\]
    Hence, by Lemma~\ref{lem:DP} with the event $\calE = \calE_V$, 
    \[ \Exp{Q,\A'}{\tau_\alpha \wedge \frac{d_l}{\KLinf(Q,\mathcal P)} } \KL(Q, P) \ge \kl\lrp{ Q_{\A'}[\calE_V] , P_{\A'}[\calE_V]} \ge \kl(0.5 - \epsilon,\alpha), \]
    where the last inequality follows since for $x,y\in (0,1)$ such that $x\ge y$, $\kl(x,y)$ is monotonically increasing in $x$ and decreasing in $y$. Finally, since the above inequality is true for all $P\in\mathcal P$, taking infimum over $P\in\calP$ on both the sides, we get
    \[ \Exp{Q,\A'}{\tau_\alpha \wedge \frac{d_l}{\KLinf(Q,\mathcal P)} } \KLinf(Q,\mathcal P) \ge \kl(0.5 - \epsilon, \alpha). \]
    Since the left hand side above is at most $d_l$, and hence, we get $d_l \ge \kl(0.5 - \epsilon, \alpha)$, which is a contradiction to the condition on the constant $d_l$ in the lemma statement. Hence,~\eqref{eq:intineq1} holds. 
    
    Finally, consider the following inequalities
    \begin{align*}
        Q\lrs{\calE(\Delta_{Q, \mathcal P})}
        &= Q\lrs{ \calE_U, \calE(\Delta_{Q,\mathcal P}) } \\
        &\overset{(a)}{\ge} Q\lrs{\calE_U} - Q\lrs{ \calE^c(\Delta_{Q,\mathcal P}) } \\
        &\overset{(b)}{\ge} Q\lrs{\calE_U} - Q\lrs{ \tau_\alpha \le \frac{d_l}{\KLinf(Q,\mathcal P)} } - Q\lrs{ \tau_\alpha \ge \frac{1}{\epsilon} \Exp{Q}{\tau_\alpha} } \\
        &\overset{(c)}{\ge} 1-0.5 + \epsilon - \epsilon \\ 
        &= 0.5,
    \end{align*}
    where $(a)$ follows since $P(A\cap B) \ge P(A) - P(B^c)$, $(b)$ from an application of union bound, and $(c)$ using~\eqref{eq:intineq1}, Markov's inequality, and the fact that $\tau_\alpha$ is a power-one test.

\section{Technical details for Section~\ref{sec:ub}}
\subsection{Proof for Theorem~\ref{generalupperboundcond2}} \label{app:generalupperboundcond2proof}
For $\alpha \in (0,1/2)$, let  
\[ \tau^{(2)}_\alpha := \min\lrset{n\ge\log
(1/\alpha): n\KLinf(F_n,\calP)  \ge C\lrp{ \ln\ln(e+n)+ \log({1}/{\alpha}) } }. \]
$\alpha$-correctness of $\tau^{(2)}_\alpha$ follows from Condition~\ref{cond1}. Further, if $\liminf_n \KLinf(F_n, \calP) > 0$ $Q$-a.s. for every $Q\in\mathcal Q$, then $\tau^{(2)}_\alpha$ has power one against $\mathcal Q$. 

Next, let $k_0$ be as in Condition~\ref{cond2}, and consider $Q_k$ with $k\ge k_0$. We now establish a high probability bound on $\tau^{(2)}_\alpha$. We will then use Theorem~\ref{th:hptoexp} to prove the existence of a desired stopping rule.  To this end, consider  the following inequalities:  
\begin{align*}
    \tau^{(2)}_\alpha 
    &= \min\lrset{ n \ge \log(1/\alpha): \sqrt{n \KLinf(F_n,\calP)} \ge \sqrt{ C\lrp{\ln \ln (e+n) + \log({1}/{\alpha})} } }\\
    &\le \min\lrset{ 2^\ell \ge \log(1/\alpha):~ \sqrt{2^\ell \KLinf(F_{2^\ell},\calP)} \ge \sqrt{ C\lrp{\ln \ln (e+2^\ell) + \log({1}/{\alpha})}} }.\numberthis\label{eq:stbound1}
\end{align*}
Next, for constant $C'$ and function $H$ from Condition~\ref{cond2}, consider the set 
\begin{align*} 
\mathcal G_{\alpha} &:= \Bigg\{\forall \ell\ge 1,~\sqrt{\KLinf(F_{2^\ell},\calP)} - \sqrt{\KLinf(Q_k,\calP)} \\
& \qquad \qquad \qquad \qquad \ge  - \sqrt{\frac{H(Q_k)}{2^\ell}} \sqrt{C'\lrp{\log(1/\alpha) + \log(\pi^2/6) + \log \ell^2}} \Bigg\}. 
\end{align*}
Using Condition~\ref{cond2} and Lemma~\ref{lem:DH}, we have $Q_k[\mathcal G_{\alpha}] \ge 1- \alpha$. Using this in~\eqref{eq:stbound1}, on $\mathcal G_\alpha$,
\begin{align*}
    \tau^{(2)}_\alpha 
    &\le \min\left\{ 2^\ell\ge \log(1/\alpha):~ 2^\frac{\ell}{2}\sqrt{\KLinf(Q_k,\calP)} \ge \sqrt{C\lrp{\ln \ln (e+2^\ell) + \log({1}/{\alpha})} } \right. \\
    &\qquad\qquad\qquad\qquad\qquad\qquad \left.+ \sqrt{H(Q_k)} \sqrt{C'\lrp{\log\tfrac{1}{\alpha} + \log(\pi^2\ell^2/6)}}  \right\}\\
    &\le \min\left\{ 2^\ell \ge \log(1/\alpha):~ 2^\ell \KLinf(Q_k,\calP) \ge  2C\lrp{ \ln \ln (e+2^\ell) + \log({1}/{\alpha}) }\right. \\
    &\qquad\qquad\qquad\qquad\qquad\qquad \left. +2C'{H(Q_k)} \lrp{\log\tfrac{1}{\alpha}+\log(\pi^2\ell^2/6)} \right\}\\    
    &\le \min\left\{ 2^\ell \ge \log(1/\alpha):~ 2^\ell \KLinf(Q_k,\calP) \ge 2C\lrp{ \ln \ln (e+2^\ell) + \log({1}/{\alpha}) + 1 }  \right. \\
    &\qquad\qquad\qquad\qquad\qquad\qquad \left. + 4C'{H(Q_k)} \lrp{\log\tfrac{1}{\alpha}+ 1 + \log\log(e+2^\ell)} \right\}\\    
    &\le \min\left\{ 2^\ell\ge\log(1/\alpha):~ 2^\ell  \ge \frac{4\lrp{2C'H(Q_k) \vee C}}{\KLinf(Q_k,\calP)} \lrp{\log\tfrac{1}{\alpha} + 1 + \log\log(e+2^\ell) }\right\},
\end{align*}
Lemma~\ref{lem:boundontaualpha} (with $\gamma = 0.5$) then implies that pathwise on $\mathcal G_\alpha$, the following holds: 
\begin{align*}
\tau^{(2)}_\alpha & \le 2 + \frac{32\lrp{2C'H(Q_k) \vee C}}{ \KLinf(Q_k,\calP)} \log\log \lrp{e+\frac{4C'H(Q_k) \vee C}{\KLinf(Q_k,\calP)}} \\
&\qquad\qquad\qquad + \frac{8\lrp{2C'H(Q_k) \vee C} }{\KLinf(Q_k,\calP)} (1+\log(1/\alpha)).
\end{align*}
Let right hand side of the above bound be $T(\alpha, Q_k)$. Then, for fixed $Q_k$ and $\alpha' < \alpha$, we have
\[ \frac{T(\alpha', Q_k)}{\log\tfrac{1}{\alpha'}} \le \frac{T(\alpha, Q_k)}{\log\tfrac{1}{\alpha}}. \]
Theorem~\ref{th:hptoexp} now gives $\tau_\alpha$ (using a meta algorithm that uses $\tau^{(2)}_\alpha$ as a sub routine) such that  for all $Q_k$ with $k\ge k_0$, 
\begin{align*} 
\Exp{Q_k}{\tau_\alpha} &\!\le\!  \frac{4(1-\alpha)}{(1-2\alpha)^2}\!\left( 2 + \frac{32\lrp{2C'H(Q_k) \vee C}}{\KLinf(Q_k,\calP)} \log\log \lrp{e+\frac{4C'H(Q_k) \vee C}{\KLinf(Q_k,\calP)}}\right. \\
&\qquad\qquad\qquad + \left. \frac{8\lrp{2C'H(Q_k) \vee C} }{ \KLinf(Q_k,\calP)} (1+\log(1/\alpha))\right).\numberthis \label{eq:sctau2alpha}
\end{align*}
Finally, since $\limsup_{k\to\infty} H(Q_k) < \infty$, dividing the above bound by $$\KL^{\text{-1}}_{\operatorname{inf}}(Q_k, \calP) \log\log\KL^{\text{-1}}_{\operatorname{inf}}(Q_k, \calP),$$ and taking limit as $k\rightarrow \infty$, we get~\eqref{eq:loglogub}, proving the theorem.

\subsection{General upper bound theorem}\label{app:generalUB2}
In this section, we present a generalization of Theorem~\ref{generalupperboundcond2} with its proof. 
\begin{theorem}\label{app_generalupperboundcond2}
    Let $\calP$ and $\mathcal Q$ be arbitrary subsets of probability measures, and $\{Q_k\}_k$ be a sequence in $\mathcal Q$ such that $\KLinf(Q_k,\calP) \to 0$ as $k\to \infty$. Let $\lrset{M_n}_{n}$ be non-negative stochastic processes such that the following hold:
    \begin{enumerate}
    \item\label{app_cond1} There exists a constant $C> 0$ such that for every $\alpha \in (0,1/2)$, 
    \[ \sup_{P\in \calP}~ P\lrs{\exists n \ge \log\tfrac{1}{\alpha} : n M_n \ge C\lrp{\log\tfrac{1}{\alpha} + \ln\ln(e+n)} } \le \alpha.\]

    \item\label{app_cond2} There exist positive constants $k_0$, $C'$, and functions $G: M_1(\R) \rightarrow \R_+$ and $H: M_1(\R) \rightarrow \R_+$ such that for all $n\ge 1$, and every $\alpha \in (0,1)$, 
    \begin{align*} \sup\limits_{k\ge k_0}~ &Q_k\left[\sqrt{M_{n}} - \sqrt{G(Q_k)} \le - \sqrt{H(Q_k)}\sqrt{\frac{C'\log(1/\alpha)}{n}} \right] \le \alpha.  \end{align*}
    \end{enumerate}
    Then, for every fixed $\alpha\in (0,1/2)$, there exists $\tau_\alpha$ that is $\alpha$-correct under $\calP$. It has power one against $\mathcal Q$  if  $\liminf_n M_n > 0$, $Q$-a.s. for every $Q \in \mathcal Q$. Additionally, if $G$ and $H$ satisfy the following:
    \[ \limsup\limits_{k\rightarrow \infty} H(Q_k) < \infty \qquad \text{and} \qquad \limsup\limits_{k\rightarrow \infty}~ \frac{G^{\operatorname{-1}}(Q_k)}{\operatorname{KL^{-1}_{inf}}(Q_k,\calP)}   < \infty, \numberthis\label{eq:app_cond3} \]
    then 
    \[ \limsup\limits_{k\rightarrow \infty} ~\frac{\Exp{Q_k}{\tau_\alpha}}{\operatorname{KL^{-1}_{inf}}(Q_k,\calP) \log\log \operatorname{KL^{-1}_{inf}}(Q_k, \calP)} < \infty. \numberthis\label{eq:app_loglogub}\]
\end{theorem}

\begin{proof}
Let 
\[ \tau^{(2)}_\alpha := \min\lrset{n \ge \log\tfrac{1}{\alpha}: n M_n \ge C\lrp{\log\tfrac{1}{\alpha} + \ln \ln (e+n)} }. \]
Observe that $\alpha$-correctness of $\tau^{(2)}_\alpha$ follows from Condition~\ref{app_cond1} in Theorem~\ref{app_generalupperboundcond2}. Next, let $k_0$ be as in Condition~\ref{app_cond2}, and fix $Q_k$ with $k\ge k_0$. We now establish a high-probability bound on $\tau^{(2)}_\alpha$. We will then use Theorem~\ref{th:hptoexp} to conclude existence of the desired $\tau_\alpha$ (which will be a meta algorithm using $\tau^{(2)}_\alpha$, defined above, as a sub-routine). To this end, consider the following inequalities: 
\begin{align*}
    \tau^{(2)}_\alpha 
    &= \min\lrset{ n \ge \log({1}/{\alpha}): \sqrt{n M_n} \ge \sqrt{ C\lrp{\ln \ln (e+n) + \log({1}/{\alpha})} } }\\
    &\le \min\lrset{ 2^\ell \ge \log({1}/{\alpha}):~ \sqrt{2^\ell M_{2^\ell}} \ge \sqrt{ C\lrp{\ln \ln (e+2^\ell) + \log({1}/{\alpha})}} }.\numberthis\label{eq:app_stbound1}
\end{align*}
Next, for constant $C'$ and  functions $G$ and $H$ from Condition~\ref{cond2}, consider the set 
\begin{align*} 
\mathcal G_{\alpha} &:= \Bigg\{\forall \ell\ge 1,~\sqrt{M_{2^\ell}} - \sqrt{G(Q_k)} \\
& \qquad \qquad \qquad \qquad \ge  - \sqrt{\frac{H(Q_k)}{2^\ell}} \sqrt{C'\lrp{\log(1/\alpha) + \log(\pi^2/6) + \log\ell^2}} \Bigg\}. 
\end{align*}
Using Condition~\ref{app_cond2} and Lemma~\ref{lem:DH}, we have $Q_k[\mathcal G_{\alpha}] \ge 1- \alpha$. Using this in~\eqref{eq:app_stbound1}, on $\mathcal G_{\alpha}$,
\begin{align*}
    \tau^{(2)}_\alpha 
    &\le \min\left\{ 2^\ell \ge \log(1/\alpha):~ 2^\frac{\ell}{2}\sqrt{G(Q_k)} \ge  \sqrt{C\lrp{\ln \ln (e+2^\ell) + \log({1}/{\alpha})} } \right. \\
    &\qquad\qquad\qquad\qquad\qquad\qquad \left. + \sqrt{H(Q_k)} \sqrt{C'\lrp{\log\tfrac{1}{\alpha} + \log\lrp{\pi^2\ell^2/6}}} \right\}\\
    &\le \min\left\{ 2^\ell \ge \log(1/\alpha):~ 2^\ell G(Q_k) \ge 2C\lrp{ \ln \ln (e+2^\ell) + \log({1}/{\alpha}) } \right. \\
    &\qquad\qquad\qquad\qquad\qquad\qquad \left. +2C'{H(Q_k)} \lrp{ \log\tfrac{1}{\alpha} + \log(\pi^2\ell^2/6)}  \right\}\\   
    &\le \min\left\{ 2^\ell \ge \log(1/\alpha):~ 2^\ell G(Q_k) \ge 2C\lrp{\ln \ln (e+2^\ell) + \log({1}/{\alpha}) + 1} \right. \\
    &\qquad\qquad\qquad\qquad\qquad\qquad \left. + 4C'{H(Q_k)} \lrp{ \log\tfrac{1}{\alpha} + \log\log(e+2^\ell) + 1}  \right\}\\    
    &\le \min\left\{ 2^\ell\ge \log({1}/{\alpha}):~ 2^\ell  \ge \frac{4\lrp{2C'H(Q_k) \!\vee\! C}}{G(Q_k)}\lrp{\log\tfrac{1}{\alpha} + 1+ \log\log(e+2^\ell) }\right\}.
\end{align*}
Lemma~\ref{lem:boundontaualpha} (with $\gamma = 0.5$) then implies that pathwise on $\mathcal G_\alpha$, the following holds: 
\begin{align*}
\tau^{(2)}_\alpha & \le  1 + \frac{32\lrp{2C'H(Q_k) \vee C}}{G(Q_k)} \log\log \lrp{e+\frac{4C'H(Q_k) \vee C}{G(Q_k)}} \\
&\qquad\qquad\qquad + \frac{8\lrp{2C' H(Q_k) \vee C}}{G(Q_k)} \lrp{1+\log(1/\alpha)}.
\end{align*}
Let the right-hand side of the above bound be $T(\alpha, Q_k)$. Then, for fixed $Q_k$ and $\alpha' < \alpha$, we have ${T(\alpha', Q_k)}/{\log\tfrac{1}{\alpha'}} \le {T(\alpha, Q_k)}/{\log\tfrac{1}{\alpha}}$. Using Theorem~\ref{th:hptoexp}, there exists $\tau_\alpha$ such that for all $Q_k$ with $k\ge k_0$,
\begin{align*} 
\Exp{Q_k}{\tau_\alpha} &\!\le\!  \frac{4(1-\alpha)}{(1-2\alpha)^2}\!\left( 1 +  \frac{8\lrp{2C'H(Q_k) \vee C}}{G(Q_k)}\lrp{1+\log(1/\alpha)}\right. \\
&\left. + \frac{32\lrp{2C'H(Q_k) \vee C}}{G(Q_k)} \log\log\lrp{ e+ \frac{4C' H(Q_k) \vee C}{G(Q_k)}} \right).\numberthis \label{eq:app_sctau2alpha}
\end{align*}
Finally, from the conditions in~\eqref{eq:app_cond3}, 
\[\limsup\limits_{k\rightarrow\infty}~ \frac{G^{\text{-1}}(Q_k) \log\log (e + G^{\text{-1}}(Q_k))}{ \operatorname{KL^{-1}_{inf}}(Q_k, \calP) \log\log \operatorname{KL^{-1}_{inf}}(Q_k, \calP)}  < \infty.  \numberthis\label{eq:app_lnlnconvcond}\]
Dividing the bound obtained in~\eqref{eq:app_sctau2alpha} for $Q_k$ by $\KL^{\text{-1}}_{\operatorname{inf}}(Q_k, \calP) \log\log \KL^{\text{-1}}_{\operatorname{inf}}(Q_k, \calP)$, taking limits as $k\rightarrow \infty$, and using the conditions in~\eqref{eq:app_cond3} and~\eqref{eq:app_lnlnconvcond}, we get~\eqref{eq:app_loglogub}, completing the proof. 
\end{proof}

\subsection{Integrable convergence condition holds on examples in Section~\ref{sec:ub}}
In this section, we show that several e-processes studied in Section~\ref{sec:ub} satisfy the integrable convergence condition (Assumption~\ref{ass:warm-start}). Recall that to establish the integrable convergence condition, we need to prove that, for every $\eta > 0$ and every $Q \in \mathcal Q$, 
    \[ \Exp{Q}{N_\eta(Q)} = \sum\limits_{t=1}^\infty Q\lrs{ N_\eta(Q) \ge t } < \infty. \numberthis\label{eq:warmstartcondition} \]
In the following lemmas, we will verify this condition for the different e-processes. 

\begin{lemma}\label{lem:warmstart_cor:1}
    The e-process from Corollary~\ref{cor:1} satisfies the integrable convergence condition.
\end{lemma}
\begin{proof}
    Fix $\eta > 0$ and $Q = N(m,1)$. Recall from~\eqref{eq:klinfsg} that for the setting in Corollary~\ref{cor:1}, $0 < \KLinf(Q,\calP)\le \KL(Q, N(0,1))= m^2/2$, and   $E^*_Q(X) = e^{mX - m^2/2}$. With these, bound the summand in~\eqref{eq:warmstartcondition} as below:
    \begin{align*}
        Q\lrs{ N_\eta(Q) \ge t } &\le Q\lrs{ \exists k \ge t-1: \frac{1}{k}\sum\limits_{i=1}^k \log E^*_Q(X_i) < (1-\eta)\KLinf(Q,\calP) }\\
        &\le Q\lrs{ \exists k \ge t-1: \frac{1}{k}\sum\limits_{i=1}^k \lrp{mX_i - \frac{m^2}{2}} < \frac{m^2}{2} - \eta\KLinf(Q,\calP) }\\
        &\le \sum\limits_{k=t-1}^\infty  Q\lrs{ \frac{1}{k}\sum\limits_{i=1}^k \lrp{X_i - {m}} < - \frac{\eta}{m}\KLinf(Q,\calP) } \tag{$m > 0$}\\
        &\le \sum\limits_{k=t-1}^\infty e^{-ck\frac{\eta^2}{m^2}\operatorname{KL}^2_{\operatorname{inf}}(Q,\calP)},
    \end{align*}
    for some constant $c > 0$, where the last inequality follows from the Gaussian Chernoff bound. Summing the above bound over $t$, we get~\eqref{eq:warmstartcondition}.  
\end{proof}

\begin{lemma}\label{lem:warmstart_cor:2}
    The e-process from Corollary~\ref{cor:2} satisfies the integrable convergence condition.
\end{lemma}
\begin{proof}
    Fix $\eta > 0$ and $Q = N(m_Q,1)$, for some $m_Q > 0$. As in~\eqref{eq:klinfsg}, for the setting in Corollary~\ref{cor:2} as well, $0 < \KLinf(Q,\calP)\le \KL(Q, N(0,1))= m^2_Q/2$. Further, recall that, in this setting, the log-e-process equals
    \[ \log E_k =   \frac{1}{2(k+1)} \lrp{\sum\limits_{i=1}^k X_i}^2  -\frac{1}{2}\log\lrp{k+1} \numberthis\label{eq:E_n} \]

    Next, fix $2\eta \KLinf(Q,\calP) > \epsilon > 0$. There exists a constant $T_\epsilon$ such that
    \[ \lrp{\frac{k+1}{k^2}}\log(k+1)+ \frac{m^2_Q}{k} - 2\eta\lrp{1+\frac{1}{k}}\KLinf(Q,\calP) \le -\epsilon, ~\text{ for all } k \ge T_{\epsilon}.\numberthis\label{eq:Teps} \]
    Then,
    \[ \Exp{Q}{N_\eta(Q)} \le T_{\epsilon} + \sum\limits_{t= 1+T_{\epsilon}}^\infty Q\lrs{N_\eta(Q) \ge t}. \numberthis\label{eq:warmstartbound1} \]
    To bound the summand above, consider the following: 
    \begin{align*}
        Q\lrs{ N_\eta(Q) \ge t } 
        &\le Q\lrs{ \exists k \ge t-1: \frac{\lrp{\sum\limits_{i=1}^k X_i}^2}{2k(k+1)}   -\frac{\log\lrp{k+1}}{2k} < (1-\eta)\KLinf(Q,\calP) }\\
        &\le \sum\limits_{k=t-1}^\infty Q\lrs{ \lrp{\Bar{X}_k}^2   - {m^2_Q} <  -\epsilon }\\
        &\le \sum\limits_{k=t-1}^\infty Q\lrs{ |\Bar{X}_k| \le \sqrt{m^2_Q - \epsilon} }\\
        &\le \sum\limits_{k=t-1}^\infty Q\lrs{ \Bar{X}_k - m_Q \le \sqrt{m^2_Q - \epsilon} - m_Q }, \numberthis\label{eq:bound1}
    \end{align*}
    where $\Bar{X}_k := (\sum_{i\le k} X_i) / k$. Without loss of generality, we consider $\epsilon < m^2_Q$ as otherwise, the right hand side above is bounded by $0$. Next, since $X_i \sim N(m_Q, 1)$, we have $\Bar{X}_k \sim N(m_Q, 1/k)$. Using this to further bound the right hand side in~\eqref{eq:bound1}, we get
    \begin{align*}
        Q\lrs{N_\eta(Q) \ge t} &\le \sum\limits_{k=t-1}^\infty 0.5 e^{-c k \lrp{\sqrt{m^2_Q-\epsilon} - m_Q}^2},
    \end{align*}
    for some constant $c > 0$, where the last bound follows from Hoeffding's inequality (or Chernoff bound).
    Using this in~\eqref{eq:warmstartbound1}, we get
    \[ \Exp{Q}{N_\eta(Q)} \le T_\epsilon + \sum\limits_{t=1+T_\epsilon}^\infty\sum\limits_{k=t-1}^\infty 0.5 e^{-c' k} < \infty, \]
    for a positive constant $c'$, completing the proof.
\end{proof}

\begin{lemma}\label{lem:warmstart_cor:3}
    The e-process from Corollary~\ref{cor:3} satisfies the integrable convergence condition. 
\end{lemma}
\begin{proof}
    First, observe that for $Q= \delta_0$, $F_n = \delta_0$ deterministically, and 
    \[ \KLinf(\delta_0,m_0) = -\log(1-m_0) > 0. \]
    Therefore, in this setting, 
    \[ \frac{1}{n} \log E_n = -\log(1-m_0) - \frac{\log(n+1) + 1}{n}, \]
    deterministically. So for every $\eta > 0$, $N_\eta(\delta_0)$ is finite and deterministic. 
    
    Henceforth, we consider $Q \ne \delta_0$, and hence, $m_Q > 0$. Fix $\eta > 0$, $Q \in \mathcal Q$, and $\epsilon \in (0, e^{\KLinf(Q,m_0)}+e^{-\KLinf(Q,m_0)}\tfrac{1-m_Q}{1-m_0} \wedge \eta \KLinf(Q,m_0))$.  Then, there exists a constant $T_\epsilon$ such that
    \[ \tfrac{1}{k} \lrp{1+\log(k+1)} - \eta \KLinf(Q, m_0) \le -\epsilon, ~\text{ for all } k \ge T_{\epsilon},\numberthis\label{eq:Teps1} \]
    and we get a bound as in~\eqref{eq:warmstartbound1}. Thus, it suffices to bound the summand in~\eqref{eq:warmstartbound1} as following for $t\ge T_\epsilon$: 
    \begin{align*}
        &Q\lrs{ N_\eta(Q) \ge t } \\
        &\le Q\lrs{ \exists k \ge t-1: ~ \KLinf(F_k, m_0) - \tfrac{1}{k} \lrp{1+ \log(k+1)} < (1-\eta) \KLinf(Q,m_0) },
    \end{align*}
    where we used the definition of $N_\eta(Q)$ and that $\log E_k = k\KLinf(F_k, m_0) - \log(k+1) - 1$. In the following, we will bound the right hand side above, and show that, when summed across $t$, it is still finite. Consider the following inequalities:
    \begin{align*}
        &Q\big[\exists k \!\ge\! t\!-\!1:\KLinf(F_k, m_0) \!-\! \tfrac{1}{k}\!\lrp{1\!+\!\log(k+1)} < (1\!-\!\eta) \KLinf(Q,m_0) \big]\\
        &\qquad\le Q\lrs{ \exists k \ge t-1: ~ \KLinf(F_k, m_0) - \KLinf(Q,m_0) \le -\epsilon } \tag{From~\eqref{eq:Teps1}}\\
        &\qquad\le \sum\limits_{k=t-1}^\infty  Q\lrs{\KLinf(F_k, m_0) - \KLinf(Q,m_0) \le -\epsilon }\tag{Union bound}\\
        &\qquad\le \sum\limits_{k=t-1}^\infty e^{-\frac{k\epsilon^2}{8 C}},\tag{Lemma~\ref{lem:subexpconcKLinf}}
    \end{align*}
    where $C > 0$. Using the above bound in~\eqref{eq:warmstartbound1}, we get
    \[ \Exp{Q}{N_\eta(Q)} \le T_\epsilon + \sum\limits_{t=1+T_\epsilon}^\infty \sum\limits_{k=t-1}^\infty e^{-\frac{k\epsilon^2}{8C}} < \infty,  \] 
    completing the proof.
\end{proof}

\begin{lemma}\label{lem:warmstart_cor:4}
    The e-process from Corollary~\ref{cor:4} satisfies the integrable convergence condition. 
\end{lemma}
\begin{proof}
    First, fix $\eta > 0$, $Q \in \mathcal Q$, and $\epsilon \in (0,e^{\KLinf(Q,\calP)} \wedge \eta \KLinf(Q,\calP))$.  Then, there exists a constant $T_\epsilon$ such that
    \[ \tfrac{K+1}{k} \lrp{1+\log(k+1)} - \eta \KLinf(Q, \calP) \le -\epsilon, ~\text{ for all } k \ge T_{\epsilon},\numberthis\label{eq:Teps2} \]
    and we get a bound as in~\eqref{eq:warmstartbound1}. Thus, it suffices to bound the summand in~\eqref{eq:warmstartbound1} as following for $t\ge T_\epsilon$: 
    \begin{align*}
        &Q\lrs{ N_\eta(Q) \ge t }\\
        &\le Q\lrs{ \exists k \ge t-1: ~ \KLinf(F_k, \calP) - \tfrac{K+1}{k}\lrp{1+\log(k+1)} < (1-\eta) \KLinf(Q,\calP) },
    \end{align*}
    where we used the definition of $N_\eta(Q)$ and that $\log E_k \ge  k\KLinf(F_k, \calP) - (K+1)\log(k+1) - 1$. In the following, we will bound the right hand side above, and show that, when summed across $t$, it is still finite. Consider the following inequalities:
    \begin{align*}
        &Q\big[ \exists k \!\ge\! t\!-\!1: \KLinf(F_k, \calP) \!-\! \tfrac{K+1}{k}\!\lrp{1+ \log(k+1)} \!<\! (1\!-\!\eta) \KLinf(Q,\calP) \big]\\
        &\quad\le Q\lrs{ \exists k \ge t-1: \KLinf(F_k, \calP) \!-\! \KLinf(Q,\calP) \le -\epsilon } \tag{From~\eqref{eq:Teps2}}\\
        &\quad\le \sum\limits_{k=t-1}^\infty  \!Q\lrs{\KLinf(F_k, \calP) - \KLinf(Q,\calP) \le -\epsilon }\tag{Union bound}\\
        &\quad\le \sum\limits_{k=t-1}^\infty e^{-\frac{k\epsilon^2}{8 C'}},\tag{Lemma~\ref{lem:subexpconcKLinf_finitelyconstrained}}
    \end{align*}
    for a constant $C' > 0$. Using the above bound in~\eqref{eq:warmstartbound1}, we get
    \[ \Exp{Q}{N_\eta(Q)} \le T_\epsilon + \sum\limits_{t=1+T_\epsilon}^\infty \sum\limits_{k=t-1}^\infty e^{-\frac{k\epsilon^2}{8C'}} < \infty,  \] 
    completing the proof.
\end{proof}

\subsection{More examples: Composite one-sided sub-Gaussian null, composite Gaussian alternative} \label{app:onesidedcompnullpcompalt}
In this section, we demonstrate the tightness of our lower bounds in Theorems~\ref{th:lbsmallalpha} and~\ref{th:lowerbound2} for testing the null, which is the one-sided 1-sub-Gaussian class considered in  Section~\ref{sec:ub_nonpcompnullparampointalt}, against a composite alternative consisting of all Gaussian distributions with a positive mean and unit variance. Formally, 
\[ \mathcal P = \lrset{P: \Exp{P}{e^{\theta X - \frac{\theta^2}{2}}}\le 1, \quad \forall \theta\ge 0} \quad\text{and}\quad \mathcal Q = \{ N(m,1): m > 0 \}. \numberthis\label{app_eq:PQ2} \]
Note that this example is close to the setting in Section~\ref{sec:ub_nonpcompnullpcompalt}. However, our null is much larger. In particular, the $\calP$-e-processes considered in Section~\ref{sec:ub_nonpcompnullpcompalt} are no longer $\calP$-e-processes. 

Recall from Section~\ref{sec:ub_nonpcompnullparampointalt} that $m_P \le 0$ for all $P\in\cal P$. Moreover, for any $Q\in\mathcal Q$ with mean $m_Q > 0$, we have $\KLinf(Q,\mathcal P) \le m^2_Q/2$. In the following corollary of Theorem~\ref{th:generalasupperboundcond}, we present a test that is optimal for every (unknown) $Q\in\mathcal Q$, in the $\alpha\rightarrow 0$ limit. 
\begin{corollary}\label{app_cor:2}
    There exists an $\alpha$-correct, power-one stopping rule $\tau_\alpha$ for testing $\calP$ vs. $\mathcal Q$ as defined in~\eqref{app_eq:PQ2}, such that for every $Q\in\mathcal Q$,
    \[ Q\lrs{\limsup\limits_{\alpha \rightarrow 0}~ \frac{\tau_\alpha}{\log({1}/{\alpha})} \le \frac{1}{\KLinf(Q,\mathcal P)} } = 1.\]
\end{corollary}
\begin{proof} For $\theta > 0$, define
$E'_n(\theta) := \exp\{\theta \sum\nolimits_{i=1}^n X_i - \frac{n\theta^2}{2}\}$, and  consider a mixture of $E'_n(\cdot)$, mixed  with a  truncated and scaled $N(0,1)$ prior on $\theta \in [0,\infty)$, that is, 
\[ E_n := \frac{2}{\sqrt{2\pi}}\int\limits_{0}^\infty E'_n(\theta) e^{-\frac{\theta^2}{2}} d\theta = \frac{2}{\sqrt{n+1}}~ e^{\frac{n^2 (\hat{m}_n)^2}{2(n+1)}} \Phi\lrp{\frac{\sum_{i\le n} X_i}{\sqrt{n+1}}}, \]
where $\hat{m}_n \!:= \!\frac{1}{n}\sum\nolimits_{i=1}^n X_i$, and $\Phi(\cdot)$ is the cumulative density function of $N(0,1)$. $E_n$ is a $\calP$-e-process that satisfies Assumption~\ref{assmp:E}, that is, for every $Q\in\mathcal Q$ satisfies $Q$-a.s.:
\begin{align*}
    \liminf\limits_{n\rightarrow \infty}~\lrp{\frac{1}{n} \log E_n} = \frac{m^2_Q}{2} \ge \KLinf(Q,\mathcal P).
\end{align*}
The bound now follows from Theorem~\ref{th:generalasupperboundcond}. 
\end{proof}
We now present a sequential test that incurs an iterated-logarithmic complexity along every sequence $Q_n \in\mathcal Q$ satisfying $\KLinf(Q_n,\calP) \rightarrow 0$ as $n\rightarrow \infty$, thus achieving the lower bound in Theorem~\ref{th:lowerbound2} along every directly-approaching sequence of alternatives. 
\begin{proposition}
    For every fixed $\alpha\in (0,1/2)$, there exists a power-one, $\alpha$-correct stopping rule $\tau_\alpha$ to test $\calP$ vs. $\mathcal Q$, as defined in~\eqref{app_eq:PQ2}, which satisfies the following for a universal constant $C$ and every sequence $Q_n\in\mathcal Q$ such that $\KLinf(Q_n,\calP) \rightarrow 0$ as $n \rightarrow \infty$: 
    \[ \limsup\limits_{n\rightarrow \infty}~\frac{\Exp{Q_n}{{\tau}_\alpha}}{\operatorname{KL^{-1}_{inf}}(Q_n,\mathcal P) \log\log\operatorname{KL^{-1}_{inf}}(Q_n,\mathcal P) } \le  \frac{C(1-\alpha)}{(1-2\alpha)^2}.  \]
\end{proposition}
\begin{proof}
For suitable nonnegative constants $c_1, c_2$, and $c_3> 1$, and for each $n\in\N$,  define
\[ L_n(\alpha) := \frac{1}{n}\sum\limits_{i=1}^n X_i 
 - c_1 U(n,\alpha), \quad \text{where}\quad U(n,\alpha)=\sqrt{\frac{\log\log(e+c_2 n) + \log\frac{c_3}{\alpha}}{n}}. \]
Observe that in $\cal P$, the upper tails of the distributions are controlled, while the lower tails can be arbitrary. In addition, recall that an assumption on the upper tail of $P$ gives a lower confidence sequence for $m_P$. In fact, from \citet[Equation 1.2]{howard2021time}, we have for any $P\in\mathcal P$ with mean $m_P$ ($\le 0$) that  
 \[ P\lrs{\forall n, ~m_P \ge L_n(\alpha)} \ge 1-\alpha. \]
Let $\tilde{\tau}_\alpha := \min\lrset{n: L_n(\alpha) > 0}$. Then, $\tilde{\tau}_\alpha$ is $\alpha$-correct since
 \[ \sup\limits_{P\in\mathcal P}~P\lrs{\tilde{\tau}_\alpha < \infty} = \sup\limits_{P\in\mathcal P}~P\lrs{\exists n: L_n(\alpha) > 0} \le \sup\limits_{P\in\mathcal P}~P\lrs{\exists n: L_n(\alpha) > m_P} \le \alpha. \]
In addition, $\tilde{\tau}_\alpha$ has power-one since 
\[ \inf\limits_{Q\in\mathcal Q}~Q[\tilde{\tau}_\alpha < \infty] = \inf\limits_{Q\in\mathcal Q}~Q\lrs{\exists n: L_n(\alpha) > 0} \overset{(a)}{=} 1, \]
where $(a)$ follows since by SLLN, under $Q$, $L_n(\alpha)\rightarrow m_Q > 0$ a.s. Now, consider the event
\[ \mathcal E(\alpha', m) := \lrset{\forall n\in \N, ~ \abs{\frac{1}{n}\sum\limits_{i=1}^n X_i - m} \le c_1 U(n,\alpha') }. \]
Again, from~\citet[Equation 1.2]{howard2021time}, it follows that for all $Q\in\mathcal Q$,  $Q[\mathcal E(\alpha',m_Q)] \ge 1-\alpha'$. Then, on $\mathcal E(\alpha, m_Q)$, 
\begin{align*} 
    \tilde{\tau}_\alpha
    &= \min\lrset{n :   c_1 U(n,\alpha) < \frac{1}{n}\sum_i X_i}\\ 
    &\le \min\lrset{n : 2 c_1 U(n,\alpha)  < m_Q}\\
    &\le \min\lrset{n : n \ge \frac{4 c^2_1}{m^2_Q} \log\log(e+c_2 n) + \frac{4 c^2_1}{m^2_Q} \log\lrp{\frac{c_3}{\alpha}} }.
\end{align*}

From Lemma~\ref{lem:boundontaualpha} (with $\gamma = 0.5$), we have, 
\[\tilde{\tau}_\alpha \le 1 + \frac{32 c^2_1}{ m^2_Q} \log \log \lrp{e+\frac{8 c^2_1 c_2}{m^2_Q}} + \frac{8 c^2_1}{m^2_Q} \log\lrp{\frac{c_3}{\alpha}}, \quad \text{pathwise on }\mathcal E(\alpha, m_Q).\]
Letting the right hand side of the above bound be $T(\alpha, Q)$, we see that for $\alpha' < \alpha$, 
\[ \frac{T(\alpha',Q)}{\log\tfrac{1}{\alpha'}} \le \frac{T(\alpha, Q)}{\log\tfrac{1}{\alpha}}. \]
Theorem~\ref{th:hptoexp} shows the existence of $\tau_\alpha$ (as a meta algorithm using $\tilde{\tau}_\alpha$ as a sub-routine) such that
\[ \Exp{Q}{\tau_\alpha} \le \frac{4(1-\alpha)}{(1-2\alpha)^2}\lrp{1 + \frac{32 c^2_1}{m^2_Q} \log \log \lrp{e+\frac{8 c^2_1 c_2}{m^2_Q}} + \frac{8 c^2_1}{m^2_Q} \log\lrp{\frac{c_3}{\alpha}}}.\]
The bound in the proposition easily follows from the above, along with the fact that $\KLinf(Q,\calP) \le \frac{m^2_Q}{2}$, completing the proof. 
\end{proof}

\section{High-probability to expected sample complexity}
In this section, we recall a meta-algorithm due to \cite{chen2015optimal}, which converts an $\alpha$-correct, power-one sequential test with a high-probability bound on the stopping time to one with a similar bound in expectation. 

\begin{theorem}[\text{\citet[Theorem H.5]{chen2015optimal}}]\label{th:hptoexp}
For every $\alpha \in (0,0.5)$, suppose that we have a family $\{\tau_\alpha\}$ of $\alpha$-correct power-one sequential tests, and for every $Q\in\mathcal Q$, events $\calE_{\alpha,Q}$ with $Q[\calE_{\alpha,Q}] \ge 1- \alpha$, such that 
\[  Q\lrs{\tau_\alpha < \infty ~|~ \calE_{\alpha,Q}} = 1, \quad\text{and}\quad \Exp{Q}{\tau_\alpha|\calE_{\alpha,Q}} \le T(\alpha,Q),\]
where the bound $T(\alpha,Q)$ satisfies the following: there exists $\alpha_0 < 0.5$ such that for all $0 <\alpha' < \alpha < \alpha_0$, and all $Q\in\mathcal Q$,
\[ \frac{T(\alpha',Q)}{\log({1}/{\alpha'})} \le \frac{T(\alpha, Q)}{\log({1}/{\alpha})}. \numberthis\label{eq:CL-monotone}\] 
Then, for every $\alpha \in (0,\alpha_0)$, there exists an $\alpha$-correct, power-one sequential test (denoted by $\tilde{\tau}_\alpha$) such that 
$$\Exp{Q}{\tilde{\tau}_\alpha} \le \frac{4(1-\alpha)}{(1-2\alpha)^2}T(\alpha,Q), \qquad  \forall Q\in\mathcal Q.$$ 
\end{theorem}

The proof of Theorem~\ref{th:hptoexp} proceeds by designing a meta-algorithm that simulates multiple copies of the given sequential test with different error probabilities in parallel. At a high level, this meta-algorithm stops and rejects the null at the first time when one of the ingredient tests stops. In the following, we give the construction of the meta-algorithm and a proof for its sample complexity bound for completeness.

\begin{proof}
Consider an algorithm $\A$ that simulates $\tau_{\alpha_i}$ for $i\in\N$ with $\alpha_i= \frac{\alpha}{2^i}$, as described next. In step $r\in\N$, $\tau_{\alpha_i}$ generates an independent sample only if $2^{i-1}$ divides $r$. For concreteness, the test $\tau_{\alpha_1}$ runs (and generates samples) only at steps $1$, $2$, $3$, $4$, $\dots$, test $\tau_{\alpha_2}$ at steps $2$, $4$, $6$, $\dots$, etc. If at any time $r$, multiple tests $\tau_{\alpha_i}$ generate samples, then $\A$ generates an independent sample for each of them and feeds it to them. Thus, $\A$ can generate multiple samples at each step. $\A$ stops at the first time that any of the ingredient tests $\tau_{\alpha_i}$ stops, and rejects the null. Let $\tilde{\tau}_\alpha$ denote the stopping time of $\A$, or equivalently, the total number of samples generated by $\A$ before it stops. 

\paragraph*{$\alpha$-correctness of $\A$}  Clearly, $\A$ is a power-one sequential test that satisfies
\begin{align*}
    \sup\limits_{P\in\mathcal P} P\lrs{ \tilde{\tau}_\alpha < \infty } \le \sup\limits_{P\in\mathcal P} P\lrs{\exists i \in \N: \tau_{\alpha_i} < \infty} \le \sup\limits_{P\in\mathcal P} \sum\limits_{i=1}^\infty P\lrs{\tau_{\alpha_i} < \infty}\le \sum\limits_{i=1}^\infty \frac{\alpha}{2^i} = \alpha. 
\end{align*}

\paragraph*{Sample complexity of $\A$ } From our assumptions about the base sequential tests $\tau_{\alpha_i}$, there exist events $\mathcal E_i$ such that under the unknown data-generating distribution $Q[\mathcal E_i] \ge 1-\alpha_i$. Let $\mathcal F_i = (\cap_{j < i} \mathcal E^c_j) \cap \mathcal E_i$. 
Clearly, $\lrset{\mathcal F_i}$ forms a partition with
\[ Q\lrs{\mathcal F_i} \le \prod\limits_{j=1}^{i-1}\frac{\alpha}{2^j}  \le \alpha^{i-1}.\]

Moreover, since $\A$ feeds each $\tau_{\alpha_i}$ with independent samples, we have 
$$\Exp{Q}{\tau_{\alpha_i} | \mathcal F_i} = \Exp{Q}{\tau_{\alpha_i} | \mathcal E_i} \le  T(\alpha_i, Q),$$
where the last inequality follows from the assumption on $\tau_{\alpha_i}$. Next, observe that for $j\ne i$, for each sample fed to $\tau_{\alpha_i}$  $\A$ feeds $2^{i-j}$ samples to $\tau_{\alpha_j}$. Thus, 
\[ \Exp{Q}{\tilde{\tau}_\alpha|\mathcal F_i} \le \Exp{Q}{\tau_{\alpha_i} | \mathcal F_i}\sum\limits_{j=1}^\infty 2^{i-j} \le  T\lrp{\alpha_i, Q}\sum\limits_{j=1}^\infty 2^{i-j} \le T\lrp{\alpha_i, Q} 2^i. \]

Now, consider the following
\begin{align*}
    \Exp{Q}{\tilde{\tau}_\alpha} = \sum_i\Exp{Q}{\tilde{\tau}_{\alpha} | \mathcal F_i} Q[\mathcal F_i]&\le \sum_i T\lrp{\alpha_i, Q} 2^i \alpha^{i-1} \\
    &= 2\sum\limits_i T(\alpha_i, Q)(2\alpha)^{i-1}\\
    &\le 2\sum\limits_i \frac{\log\frac{1}{\alpha_i}}{\log \frac{1}{\alpha}}T(\alpha, Q)(2\alpha)^{i-1}\\
    &= 2\sum\limits_i \lrp{\frac{\log\frac{1}{\alpha} + i \log 2}{\log\frac{1}{\alpha}}}T(\alpha, Q)(2\alpha)^{i-1}\\
    &\le 2\sum\limits_i \lrp{1+i}T(\alpha, Q)(2\alpha)^{i-1}\\
    &\le  \frac{4(1-\alpha)}{(1-2\alpha)^2} T(\alpha,Q),
\end{align*}
completing the proof.
\end{proof}

\section{$\KLinf$ for finitely generated class}\label{sec:KLinfdual}
In this section, we prove certain properties of $\KLinf$ for the finitely-constrained class studied in Section~\ref{sec:ub_finitelygenhypo}, including the alternative formulation from Theorem~\ref{th:KLinfdual}. Recall that for any $Q\in M_1(\calX)$, $\KLinf(Q, \mathcal P)$ is the optimal value of the following optimization problem: 
\begin{align*}
    \inf\limits_{P\in M_+(\mathcal X)} \qquad &\KL(Q,P)\\
    \text{s.t.} \qquad & \Exp{P}{\phi_i(X)} \le 0, \quad \forall i \in [K-1]\\
    & \Exp{P}{\phi_K(X)} =0\\ 
    & \Exp{P}{1} = 1,  
\end{align*}
where, for $P\in M_+(\calX)$, we let 
\[ \KL(Q,P) := \begin{cases}
    \int\limits_{\calX} \log\lrp{\frac{dQ}{dP}} dQ, & Q\ll P,\\
    +\infty, & Q\not\ll P.
\end{cases} \]

\subsection{Proof of Theorem~\ref{th:KLinfdual}}\label{app:finitelyconstraineddual}
\begin{proof}
Since, by assumption, $\KLinf(Q,\calP)<\infty$, $\calP$ is compact, and $\KL(Q,\cdot)$ is lower semi-continuous, there exists some $P_0\in\mathcal P$ such that $\KL(Q,P_0)<\infty$, and hence $Q\ll P_0$. Consequently, every $\mathcal P$-quasi-sure inequality also holds $Q$-almost surely.

We first relate $\KLinf(Q,\calP)$ to optimization over $\mathcal P$-e-variables. Let $C_b(\R)$ denote the set of bounded continuous functions on $\R$. For $g\in C_b(\R)$ and $P\in\mathcal P$, define
\[
    \Psi(g,P)
    :=
    \mathbb E_Q[g(X)]
    -\log \mathbb E_P[e^{g(X)}].
\]
Recall, by Donsker--Varadhan variational representation of $\KL(Q,P)$ (eg., \citet[Equation 1.1]{dupuis2022formulation}), 
\[
    \KL(Q,P)
    =
    \sup_{g\in C_b(\R)} \lrset{\mathbb E_Q[g(X)]
    -\log \mathbb E_P[e^{g(X)}]}.
\]
Next, observe that for every fixed $g$, the map $P\mapsto\Psi(g,P)$ is continuous under weak convergence and convex in $P$. Moreover, for every fixed $P$, the map $g\mapsto\Psi(g,P)$ is concave and continuous in the supremum norm. Finally, since $\mathcal P$ is convex and compact in the weak topology, Sion's minimax theorem \citep{sion1958general} therefore gives
\begin{align}
    \KLinf(Q,\calP)
    &=
    \inf_{P\in\mathcal P}~\sup_{g\in C_b(\R)}~ \Psi(g,P)
    \nonumber\\
    &=
    \sup_{g\in C_b(\R)}~\inf_{P\in\mathcal P}~
    \Psi(g,P)
    \nonumber\\
    &=
    \sup_{g\in C_b(\R)}
    \left\{
        \mathbb E_Q[g(X)]
        -
        \log\sup_{P\in\mathcal P}\mathbb E_P[e^{g(X)}]
    \right\}.
    \label{eq:robust-DV}
\end{align}

Fix $g\in C_b(\R)$ and define
\[
    c_g
    :=
    \sup_{P\in\mathcal P}\mathbb E_P[e^{g(X)}],
    \qquad
    E_g(x)
    :=
    \frac{e^{g(x)}}{c_g}.
\]
Clearly, $E_g$ is a $\mathcal P$-e-variable. By
\citet[Theorem~3.1 and Lemma~3.2]{larsson2025variables}, there exists
$(\lambda_g,\kappa_g)\in\Pi\cap\mathcal K$ such that
\[
    E_g
    \le
    S^{\lambda_g,\kappa_g}
    \qquad \mathcal P\text{-quasi surely}.
\]
As observed above, this inequality also holds $Q$-almost surely. Therefore,
\begin{align*}
    \sup_{(\lambda,\kappa)\in\Pi\cap\mathcal K}
    \mathbb E_Q\!\left[\log S^{\lambda,\kappa}(X)\right]
    &\ge
    \mathbb E_Q\!\left[
        \log S^{\lambda_g,\kappa_g}(X)
    \right] \\
    &\ge
    \mathbb E_Q[\log E_g(X)] \\
    &=
    \mathbb E_Q[g(X)]
    -
    \log\sup_{P\in\mathcal P}\mathbb E_P[e^{g(X)}].
\end{align*}
Taking the supremum over $g\in C_b(\R)$ and using
\eqref{eq:robust-DV} yields
\begin{equation}
    \sup_{(\lambda,\kappa)\in\Pi\cap\mathcal K}
    \mathbb E_Q\!\left[\log S^{\lambda,\kappa}(X)\right]
    \ge \KLinf(Q,\calP).
    \label{eq:dual-lower}
\end{equation}

We next prove the reverse inequality. Fix
$(\lambda,\kappa)\in\Pi\cap\mathcal K$. By construction, $S^{\lambda,\kappa}$ is a $\mathcal P$-e-variable, so
\[
    \sup\limits_{P\in\calP}\mathbb E_P[S^{\lambda,\kappa}]\le 1.
\]
Fix $P\in\mathcal P$. If $\KL(Q,P)=\infty$, then
$\mathbb E_Q[\log S^{\lambda,\kappa}]\le\KL(Q,P)$ holds trivially. Otherwise, let
$R:=dQ/dP$. If $S^{\lambda,\kappa}=0$ with positive $Q$-probability, then
$\mathbb E_Q[\log S^{\lambda,\kappa}]=-\infty$ and again, there is nothing to prove. Otherwise,
\begin{align*}
    \mathbb E_Q[\log S^{\lambda,\kappa}]
    -
    \KL(Q,P)
    &=
    \mathbb E_Q\!\left[
        \log\frac{S^{\lambda,\kappa}}{R}
    \right] \\
    &\le
    \log\mathbb E_Q\!\left[\frac{S^{\lambda,\kappa}}{R}\right] \tag{Jensen's inequality} \\
    &=
    \log\int_{\{R>0\}} S^{\lambda,\kappa}\,dP \\
    &\le
    \log\mathbb E_P[S^{\lambda,\kappa}]
    \le 0.
\end{align*}
Hence
\[
    \mathbb E_Q[\log S^{\lambda,\kappa}]
    \le
    \KL(Q,P),
    \qquad \forall P\in\mathcal P.
\]
Taking the infimum over $P\in\mathcal P$ gives
\[
    \mathbb E_Q[\log S^{\lambda,\kappa}]
    \le
    \KLinf(Q,\mathcal P).
\]
Since $(\lambda,\kappa)\in\Pi\cap\mathcal K$ was arbitrary,
\begin{equation}
    \sup_{(\lambda,\kappa)\in\Pi\cap\mathcal K}
    \mathbb E_Q\!\left[\log S^{\lambda,\kappa}(X)\right]
    \le \KLinf(Q,\calP).
    \label{eq:dual-upper}
\end{equation}
Combining \eqref{eq:dual-lower} and \eqref{eq:dual-upper} proves the desired
representation.

We also note that the supremum in the alternative representation above is attained. To see this, first recall that  $\Pi\cap\mathcal K$ is compact \citep[Lemma~7.2]{larsson2025variables}. Next, using the reverse Fatou's lemma \citep{williams1991probability}, we will conclude that the map
\[
    (\lambda,\kappa)
    \longmapsto
    \mathbb E_Q[\log S^{\lambda,\kappa}(X)]
\]
is upper-semicontinuous, and hence,
attains the maximum (possibly non-unique) on the compact set $\Pi\cap\mathcal K$. We now verify the conditions for the reverse Fatou's lemma. 

To this end, let
$(\lambda_n,\kappa_n)\to(\lambda,\kappa)$ in $\Pi\cap\mathcal K$ (existence of such a sequence is guaranteed by compactness of $\Pi \cap \calK$). Then, pointwise, for every $x\in \calX$, 
\[
    \limsup_{n\to\infty}
    \log S^{\lambda_n,\kappa_n}(x)
    \le
    \log S^{\lambda,\kappa}(x),
\]
where we use the convention $\log 0=-\infty$ (the map $x\to \log x$ is upper-semicontinuous on $[0, \infty)$). Moreover, the compactness of
$\Pi\cap\mathcal K$ gives a constant $M<\infty$ such that
\[
    \log S^{\lambda,\kappa}(x)
    \le
    S^{\lambda,\kappa}(x)
    \le
    1+M\sum_{i=1}^K |\phi_i(x)|
\]
for every $(\lambda,\kappa)\in\Pi\cap\mathcal K$, $Q$-almost surely. The right-hand side is $Q$-integrable by assumption, showing that the maximum is attained. 

Finally, let $S^{\lambda_\star,\kappa_\star}$ denote any maximizer above. Then, for any $\mathcal P$-e-variable $E$, \citet[Theorem~3.1 and Lemma~3.2]{larsson2025variables} yields some $(\lambda,\kappa)\in\Pi\cap\mathcal K$ such that
$E\le S^{\lambda,\kappa}$ $Q$-almost surely. Hence
\[
    \mathbb E_Q[\log E]
    \le
    \mathbb E_Q[\log S^{\lambda,\kappa}]
    \le
    \mathbb E_Q[\log S^{\lambda_\star, \kappa_\star}].
\]
Hence, $S^{\lambda_\star, \kappa_\star}$ is also log-optimal. Since $\mathbb E_Q[\log S^{\lambda_\star, \kappa_\star}]=\KLinf(Q,\calP)<\infty$, we have $S^{\lambda_\star, \kappa_\star}>0$,  $Q$-almost surely. Thus, by \citet[Proposition~2.4]{larsson2024numeraire}, $S^{\lambda_\star, \kappa_\star}$ is the numeraire e-variable. In particular, applying the numeraire property to the constant e-variable $1$ gives
\begin{equation*}
    \mathbb E_Q\!\left[\frac{1}{S^{\lambda_\star, \kappa_\star}(X)}\right]
    \le 1,
\end{equation*}
completing the proof.
\end{proof}

\subsection{A concentration result}
\begin{lemma}\label{lem:subexpconcKLinf_finitelyconstrained}
    Consider $\calP$ and $\mathcal Q$ as in Section~\ref{sec:ub_finitelygenhypo}. Let $F_n$ denote the empirical distribution with $n$ i.i.d.\ samples from $Q \in \mathcal Q$ such that $0 < \KLinf(Q,\calP) < \infty$. Then there exists a constant $C' > 0$ such that for every $\epsilon > 0$, 
    \[ Q\lrs{ \KLinf(F_n,\calP) - \KLinf(Q,\calP) \le -\epsilon } \le e^{-n\epsilon^2/\lrp{4 (\epsilon + C')}}.\]
    Further, for $\epsilon \in (0,C')$,
    \[ Q\lrs{ \KLinf(F_n, \calP) - \KLinf(Q,\calP) \le -\epsilon } \le e^{-\frac{n\epsilon^2}{8 C'}}. \]
\end{lemma}
\begin{proof}
First note that Theorem~\ref{th:KLinfdual} applies to $F_n$ for every realization since for each observed $X_j$, choose $P_{X_j} \in \calP$ such that $P_{X_j}(\{X_j\}) > 0$ (existence guaranteed by Assumption~\ref{assmp:compactnull}), and let 
\[ \bar{P}_n := \frac{1}{n} \sum_j P_{X_j}. \]
Convexity of $\calP$ (Assumption~\ref{assmp:compactnull}) gives that $\bar{P}_n \in \calP$ and $F_n \ll \bar{P}_n$, and hence, $\KLinf(F_n, \calP) \le \KL(F_n,\bar{P}_n) < \infty$.
    
    For $x\in\R$, let $[\Phi(x)] = [\phi_1(x) ~\dots~ \phi_K(x)~ -\phi_K(x)]^T$ denote the column vector of the constraint functions evaluated at $x$. Further, for $({\bar{ \lambda}}, {\bar{ \kappa}}) \in \Pi \cap \calK$, let $[{\bar{ \lambda}~\bar{\kappa}}][\Phi(x)] = \sum_i \lambda_i \phi_i(x) + \kappa_1 \phi_K(x) - \kappa_2 \phi_K(x)$.  Let $(\bar{\lambda}^*, \bar{\kappa}^*)$ be the maximizers in $\KLinf(Q,\calP)$ (Theorem~\ref{th:KLinfdual}). Consider the following inequalities:
    \begin{align*}
        &Q\lrs{ \KLinf(F_n, \calP) - \KLinf(Q, \calP) \le -\epsilon }\\
        &= Q\lrs{ \max\limits_{(\bar{ \lambda}, \bar{ \kappa})\in \Pi\cap \calK}~ \frac{1}{n} \sum\limits_{i=1}^n \log(1+{[{\bar{\lambda} ~ \bar{\kappa}}]}[\Phi(X_i)]) - \Exp{Q}{\log(1+[\bar{\lambda}^*~\bar{\kappa}^*][\Phi(X)])} \le -\epsilon  }\\
        &\le Q\lrs{  \frac{1}{n} \sum\limits_{i=1}^n \log(1+[\bar{\lambda}^*~ \bar{\kappa}^*][\Phi(X_i)]) - \Exp{Q}{\log(1+[\bar{\lambda}^*~ \bar{\kappa}^*][\Phi(X)])} \le -\epsilon  }.
    \end{align*}
    Let 
    \[ Y_i := \log(1+[\bar{\lambda}^*~\bar{\kappa}^*][\Phi(X_i)]) - \Exp{Q}{\log(1+[\bar{\lambda}^*~\bar{\kappa}^*][\Phi(X)])}. \]
    Clearly, $\Exp{Q}{Y_i} = 0$. Further, 
    \[ \Exp{Q}{e^{Y_i}}  \le e^{-\KLinf(Q,\calP)}\Exp{Q}{1+[\bar{\lambda}^*~\bar{\kappa}^*][\Phi(X)]}\] and
    \[\Exp{Q}{e^{-Y_i}} = e^{\KLinf(Q,\calP)} \Exp{Q}{\frac{1}{1+[\bar{\lambda}^*~\bar{\kappa}^*][\Phi(X)]}} \le e^{\KLinf(Q,\calP)}, \]
    where the last inequality follows from optimality of $[\bar{\lambda}^*~\bar{\kappa}^*]$ and Theorem~\ref{th:KLinfdual}. Using Lemma~\ref{lem:bernstein}, for $\epsilon > 0$, we get the first bound with 
    \[ C' = e^{-\KLinf(Q,\calP)}\Exp{Q}{1+[\bar{\lambda}^*~\bar{\kappa}^*][\Phi(X)]} +  e^{\KLinf(Q,\calP)}.\] The second inequality follows from the previous one and the condition on $\epsilon$. 
\end{proof}


\section{Additional supporting results}
\begin{lemma}\label{lem:bernstein} Let $Y_1, \dots, Y_n$ be i.i.d.\ samples from a distribution $Q$ such that $\Exp{Q}{Y_i} = 0$, $\Exp{Q}{e^{Y_i}} \le B_+$, and $\Exp{Q}{e^{-Y_i}} \le B_-$. Then, for any $\epsilon > 0$, 
\[ Q\lrs{ \frac{1}{n} \sum\limits_{i=1}^n Y_i \le -\epsilon } \le e^{-\frac{n\epsilon^2}{4(B_+ + B_-+ \epsilon)}}. \]
\end{lemma}
\begin{proof}
    Since $\Exp{Q}{e^{Y_i}} \le B_+$ and $\Exp{Q}{e^{-Y_i}} \le B_-$, we have
    \[ \Exp{Q}{e^{|Y_i|}} \le \Exp{Q}{e^{Y_i}} + \Exp{Q}{e^{-Y_i}} \le B_+ + B_-.  \]
    Further, for $\theta \in \R$, $i\ge 1$,  and $j\ge 2$, we have
    \begin{align*}
        \Exp{Q}{|\theta Y_i|^j} &= |\theta|^j \Exp{Q}{|Y_i|^j} \le |\theta|^j j! \Exp{Q}{e^{|Y_i|}} \le (B_+ + B_-) j! |\theta|^j.  \numberthis\label{eq:momentbound}
    \end{align*}
    Finally, consider the following for $\theta \in (-1,1)$,
    \begin{align*}
        \Exp{Q}{e^{\theta Y_i}} = \sum\limits_{j=0}^\infty \frac{\Exp{Q}{(\theta Y_i)^j}}{j!} &= 1 + \sum\limits_{j=2}^\infty \frac{\Exp{Q}{(\theta Y_i)^j}}{j!}\\
        &\le 1 + \sum\limits_{j=2}^\infty \frac{\Exp{Q}{|\theta Y_i|^j}}{j!}\\
        &\le 1 + \sum\limits_{j=2}^\infty {(B_+ + B_-)|\theta|^j}\tag{Using~\eqref{eq:momentbound}}\\
        &= 1 + (B_+ + B_-) \frac{\theta^2}{1-|\theta|}\\
        &\le  e^{(B_+ + B_-)\frac{\theta^2}{1-|\theta|}} \tag{$1+x \le e^x$, for all $x\in \R$},
    \end{align*}
    which further implies
    \[ \log \Exp{Q}{e^{\theta Y_i}} \le (B_+ + B_-) \frac{\theta^2}{1-|\theta|}, \quad \text{ for } \theta \in [-1,1]. \numberthis\label{eq:MGFBound} \]

Next, let $\epsilon > 0$ and $\theta \in (0,1)$. Using the above bound on the log moment generating function,
\begin{align*}
    Q\lrs{ \frac{1}{n} \sum\limits_{i=1}^n Y_i \le -\epsilon } &\le Q\lrs{ -\theta \sum\limits_{i=1}^n Y_i \ge n\theta \epsilon } \le e^{-n \lrp{ \theta \epsilon - (B_+ + B_-)  \frac{\theta^2}{1-|\theta|}} }. \numberthis\label{eq:tailbound}
\end{align*}
Let  
\[ \theta^* = \argmax\limits_{\theta \in [0,1]}~ f(\theta), \qquad \text{ where } \qquad f(\theta) =  \lrs{\theta \epsilon - (B_+ + B_-) \frac{\theta^2}{1-\theta}}.\]
It is easy to show that
\[ \theta^* = 1 - \frac{1}{\sqrt{1 + \frac{\epsilon}{B_+ +  B_-}}},  \]
which lies in $(0,1)$. Using this,
\begin{align*}
    f(\theta^*) 
    = \epsilon \theta^* - (B_+ + B_-) \frac{{\theta^*}^2}{1-\theta^*} &= \frac{\epsilon^2}{B_+ + B_-} {\lrp{1 + \sqrt{1 + \frac{\epsilon}{B_+ + B_-}}}^{-2}}\\
    &\ge \frac{\epsilon^2}{4 (B_+ + B_- + \epsilon)},
\end{align*}
where the above inequality holds for all $\epsilon > 0$. Since the bound in~\eqref{eq:tailbound} holds for all $\theta\in (0,1)$, we get the following by optimizing the bound over $\theta$:
\[ Q\lrs{ \frac{1}{n} \sum\limits_{i=1}^n Y_i \le -\epsilon } \le e^{-\frac{n\epsilon^2}{4(B_+ + B_- + \epsilon)}}, \]
completing the proof. 
\end{proof}

\begin{lemma}\label{lem:DH}
    For a probability measure $Q$, let $\{M_n\}_n$ be a stochastic process, and $G(Q) > 0$, $H(Q) > 0$ be such that for all $n\in \N$ and for all $\alpha \in (0,1)$,  
    \[ Q\lrs{\sqrt{M_n} \le \sqrt{G(Q)} - \sqrt{H(Q)}\sqrt{\frac{C'\log(1/\alpha)}{n}}} \le \alpha. \]
    Then,
    \[ Q\lrs{ \exists \ell\ge 1: \sqrt{M_{2^\ell}} \le \sqrt{G(Q)}  - \sqrt{H(Q)}\sqrt{\frac{ C'\lrp{\log(\pi^2/6) + \log(1/\alpha) + 2\log \ell} }{2^\ell}} } \le \alpha. \]
\end{lemma}
\begin{proof}
    Define the event,
    \[ \mathcal B_\ell := \lrset{\sqrt{M_{2^\ell}} \le \sqrt{G(Q)}  - \sqrt{H(Q)}\sqrt{\frac{ C'\log \tfrac{\pi^2\ell^2}{6\alpha} }{2^\ell}} }. \]
    Then, from the given condition, for all $\ell\ge 1$, $Q[\mathcal B_\ell] \le \tfrac{6\alpha}{\pi^2\ell^2}$. The conclusion now follows by taking a union bound over $\ell \ge 1$.
\end{proof}

\begin{lemma}\label{lem:boundontaualpha} 
Let $c_1$ and $c_2$ be positive constants and let $f:[0,1]\rightarrow \R_+$. For $\alpha \in (0,1)$, define
    \[ \tau_\alpha := \min \lrset{ n : n \ge \frac{1}{c_1} \log \log (e+c_2 n) + \frac{f(\alpha)}{c_1} }. \]
    Then, for every $\gamma \in (0,1)$, 
    \[ \tau_\alpha \le 2 + \frac{4}{\gamma c_1} \log \log \lrp{e+\frac{c_2}{\gamma c_1}} + \frac{f(\alpha)}{(1-\gamma) c_1}. \]
\end{lemma}
\begin{proof}
    First, notice that the function $t \to (\log\log(e+c_2 t))/t$ is non-increasing on $(0,\infty)$, that is, 
    \[ \frac{\log\log(e+c_2s)}{s} \le \frac{\log\log(e+c_2 t)}{t}, \qquad \forall s > t > 0. \]
    We will use this later. Next, from the definition of $\tau_\alpha$, we have
    \begin{align*}
        \tau_\alpha 
        &\le 1 + \max\lrset{ n: n   \le \frac{1}{c_1} \log \log (e+c_2 n) + \frac{f(\alpha)}{c_1} }\\
        &= 1 + \max\lrset{ n: n  - \frac{1}{c_1}\log \log (e+c_2 n) \le \frac{f(\alpha)}{c_1} }\\
        &\le 1 + N_\gamma + \max\lrset{ n : (1-\gamma) n \le \frac{f(\alpha)}{c_1}} \\
        &\le 1 + N_\gamma + \frac{f(\alpha)}{c_1 (1-\gamma)}, \numberthis\label{eq:nbound1}
    \end{align*}
    where
    \begin{align*} 
        N_\gamma 
        :&= \min\lrset{n : (1-\gamma) n \le n  - \frac{1}{c_1}\log\log(e+c_2 n)}\\
        &= \min\lrset{n :  \frac{1}{c_1}\log \log(e+c_2n) \le  n \gamma }.
    \end{align*}
    It is easy to verify that  
    $$t_0 := \frac{4}{c_1\gamma} \log\log\lrp{e+\frac{c_2}{c_1\gamma}} $$
    is feasible for the inequality definining $N_\gamma$, and hence, we get $N_\gamma \le 1 + t_0$. Using this bound on $N_\gamma$ in~\eqref{eq:nbound1}, we get the desired bound on $\tau_\alpha$.
\end{proof}

\end{appendix}
%
%

\end{document}